\input amstex
\documentstyle{amsppt}
\magnification=1200
\hyphenation{comple-xes}
\hoffset=-0.5pc

\define\MZ{{\bold Z}}

\define\Me{{\bold e}}

\define\Ext{\text{\rm Ext}}
\define\Hom{\text{\rm Hom}}
\define\Tor{\text{\rm Tor}}
\define\Ho{\text{\rm H}}

\define\hargl#1{${h_{#1}(\ell)}$}
\define\Mhargl#1{h_{#1}(\ell)}

\define\kargl#1{${k_{#1}(\ell)}$}
\define\Mkargl#1{k_{#1}(\ell)}

\define\qargl#1{${q_{#1}(\ell)}$}
\define\Mqargl#1{q_{#1}(\ell)}

\define\ExtBQV#1{${\Ext^{#1}_B(Q,V)}$}
\define\MExtBQV#1{\Ext^{#1}_B(Q,V)}

\define\TorBQV#1{${\Tor_{#1}^B(Q,V)}$}
\define\MTorBQV#1{\Tor_{#1}^B(Q,V)}

\define\ExtQ#1#2#3{${\Ext_Q^{#1}({#2},{#3})}$}
\define\MExtQ#1#2#3{\Ext_Q^{#1}({#2},{#3})}

\define\TorQ#1#2#3{${\Tor^Q_{#1}({#2},{#3})}$}
\define\MTorQ#1#2#3{\Tor^Q_{#1}({#2},{#3})}

\define\E#1#2#3#4#5{${E^{#1}_{#2}({#3},{#4},{#5})}$}
\define\ME#1#2#3#4#5{E^{#1}_{#2}({#3},{#4},{#5})}

\define\EEzeroq#1{${E^{0,q}_{#1}}$}
\define\MEEzeroq#1{E^{0,q}_{#1}}

\define\ETzeroq#1{${E_{0,q}^{#1}}$}
\define\METzeroq#1{E_{0,q}^{#1}}

\define\EET#1#2{${E^{#1}_{#2}}$}
\define\MEET#1#2{E^{#1}_{#2}}

\define\EExt#1#2#3#4{${\Ext^{#1}_{#2}({#3},{#4})}$}
\define\MEExt#1#2#3#4{\Ext^{#1}_{#2}({#3},{#4})}

\define\TTor#1#2#3#4{${\Tor^{#1}_{#2}({#3},{#4})}$}
\define\MTTor#1#2#3#4{\Tor^{#1}_{#2}({#3},{#4})}

\define\paragraph{\smallskip\par}

\define\MNorm#1#2{1+{#1}+\cdots+{#1}^{{#2}-1}}

\define\sddata#1#2#3#4#5{
$$
 ( #1
\text{
\vbox 
to 1.15 pc
{
  \hbox {$@>{\,\,\, #2\,\,\,}>>$}
  \vskip-.7pc
  \hbox {$@<<{\,\,\, #3\,\,\,}<$}
                }
      }
 #4,#5 )
$$}

\define\Msddata#1#2#3#4#5{
 ( #1
\text{
\vbox 
to 1.15 pc
{
  \hbox {$@>{\,\,\, #2\,\,\,}>>$}
  \vskip-.7pc
  \hbox {$@<<{\,\,\, #3\,\,\,}<$}
                }
      }
 #4,#5 )
}

\nologo
\vsize=57.2truepc
\hsize=38.5truepc
\spaceskip=.5em plus.25em minus.20em
\topmatter
\title{Poisson cohomology and quantization\\
J. f\"ur die Reine und Angew. Math. 408 (1990), 57-113
}\endtitle
\author Johannes Huebschmann{\dag}\endauthor
\affil 
Universit\'e de Lille 1, UFR de Math\'ematiques
\\
UMR CNRS 8524, Laboratoire Paul Painlev\'e
\\
Labex CEMPI (ANR-11-LABX-0007-01)
\\
59655 VILLENEUVE D'ASCQ Cedex/France
\endaffil
\address{\newline
Universit\'e de Lille 1, UFR de Math\'ematiques
\newline
UMR CNRS 8524, Laboratoire Paul Painlev\'e
\newline
59655 VILLENEUVE D'ASCQ Cedex/France
\newline
Johannes,Huebschmann\@math.univ-lille1.fr}
\endaddress
\thanks{{\dag} Supported
 by the Deutsche Forschungsgemeinschaft
 under a Heisenberg grant 1985-1990}
\endthanks
\keywords{Poisson manifolds, Poisson algebras,
Poisson cohomology,
connection, curvature, extensions of Lie algebras,
quantization}
\endkeywords
\subjclass
\nofrills {{\rm 2010}{\it Mathematics Subject Classification}.\usualspace}
17B55, 17B56, 17B63, 17B65, 17B66, 17B83, 53C05, 53D50, 70H99
81Q99, 81S10
\endsubjclass
\abstract{
Let
$R$ be a commutative ring, and let
$(A,\{\cdot,\cdot\})$ 
be a Poisson algebra
over $R$.
We construct a structure of
an $(R,A)$-Lie algebra in the sense of Rinehart
on the $A$-module of K\"ahler differentials
of $A$
depending naturally on $A$ and $\{\cdot,\cdot\}$.
This gives rise to suitable algebraic notions of
Poisson homology and cohomology
for an arbitrary Poisson algebra.
A geometric version thereof
includes the 
\lq canonical homology\rq \ 
and \lq Poisson cohomology\rq\ 
of a Poisson manifold introduced by
Brylinski, Koszul, and Lichnerowicz,
and absorbes the latter
in standard homological algebra
by
expressing them as 
Tor and Ext groups, respectively,
over a
suitable 
algebra of differential operators.
Furthermore, the Poisson structure determines a 
closed 2-form $\pi_{\{\cdot,\cdot\}}$ in the complex
computing Poisson cohomology.
This 2-form generalizes
the 
2-form $\sigma$
defining 
a symplectic structure
on a smooth manifold $N$;
moreover, 
the class 
of $\pi_{\{\cdot,\cdot\}}$
in Poisson cohomology
generalizes the class $[\sigma] \in \Ho^2_{\roman{de Rham}}(N,\bold R)$
of a symplectic structure $\sigma$ on a smooth manifold $N$
and  appears as a crucial ingredient for the construction of
suitable linear representations
of $(A,\{\cdot,\cdot\})$, viewed as a
Lie algebra;
representations of this kind 
occur in quantum theory.
To describe this class and to construct the representations,
 we relate
formal concepts of
connection and curvature
generalizing the classical ones
with extensions of Lie algebras.
We illustrate our results with 
a number of examples of Poisson algebras 
and with a quantization procedure for 
a relativistic particle with zero rest mass
and spin zero.}
\endabstract
\endtopmatter
\document
\leftheadtext{Poisson cohomology and quantization}
\rightheadtext{Johannes Huebschmann}

\bigskip
\noindent
{\smc 0. Introduction}
\medskip
\noindent
The concept of a Poisson manifold is currently
of much interest, see e.~g.
{\smc Berger} [7],
{\smc Bhaskara-Viswanath} [8],
{\smc Braconnier} [10], [11],
{\smc Brylinski} [12],
{\smc Coste-Dazord-Weinstein} [17],
{\smc Conn} [18], [19],
{\smc De Wilde-Le Compte} [21],
{\smc Gelfand-Dorfman}
\linebreak
 [27] -- [29],
{\smc Karasev} [44],
{\smc Kassel} [45],
{\smc Koszul} [52],
{\smc Lichnerowicz} [56] -- [62], 
{\smc Magri-Morosi} [66],
{\smc Magri-Morosi-Ragnisco} [67],
{\smc Mikami-Weinstein} [74],
{\smc Stasheff} [95], [96],
{\smc Tulczyjew} [97], [98],
{\smc Vinogradov-Krasil'shchik} [100],
{\smc Weinstein} [102] -- [106].
A Poisson structure
on a smooth manifold $N$ is a Lie bracket
$\{\cdot,\cdot\}$ on the (multiplicative) algebra
of smooth functions on $N$
satisfying the additional condition
$\{fg,h\} = f\{g,h\} + \{f,h\}g$.
More generally, an algebra $A$ over
a commutative ring $R$
together with a Lie bracket
$\{\cdot,\cdot\}$ on $A$
satisfying the formal analogue
of the above additional condition
is called a {\it Poisson algebra\/}.
The significance of Poisson structures in physics is classical, 
see 
{\smc Lie} [63] and
{\smc Dirac} [22], [23].
\paragraph
 For a symplectic manifold
$(N,\sigma)$, the
rule
$\{f,g\} = \sigma(X_f,X_g)$, where
$X_f$ is the Hamiltonian vector field corresponding to $f$,
defines a Poisson structure on $N$.
However, this is not the only way in which a Poisson structure 
on a manifold
arises,
see e.~g.
{\smc Weinstein} [103].
\paragraph
In [58]
{\smc Lichnerowicz} introduced what he called
\lq\lq Poisson cohomology\rq\rq \ 
of a Poisson manifold;
when the Poisson structure comes from a symplectic
one, this Poisson cohomology
coincides with de Rham cohomology
[58], cf. (3.15) below.
In 
[52]
{\smc Koszul}
introduced a notion of homology
for a Poisson manifold
which was christened
\lq\lq canonical homology\rq\rq \ 
by {\smc Brylinski} [12].
In the present paper we introduce corresponding
notions of Poisson homology and cohomology
for an arbitrary Poisson algebra $(A,\{\cdot,\cdot\})$.
We now explain briefly and informally our approach:
\paragraph
The key idea is that a Poisson structure
$\{\cdot,\cdot\}$
on an arbitrary algebra $A$ 
over a commutative ring $R$
gives rise to a structure
of an $(R,A)$-Lie algebra in the sense of {\smc Rinehart} [80]
on the $A$-module $D_A$ of K\"ahler differentials for $A$
in a natural fashion.
An $(R,A)$-Lie algebra is a Lie algebra over $R$ which 
acts on $A$ and
is also
an $A$-module
satisfying suitable compatibility conditions
which generalize the usual
properties of 
the Lie algebra
of smooth vector fields on a smooth manifold
viewed as a module over its ring of smooth functions;
these objects have been introduced by 
{\smc Herz} [37]
under the name 
\lq\lq pseudo-alg\`ebre de Lie\rq\rq\ 
and were examined by {\smc Palais} [77]
under the name \lq\lq $d$-Lie ring\rq\rq.
Any
$(R,A)$-Lie algebra $L$
gives rise to a complex $\roman{Alt}_A(L,A)$
of alternating forms
which generalizes the usual de Rham complex
of a manifold
and the usual complex
computing 
{\smc Chevalley-Eilenberg} [16]
Lie algebra cohomology.
 This observation is again due to {\smc Palais} [77].
Moreover, 
extending earlier work of
{\smc Hochschild, Kostant and Rosenberg} [39],
{\smc Rinehart\/} [80] has shown that,
when $L$ is projective as an $A$-module,
the homology of
the complex $\roman{Alt}_A(L,A)$
may be identified with
$\roman{Ext}^*_{U(A,L)}(A,A)$
over a suitably defined universal algebra
$U(A,L)$ of differential operators;
see Section 1 below for details.
In particular, when $A$ is 
the algebra of smooth functions
on a smooth manifold $N$ and $L$
the Lie algebra
of smooth vector fields on $N$, then
$U(A,L)$ is the algebra of (globally defined) differential operators
on $N$.
\paragraph
In Section 3 below
we construct, for any 
Poisson algebra $(A,\{\cdot,\cdot\})$,
a natural structure
of an $(R,A)$-Lie algebra
on the $A$-module $D_A$ of K\"ahler differentials for $A$;
we write
$D_{\{\cdot,\cdot\}}$ for the resulting $(R,A)$-Lie algebra.
We can then apply the machinery of {\smc Palais} [77] and
{\smc Rinehart} [80].
In this vein we {\it define \/}
the {\it Poisson cohomology\/}
$\Ho^*_{\roman{Poisson}}(A,\{\cdot,\cdot\};A)$
of $(A,\{\cdot,\cdot\})$
as the homology of
$\roman{Alt}_A(D_{\{\cdot,\cdot\}},A)$.
The Poisson structure $\{\cdot,\cdot\}$ then determines a natural 
closed 2-form
$\pi_{\{\cdot,\cdot\}} \in \roman{Alt}^2_A(D_{\{\cdot,\cdot\}},A)$
and hence a natural
class
$[\pi_{\{\cdot,\cdot\}}] \in \Ho^2_{\roman{Poisson}}(A,\{\cdot,\cdot\};A)$,
see (3.10) for details;
we refer to this class as the
{\it Poisson class\/} of     $(A,\{\cdot,\cdot\})$.
These notions of
Poisson homology and cohomology are entirely algebraic.
In the special case where
the ground ring $R$ is that of the reals
(or complex numbers) and
$A$  is the ring of smooth 
functions 
on a smooth finite dimensional Poisson manifold $(N,\{\cdot,\cdot\})$,
a variant of the above construction yields
a natural structure
of an $(R,A)$-Lie algebra
on the $A$-module $D_A^{\roman{geo}}$ of smooth 1-forms
on $N$;
we write
$D_{\{\cdot,\cdot\}}^{\roman{geo}}$ for the resulting $(R,A)$-Lie algebra.
Again we
can then apply the machinery of {\smc Palais} [77] and
{\smc Rinehart} [80].
We {\it define \/}
{\it geometric Poisson cohomology\/}
$\Ho^*_{\roman{Poisson}}(N,\{\cdot,\cdot\};R)$
of $(N,\{\cdot,\cdot\})$
as the homology of
$\roman{Alt}_A(D_{\{\cdot,\cdot\}}^{\roman{geo}},A)$.
The Poisson structure $\{\cdot,\cdot\}$ then determines a natural 
closed 2-form
$\pi_{\{\cdot,\cdot\}}^{\roman{geo}} \in 
\roman{Alt}^2_A(D_{\{\cdot,\cdot\}}^{\roman{geo}},A)$
and hence a natural
class
$[\pi_{\{\cdot,\cdot\}}^{\roman{geo}}] 
\in \Ho^2_{\roman{Poisson}}(N,\{\cdot,\cdot\};R)$,
see (3.12) for details;
we refer to this class as the
{\it Poisson class\/} of $(N,\{\cdot,\cdot\})$.
The complex 
$\roman{Alt}_A(D_{\{\cdot,\cdot\}}^{\roman{geo}},A)$ 
is precisely the one introduced by
{\smc Lichnerowicz} [58].
In this way
Lichnerowicz' Poisson cohomology appears as
$\roman{Ext}^*_{U(A,D_{\{\cdot,\cdot\}}^{\roman{geo}})}(A,A)$.
Furthermore, it turns out that the obvious morphism
$
D_{\{\cdot,\cdot\}} \longrightarrow
D_{\{\cdot,\cdot\}}^{\roman{geo}}
$
is one of $(R,A)$-Lie algebras
and induces an isomorphism
$
\Ho^*_{\roman{Poisson}}(N,\{\cdot,\cdot\};R)
\to
\Ho^*_{\roman{Poisson}}(A,\{\cdot,\cdot\};A)
$
on cohomology.
Hence in the smooth case there is no need to distinguish between
algebraic and geometric Poisson cohomology;
see (3.12.13) below.
Moreover,
when the Poisson structure $\{\cdot,\cdot\}$
comes from a symplectic structure $\sigma$ on $N$,
this structure induces an isomorphism
$\sigma^*
\colon
\Ho^*_{\roman{de Rham}}(N,\bold R)
\to
\Ho^*_{\roman{Poisson}}(N,\{\cdot,\cdot\};R),
$
and under this isomorphism
the class
$[\sigma] \in \Ho^2_{\roman{de Rham}}(N,\bold R)$
goes to the Poisson class
$[\pi_{\{\cdot,\cdot\}}] \in
\Ho^2_{\roman{Poisson}}(N,\{\cdot,\cdot\};R)$.
Likewise,
inspection shows that
 a suitable complex computing
$\roman{Tor}_*^{U(A,D_{\{\cdot,\cdot\}}^{\roman{geo}})}(A,A)$
is exactly
the one used by
{\smc Koszul} [52]
and
{\smc Brylinski} [12]
to define canonical homology.
It admits an obvious generalization to arbitrary
Poisson algebras;
this leads to our notion of 
(algebraic)
Poisson homology,
see Section 3 for details.
However, in the smooth case the two notions
of Poisson homology differ.
\paragraph
The concept of an $(R,A)$-Lie algebra has a geometric analogue
which is nowadays called a
{\it Lie algebroid\/},
see
{\smc Coste-Dazord-Weinstein} [17],
{\smc Mackenzie} [64],
{\smc Pradines} [78],
{\smc Weinstein} [105].
To our knowledge
the first ones to notice
that a Poisson structure 
on a smooth manifold
gives rise to a Lie bracket
on the space of its 1-forms
were
{\smc Magri and Morosi} [66], 
see also (2.2) in {\smc Magri-Morosi-Ragnisco} [67];
however, it seems that only
in
(3.1) of {\smc Weinstein} [105] 
and
(III.2.1) of {\smc Coste-Dazord-Weinstein} [17] 
is it pointed out that the bracket yields
in fact
a structure of a {\it Lie   algebroid\/}.
The relationship 
of a Poisson structure
with
the work of {\smc Palais} [77]
and {\smc Rinehart} [80]
does not seem to have been observed in the literature
so far.
\paragraph
In Section 1 of the present paper we extend some of the results
in {\smc Rinehart's} paper [80]. 
In Section 2 we introduce  formal concepts
of connection and curvature
which over a smooth manifold boil down to the usual ones.
Among others, a connection appears as a section
(of the underlying modules)
of a suitable extension of Lie algebras,
and its curvature is
a corresponding
(in general non-abelian) 2-cocycle;
for example, the Bianchi identity is then nothing else than
the cocycle condition.
In the geometric situation such a description
goes back to {\smc Atiyah} [5]
and has recently been reworked and elaborated upon by
{\smc Mackenzie} [64].
In Section 3 we introduce Poisson homology and cohomology,
relate it with the earlier notions, and
give some examples.
Section 4 deals with the problem of
constructing a linear representation
of the Lie algebra underlying a Poisson algebra
$(A,\{\cdot,\cdot\})$
in such a way that the elements of the ground ring $R$
act as scalar operators, i.~e. by the usual multiplication;
here the class 
$[\pi_{\{\cdot,\cdot\}}]$
plays a crucial role.
This is motivated by 
geometric quantization theory
{\smc I. Segal} [84], {\smc Kostant} [48], {\smc Souriau} [93],
see also 
{\smc Woodhouse} [108] and the literature there.
The problem of quantizing Poisson algebras
which are not associated with a symplectic manifold
really arises in physics,
see e.~g.
{\smc Gotay} [32], {\smc \'Sniatycki-Weinstein} [92],
and Section 5 below.
In these two papers certain singular systems are treated in analogy
with the symplectic case.
Our approach pushes the analogy further.
Below we shall show that the usual prequantization construction
carries over to arbitrary Poisson algebras,
 with the roles of the 
second integral cohomology group and the
symplectic class being
played by the 
Picard group and the
Poisson class introduced in the present paper, respectively.
Thus our approach
offers tools to handle
singular systems 
with non-trivial Poisson class.
We only mention at this stage that,
for any real Poisson algebra $(A,\{\cdot,\cdot\})$,
the usual notions of 
a polarization 
(of a symplectic structure) and of
quantizability
can be rephrased and generalized in terms
of 
the $(\bold R,A)$-Lie algebra structure on
$D_{\{\cdot,\cdot\}}$
or $D_{\{\cdot,\cdot\}}^{\roman{geo}}$
(as appropriate);
see Sections 4 and 5 for details.
For example, when 
$A$ is the Poisson algebra of smooth functions on 
a
real Poisson manifold $(N,\{\cdot,\cdot\})$ and when $H_1,\dots,H_n$
 are
smooth functions on $N$  that Poisson commute
pairwise, their differentials
$dH_1,\dots,dH_n$
generate a sub $(\bold R,A)$-Lie algebra of
$D_{\{\cdot,\cdot\}}^{\roman{geo}}$ that is isotropic with respect to
$\pi_{\{\cdot,\cdot\}}^{\roman{geo}}$.
A special case is that of a symplectic $2n$-dimensional manifold
$(N,\sigma)$ with $n$ \lq\lq independent
integrals of motion\rq\rq. We believe that our description 
in terms of differentials
is somewhat 
closer to 
the old idea of separation of variables
than the usual description in terms of
Hamiltonian vector fields.  
The obvious advantage of our description is
its applicability to situations
where arguments involving a symplectic structure
are not available.
For illustration, 
 we apply our methods to 
a relativistic particle with zero rest mass
and spin zero
in Minkowski space $Q$ in Section 5.
The resulting Poisson algebra arises from 
what is called {\it Poisson reduction\/}
of the cotangent bundle
$T^*Q$ with respect to a singular constraint
in {\smc \'Sniatycki-Weinstein} [92],
and the underlying meaningful
phase space has a singularity;
the Poisson algebra of this system
is not associated with a symplectic manifold.
\paragraph
In view of its complete generality,
we believe that our approach to Poisson structures
will prove useful
for other
singular systems and for
infinite dimensional geometrical 
theories
of the kind that have recently arisen in physics,
cf e.~g. {\smc Marsden} [70]
or {\smc Chernoff and Marsden} [15].
Also our approach may well apply to Poisson algebras
arising in different ways, cf. {\smc Berger} [7].
Furthermore, if one wants to investigate local questions
more thoroughly,
our approach can be reworked in the language of
sheafs,
cf. {\smc Berger} [7], {\smc Conn} [18], [19],
and {\smc Kamber-Tondeur} [43]. 
It is likely that our methods can then be extended
to yield among others a quantization procedure for
symplectic varieties including singular ones;
this has recently become an interesting topic of research
in view of {\smc Witten's} [107] topological quantum field theory.
We do not pursue this here 
and likewise we leave aside any analytical considerations --
which may be delicate --
since
this would only add unnecessary complications
to the formal aspects we are about to study.
\paragraph
We include a rather long
 bibliography.
We all wish to believe that we are forerunners
but it is also important not to lose contact with the past.
\paragraph
I am much indebted to
A. Weinstein and J. Stasheff for a number of most
valuable comments on a draft of the paper,
to M. Gotay for discussions about the
quantization of 
a relativistic particle with zero rest
mass and spin zero,
and to K. Mackenzie
for 
discussions
about Lie algebroids
and $(R,A)$-Lie algebras.

\bigskip
\noindent
{\bf 1. $(R,A)$-Lie algebras}
\medskip
\noindent
Let $R$ be a commutative ring,
fixed throughout;
the unadorned tensor product symbol
$\otimes$ will
always refer to the tensor product over $R$.
Recall 
that a {\it  Lie algebra\/}
$(L,\lbrack \cdot,\cdot\rbrack)$
over $R$
consists of an $R$-module $L $ and
a pairing
$
\lbrack \cdot,\cdot\rbrack
\colon
L \otimes L
\longrightarrow
L
$
of $R$-modules,
called
a {\it Lie bracket\/}, which satisfies the relations of
{\it antisymmetry\/} and {\it Jacobi identity\/}.
For $x,y \in L$,
as usual, we also write
$(ad(x))(y) =\lbrack x,y \rbrack$.
Given two  Lie algebras $L$ and $L'$,
a {\it morphism\/} $\phi \colon L \longrightarrow L'$
{\it of  Lie algebras\/} over $R$
is the obvious
thing, i.~e. it is a 
morphisms of $R$-modules which is compatible with the
Lie brackets.
\paragraph
Let $A$ be an  algebra over $R$, not necessarily with 1.
Recall that a {\it  derivation\/}  of $A$
(over $R$)
is a   morphism
$
\delta
\colon
A
\longrightarrow
A
$
of  $R$-modules
so that
$
\delta(ab) = (\delta(a))b + a\delta(b).
$
It is well known that the $R$-module 
$\roman{Der}(A)$ of 
derivations of $A$,
viewed as a submodule of
$\Hom_R(A,A)$,
and with bracket $\lbrack \cdot,\cdot \rbrack$
given by
$
\lbrack \alpha,\beta \rbrack (a) =
\alpha (\beta (a)) - 
\beta (\alpha(a)),
$
where $ \alpha, \beta \in L, \ a \in A,$
is a  Lie algebra over $R$.
If $L$ is a  Lie algebra over $R$,
as usual, 
an {\it action\/} of
$L$ on $A$
is
a morphism
$
\omega 
\colon 
L 
\longrightarrow
\roman{Der}(A)
$
of 
Lie algebras over $R$;
henceforth we shall write
$(\omega(\alpha))(a) = \alpha(a), \,\alpha \in L,\, a \in A$.
Given 
two  algebras $A$ and $A'$
over $R$, two
Lie algebras $L$ and  $L'$ 
over $R$, and actions of
$L$ and $L'$ on $A$ and $A'$, respectively,
a {\it morphism \/}
${
(\phi,\psi)
\colon
(A,L)
\longrightarrow
(A',L')
}$
({\it of actions\/})
is the
obvious thing, i.~e. it consists of a morphism
$\phi \colon
A \longrightarrow A'
$
of 
$R$-algebras
and a morphism
$\psi \colon
L \longrightarrow L'
$
of 
 Lie algebras over $R$
 so that,
for every $a \in A,\, \alpha \in L$,
${
\phi(\alpha(a)) = (\psi(\alpha))(\phi(a)).
}$
\paragraph
For 
an algebra $U$
over $R$,
the {\it associated  Lie algebra\/}
over $R$,
written $LU$ or, with an abuse of notation, just $U$,
has the same underlying $R$-module as $U$,
while 
for $u, v \in U$, as usual
the bracket $\lbrack u,v \rbrack$
is
 given by
$
\lbrack u,v\rbrack  =
uv - 
vu.
$
If $L$ is a  Lie algebra over $R$
and if $M$ is an $R$-module,
as usual, 
an {\it action\/} of
$L$ on $M$ is
a morphism
$  
\omega 
\colon 
L 
\longrightarrow
L\roman{End}_R(M)
$
of 
Lie algebras over $R$,
and $M$ is then said to be a
{\it  (left)\/} $L$-{\it module\/};
henceforth we shall write
$(\omega(\alpha))(x) = \alpha(x), \,\alpha \in L,\, x \in M$.
The precise definition of the concept of
a morphism of such structures is clear and is left to the reader.
\paragraph
Let 
now
$A$ be a commutative
$R$-algebra, not necessarily with 1.
Let $L$ be a Lie algebra over $R$,
let 
$
\mu
\colon
A \otimes L \longrightarrow L
$
be a structure of a left $A$-module on $L$ --
as usual we shall write
${\mu(a\otimes \alpha) = a\,\alpha}$ --
and 
let 
$
\omega 
\colon 
L 
\longrightarrow
\roman{Der}(A)
$
be 
an action of
$L$ on $A$.
As in {\smc Rinehart} [80] we shall refer to $L$ as an 
$(R,A)$-{\it Lie algebra\/},
provided
$$
\alignat 2
(a\,\alpha)(b) &= a\,(\alpha(b)),
\quad &&\alpha \in L, \,a,b \in A,
\tag1.1.a
\\
\lbrack 
\alpha,
a\,\beta
\rbrack
&=
a\,
\lbrack 
\alpha,
\beta
\rbrack
+
\alpha(a)\,\beta,\quad &&\alpha,\beta \in L, \,a \in A.
\tag1.1.b
\endalignat
$$
For an $R$-algebra $A$
and an $(R,A)$-Lie algebra $L$, we shall occasionally refer to the
pair $(A,L)$ as a {\it Lie-Rinehart algebra\/}.
Given 
two 
Lie-Rinehart algebras $(A,L)$ and $(A',L')$,
a {\it morphism \/}
${
(\phi,\psi)
\colon
(A,L)
\longrightarrow
(A',L')
}$
{\it of Lie-Rinehart algebras\/}
is 
a morphism of actions so that, furthermore,
$
\psi \colon
L \longrightarrow L'
$
is a morphism
of $A$-modules where
$A$ acts on $L'$ via $\phi$.
It is clear that, with this notion of morphism,
Lie-Rinehart algebras
constitute a category.
A useful example of a morphism of Lie-Rinehart algebras
will be given in (3.8.4) below.
\paragraph
An example of 
an $(R,A)$-Lie algebra
is 
the $R$-module $\roman{Der}(A)$ of derivations of 
a commutative algebra
$A$
with the obvious $A$-module structure;
here the commutativity of $A$ is crucial.
Indeed, a little thought reveals that
for a non-commutative algebra $U$ over $R$
the $R$-module of derivations of $U$ does {\it not\/}
inherit a structure of a $U$-module.
Another example is  classical:
Let $R$ be the real numbers,
let $N$ be a smooth  manifold,
let $A$ be the algebra of smooth functions 
on $N$, and let $L$ be the Lie algebra of smooth vector fields
on $N$; then,
with the obvious structures,
$L$ is an $(R,A)$-Lie algebra.
This works
for an arbitrary smooth Banach manifold,
modelled on a Banach space, 
see e.~g.
{\smc Schwartz} [83], {\smc Lang} [56].
We note that in general for a 
smooth Banach manifold
a derivation need not even locally come from a
vector field,
see e.~g. p. 105 of {\smc Schwartz} [83].
Similar examples arise in the analytic and algebraic setting.
We leave the details to the reader.

\paragraph
Given an $(R,A)$-Lie algebra $L$ and an
$R$-module $M$ having the structures of a
left $A$-module and 
that of 
a left $L$-module
$\omega
\colon
L
\to \roman{End}(M)
$,
we shall refer to $M$ as an $(A,L)$-{\it module\/}, provided
the actions are compatible, 
i.~e. for 
$\alpha
 \in L,\,a \in A,\,m \in M,$
$$
\align
(a\,\alpha)(m) &= a(\alpha(m)),
\tag1.2.a
\\
\alpha(a\,m) &=   a\,\alpha(m) + \alpha(a)\,m.
\tag1.2.b
\endalign
$$
\paragraph
For example,
let
 $A$ be the algebra of smooth functions 
on a smooth finite dimensional
manifold
$N$,
let
$L$ be the Lie algebra of 
smooth vector fields
on $N$,
let $\xi \colon E \to N$
be a smooth vector bundle
on $N$,
 and  
let $M$ be the $A$-module of smooth sections
of $\xi$;
then
an $(A,L)$-module structure on $M$ is precisely a
(linear) connection on $\xi$ with zero curvature.
We shall elaborate on this in the next Section.
\paragraph
Recall that,
for   
a Lie algebra $L$ over $R$ and  an $R$-module $M$
which is also an $L$-module,
the
$R$-multilinear
functions from $L$ into $M$
with the Cartan-Chevalley-Eilenberg differential $d$ given by 
$$
\aligned
(df)(\alpha_1,\dots\alpha_n)
&=
\quad 
(-1)^n
\sum_{i=1}^n (-1)^{(i-1)}
\alpha_i(f (\alpha_1, \dots,\widehat{\alpha_i},\dots, \alpha_n))
\\
&\phantom{=}+\quad
(-1)^n
\sum_{j<k} (-1)^{(j+k)}f(\lbrack \alpha_j,\alpha_k \rbrack,
\alpha_1, \dots,\widehat{\alpha_j},\dots,\widehat{\alpha_k},\dots,\alpha_n)
\endaligned
\tag1.3
$$
constitute a chain
complex
$\roman{Alt}_R(L,M)$
where
as usual \lq $\ \widehat {}\ $\rq \ 
indicates omission of the corresponding term.
If $L$
is
projective as an $R$-module
this chain complex
computes the usual Lie algebra cohomology
$\Ho^*(L,M)$;
we recall that the latter is defined as usual by
 $\Ho^*(L,M)=\roman {Ext}^*_{U(L,R)}(R,M)$
where $U(L,R)$ denotes the corresponding universal
algebra (= the universal enveloping algebra if $L$ is projective
as an $R$-module).
The sign $(-1)^n$ in (1.3) has been introduced 
according to 
the usual Eilenberg-Koszul
convention in differential homological algebra
for
consistency with 
what is said in a subsequent paper [40 II];
in the classical approach such a sign does not occur.
\paragraph
Likewise,
given $L$-modules
$M'$ and $M''$,
the usual formula 
$$
\alpha(x \otimes y) = \alpha(x) \otimes y + x \otimes \alpha(y),
\quad
\alpha \in L, 
\,x \in M',\, y \in M'',
\tag1.4.1
$$
endows
the tensor product $M' \otimes M''$ with a structure of an $L$-module;
 if
 $M$
is another $L$-module, 
a pairing 
$
\mu
\colon
M'
\otimes
M''
\longrightarrow
M
$
of $R$-modules
which is a morphism of $L$-modules
(with respect to (1.4.1)) will be said to be a
{\it a pairing of $L$-modules\/}.
Given such a pairing
$\mu$  of $L$-modules,
the standard shuffle multiplication of alternating maps 
given by
$$
(\alpha \wedge \beta)(x_1,\dots,x_{p+q}) =
\sum_{\sigma}\roman{sign}(\sigma)
\mu(\alpha(x_{\sigma(1)},\dots,x_{\sigma(p)})
\otimes \beta(x_{\sigma(p+1)},\dots,x_{\sigma(p+q)}))
\tag1.4.2
$$
induces a pairing
$$
\roman{Alt}_R(L,M') \otimes
\roman{Alt}_R(L,M'')
\longrightarrow
\roman{Alt}_R(L,M)
\tag1.5
$$
of chain complexes which is associative in the obvious sense;
here $\sigma$
runs through $(p,q)$-shuffles
and $\roman{sign}(\sigma)$ refers to the sign of $\sigma$.
In particular,  for an $R$-algebra $B$,
the chain complex $\roman{Alt}_R(L,B)$
inherits a structure of a differential graded algebra
which is graded commutative if   $B$ is commutative,
and this structure induces a ring structure on cohomology
as usual.
\paragraph
As before, let $A$ be a commutative $R$-algebra
and let
 $L$ be an $(R,A)$-Lie algebra.
It is not hard to see that,
as observed first by
{\smc Palais} [77],
for an $(A,L)$-module
$M$,
the 
differential on
$\roman{Alt}_R(L,M)$
passes to an
$R$-linear
differential on the
graded $A$-submodule
$\roman{Alt}_A(L,M)$
of $A$-multilinear functions;
we note that
the differential will not be
 $A$-linear unless $L$ acts trivially on $A$.
Before we proceed further we mention that
a distinction between graded $A$-algebras
and differential graded $R$-algebras will persist throughout.
We shall carry out most constructions
over $A$; however, in view of the non-trivial action of $L$ 
on $A$,
most resulting differential graded algebras will
be over the ground ring $R$ only.
\paragraph
Given $(A,L)$-modules 
$M'$ and $M''$,
a little thought reveals that
the formula
(1.4.1) endows the tensor product
$M' \otimes _A M''$
with a structure of an
$(A,L)$-module;
we refer to
$M' \otimes _A M''$
with this structure as the
{\it tensor product of\/}
$M'$ and $M''$
{\it in the category of
$(A,L)$-modules\/}.
Given $(A,L)$-modules 
$M$,\ $M'$,\ and $M''$,
a pairing
${
\mu_A
\colon
M'
\otimes_A
M''
\longrightarrow
M
}$
of $A$-modules
 which is compatible with the $L$-structures
will be said to be a 
{\it pairing of $(A,L)$-modules\/}.
Given a pairing 
${
\mu_A
\colon
M'
\otimes_A
M''
\longrightarrow
M
}$
 of $(A,L)$-modules,
let
$${
\mu = \mu_A \roman{pr}
\colon
M'
\otimes_R
M''
\longrightarrow
M'
\otimes_A
M''
\longrightarrow
M
}$$
be the indicated pairing of
$L$-modules
where \lq pr\rq\ refers to the obvious projection map;
inspection shows that under these circumstances
the 
corresponding pairing (1.5)
with respect to the present $\mu$
induces a pairing
$$
\roman{Alt}_A(L,M') \otimes_R
\roman{Alt}_A(L,M'')
\longrightarrow
\roman{Alt}_A(L,M)
\tag1.5'
$$
of chain complexes over $R$.
In particular,
${\roman{Alt}_A(L,A)}$
inherits a structure of a differential graded commutative algebra
over the ground ring $R$
(but {\it not\/} over $A$ unless $L$ acts trivially on $A$).
We now elaborate on a conceptual explanation 
of  these facts which is 
due to 
{\smc Rinehart} [80].
\paragraph
Given 
an $(R,A)$-Lie algebra $L$, its 
{\it universal object\/}
${(U(A,L),\iota_L,\iota_A)}$
is an $R$-algebra $U(A,L)$ together with a morphism
${
\iota_A
\colon 
A
\longrightarrow
U(A,L)
}$
of $R$-algebras
and
a morphism
${
\iota_L
\colon 
L
\longrightarrow
U(A,L)
}$
of Lie algebras over $R$
having the properties
$$
\iota_A(a)\iota_L(\alpha) 
= \iota_L(a\,\alpha),\quad
\iota_L(\alpha)\iota_A(a) - \iota_A(a)\iota_L(\alpha) 
= \iota_A(\alpha(a)),
$$
and
${(U(A,L),\iota_L,\iota_A)}$
is {\it universal\/} among triples
${(B,\phi_L,\phi_A)}$
having these properties.
More precisely:
\proclaim{1.6}
Given 
\paragraph
(i)\phantom{ii} another
 $R$-algebra $B$, viewed at the same time as a Lie algebra over
$R$,
\newline
\indent
(ii)\phantom{i} a morphism
$
\phi_L 
\colon
L
\longrightarrow
B
$ of Lie algebras over $R$, and
\newline
\indent
(iii) a morphism
$
\phi_A 
\colon
A
\longrightarrow
B
$
of $R$-algebras,
\paragraph
\noindent
so that, for ${\alpha \in L, a \in A}$,
$$
\align
\phi_A(a)\phi_L(\alpha) &= \phi_L(a\,\alpha),
\tag1.6.1
\\
\phi_L(\alpha)\phi_A(a) - \phi_A(a)\phi_L(\alpha) &= \phi_A(\alpha(a)),
\tag1.6.2
\endalign
$$
there is a unique morphism
${
\Phi 
\colon
U(A,L)
\longrightarrow
B
}$
of $R$-algebras
so that
$
\Phi\,\iota_A = \phi_B
$ 
and
$
\Phi\,\iota_L = \phi_L$.
\endproclaim
\paragraph
For example,
when
 $A$ is the algebra of smooth functions on a smooth
manifold $N$ and  $L$  the Lie algebra of 
smooth vector fields
on $N$, then $U(A,L)$ is the {\it algebra of 
(globally defined)
differential operators
on\/} $N$.
\paragraph
The universal property is not spelled out in
{\smc Rinehart} [80].
A universal property equivalent to the above one
is given in {\smc Malliavin} [69] where it is attributed to
{\smc Feld'man} [26].
An explicit construction for the 
$R$-algebra ${U(A,L)}$ is given 
in {\smc Rinehart} [80].
For convenience, we now give
a new alternate construction
which employs the {\smc Massey-Peterson\/} [73] algebra.
Let ${(U(R,L),\iota_L,\iota_R)}$ be the usual universal algebra
for $L$ over $R$,
for $\alpha \in L$, write
$\overline \alpha = \iota_L(\alpha)$,
and consider the algebra
$$
 A \odot U(R,L) = (A \otimes_R U(R,L), \mu),
\tag1.7
$$
whose underlying left $A$-module is the one induced from
$U(R,L)$
as indicated, and whose
multiplication on the generators
is defined by
$$
a\,\overline \alpha = a \otimes \overline \alpha,\quad
\overline \alpha\,a = a\,\overline \alpha +  \alpha (a),
\quad \alpha \in L, a \in A.
$$
The universal property of ${(U(R,L),\iota_L,\iota_R)}$
implies that this is well defined, i.~e. that
$\alpha (a)$ depends only on $\overline \alpha$ and $a$.
Furthermore,
the
Jacobi identity implies that
$\mu$ is associative, i.~e. that
${A \odot U(R,L)}$
is indeed an $R$-algebra.
However, there is a more conceptual way to understand this algebra structure:
It is 
clear that the diagonal map
$
\Delta 
\colon
L
\longrightarrow
L \oplus L
$
is a morphism of Lie algebras.
As is well known, the universal property
of
${(U(R,L),\iota_L,\iota_R)}$
yields a multiplicative extension
to a diagonal map
$\Delta$ for $U(R,L)$ 
which endows the latter
with the structure of a
cocommutative 
Hopf algebra.
Furthermore, since on the generators the diagonal map $\Delta$
is given by
$
\Delta(\overline \alpha) = 
\overline \alpha \otimes 1 +
1 \otimes \overline \alpha
$
as a Hopf algebra,
$U(R,L)$  is
primitively generated.
With the obvious structure,
$A \otimes A$ is an
$(U(R,L)
\otimes
U(R,L))
$-module,
and
the diagonal map 
on $U(R,L)$
induces on
$A \otimes A$ 
a structure of a 
left 
$U(R,L)$-module;
the above requirement that $L$ acts on $A$ by derivations
means precisely that
$A$ is an
 algebra {\it over\/} $U(R,L)$, i.~e.
that the structure map 
$
\mu \colon A \otimes A 
\longrightarrow
A
$
is a morphism of 
left $(U(R,L)$-modules.
The 
algebra 
$ A \odot U(R,L)$
is 
 the corresponding
{\smc Massey-Peterson\/}
algebra [73]; its structure map $\mu$ is given by
$$
(a \otimes u)
(b \otimes v)
=
ab \otimes uv +
\sum a\,u_i'(b)\otimes u_i''v,
\quad a,b \in A, u,v \in U(R,L),
$$
where $\Delta(u) = \sum u_i'\otimes u_i'' \in U(R,L) \otimes U(R,L)$.
Notice that the Hopf algebra structure 
of $U(R,L)$
{\it organizes\/} the requisite 
combinatorics needed to prove associativity of 
the algebra $A \odot U(R,L)$.
We 
note
that the algebra $ A \odot U(R,L)$
together with the obvious morphisms
${
\iota'_A
\colon 
A
\longrightarrow
A \odot U(R,L)
}$
of $R$-algebras
and
${
\iota_L'
\colon 
L
\longrightarrow
A \odot U(R,L)
}$
of Lie algebras over $R$
is what is called 
the {\it algebra of differential
operators of the representation\/}
$\omega \colon L \longrightarrow \roman{Der}(A)$ of $L$ in $A$
on p.175 of {\smc Jacobson} [42].
Our construction
differs from the one in [42].
\paragraph
To complete the construction of the universal object,
let
$J$
be the
right ideal in
$A \odot U(R,L)$
generated by the
elements
${ab\otimes \alpha - a \otimes b\,\alpha,\ a,b \in A,\ \alpha \in L}$,
where
the term ${b\,\alpha}$ in ${a \otimes b\,\alpha}$ refers to
the left $A$-module structure on $L$.
A straightforward calculation shows that
for
${a,b,c \in A,\ \alpha,\beta \in L}$,
$$
(c \otimes \beta)
(ab\otimes \alpha - a \otimes b\,\alpha)
=
(cab\otimes \beta \alpha - c \otimes ab \beta \alpha)
+ ca\beta(b)\otimes \alpha - c \otimes a\beta(b)\alpha,
$$
whence $J$ is a two-sided ideal in
$A \odot U(R,L)$.
Let
$$
U(A,L) = (A \odot U(R,L))/J,
\tag1.8.1
$$
and let $\iota_A$ and $\iota_L$ be the obvious morphisms.
By construction
it is then clear that
${(U(A,L),\iota_L,\iota_A)}$ has the universal property (1.6).
\paragraph
We mention in passing that when $A= R$ with trivial $L$-action
the object ${(U(R,L),\iota_L,\iota_A)}$ is the 
{\it usual universal algebra\/} of $L$ (over $R$).
\paragraph
It is obvious that there is a one-one
correspondence between
${(A,L)}$-modules and 
\linebreak
$U(A,L)$-modules;
this correspondence is in fact
an equivalence of categories.
In particular, 
the obvious $(A,L)$-module structure on $A$ mentioned above
induces
on $A$ that of a left $U(A,L)$-module;
the corresponding structure map is given by
$$
\mu
\colon
U(A,L)
\otimes
A
\longrightarrow
A,
\quad
\mu(\alpha \otimes a) = \alpha(a),
\tag1.8.2
$$
where ${\alpha \in L,\, a \in A}$.
In particular, let
$
\varepsilon
\colon
U(A,L) \longrightarrow A
$
be the obvious morphism given by 
$$
\varepsilon (a) = a,\quad
\varepsilon (a \alpha) =0,\quad
\varepsilon (\alpha a) = \alpha (a).
\tag1.8.3
$$
It is manifestly a morphism of
left $U(A,L)$-modules, but {\it not\/}
one of algebras
unless $L$ acts trivially on $A$, and
its kernel
is the left ideal in $U(A,L)$ generated by $L$.
In particular, the composite 
${\varepsilon \iota_A}$
is the identity map of $A$ whence
$\iota_A$ is injective.
Henceforth we shall identify $A$ with its image in $U(A,L)$, and we shall
not distinguish in notation between the elements of $A$ and their images in
$U(A,L)$.
Furthermore, it is clear that
given 
two 
Lie-Rinehart algebras $(A,L)$ and $(A',L')$,
a  morphism
$
(\phi,\psi)
\colon
(A,L)
\longrightarrow
(A',L')
$
of Lie-Rinehart algebras
induces a
 morphism
${
U(\phi,\psi)
\colon
U(A,L)
\longrightarrow
U(A',L')
}$
of $R$-algebras.
Hence $U(\cdot,\cdot)$ is a functor from 
the  category of Lie-Rinehart algebras
 into the
category of $R$-algebras.
We mention in passing that,
in view of the universal property of $U(A,L)$,
 $A$ also inherits a structure
of a {\it right\/}
$U(A,L)$-module
given by
$$
a \cdot \alpha = -\alpha(a),\quad a \in A, \alpha \in L.
\tag1.8.4
$$

\smallskip
\noindent
{\smc REMARK after publication.}
This construction does not work;
in general there is no such right module
$U(A,L)$-module structure.
This does not cause any difficulty, though.
The general 
right module
$U(A,L)$-module structure
is {\it never\/} used in the paper.
The special case where the right 
$U(A,L)$-module structure
arises from a Poisson structure
as explained between (3.8.6) and (3.8.7)
below 
is {\it correct\/}, that is,
the construction given there exhibits a
right 
$U(A,L)$-module structure.

\paragraph
The universal algebra $U(A,L)$ admits an
 obvious  a filtered algebra structure,
cf. {\smc Rinehart} [80],
with
$U_{-1}(A,L) = 0$
and
$U_p(A,L)$ the left $A$-submodule of
$U(A,L)$ generated by products of at most $p$ elements of
the image  ${\overline L}$ of $L$ in $U(A,L)$;
further, for $a \in A$ and $z \in U_p(A,L)$
we have
$a\,z - z \, a \in U_{p-1}(A,L)$
whence the 
inherited
left and right $A$-module structures on 
the associated graded object $E^0(U(A,L))$ coincide,
and $E^0(U(A,L))$
inherits a structure of a commutative graded $A$-algebra.
The Poincar\'e-Birkhoff-Witt
Theorem for $U(A,L)$
then assumes the following form,
cf. 
 (3.1) of {\smc Rinehart} [80],
where $S_A[L]$ denotes the symmetric $A$-algebra on $L$.
\proclaim{Theorem 1.9 [Rinehart]}
For an $(R,A)$-Lie algebra $L$ 
which 
is projective
as an $A$-module,
 the canonical
$A$-epimorphism
$
S_A[L]
\longrightarrow
E^0(U(A,L))
$
is an isomorphism of $A$-algebras.
\endproclaim
\proclaim{Corollary 1.10}
For an $(R,A)$-Lie algebra $L$ 
which 
is projective
as an $A$-module,
 the morphism
$\iota_L
\colon
L
\longrightarrow
U(A,L)
$ is injective.
\endproclaim

The usual construction of the Koszul complex
computing Lie algebra cohomology
(see e.~g. {\smc Chevalley-Eilenberg} [16])
 carries over
as well:
Let $\Lambda_A(sL)$ be the exterior Hopf algebra over $A$
on the suspension $sL$ of $L$,
where \lq\lq suspension\rq\rq\ 
means that $sL$ is $L$ except that its elements are
regraded by 1.
We shall write typical elements in the form
$$
\langle \alpha_1,\alpha_2,\dots,\alpha_n \rangle
=
(s\alpha_1)(s\alpha_2)\dots s(\alpha_n ) \in \Lambda_A(sL) ,
\quad \alpha_1,\alpha_2,\dots,\alpha_n \in L.
\tag1.11
$$
For $u \in U(A,L)$ and $\alpha_i \in L$ let
$$
\aligned
d(u\otimes \langle \alpha_1, \dots, \alpha_n \rangle)
&=
\quad\sum_{i=1}^n (-1)^{(i-1)}u\alpha_i\otimes 
\langle\alpha_1, \dots,\widehat{\alpha_i},\dots, \alpha_n \rangle
\\
&\phantom{=}+\quad
\sum_{j<k} (-1)^{(j+k)}u\otimes 
\langle \lbrack \alpha_j,\alpha_k \rbrack,
\alpha_1, \dots,\widehat{\alpha_j},\dots,
\widehat{\alpha_k},\dots \alpha_n\rangle,
\endaligned
\tag1.12
$$
where $A$ acts on the right of $U(A,L)$ by means of the embedding
$\iota_A \colon A \longrightarrow U(A,L)$.
This yields
an operator
${
d
\colon
U(A,L) \otimes _A \Lambda_A(sL) 
\longrightarrow
U(A,L) \otimes _A \Lambda_A(sL).
}$
We shall refer to
$$
K(A,L) = (U(A,L) \otimes _A \Lambda_A(sL),d),
\tag1.13
$$
 as the
{\it Koszul complex\/} for $(A,L)$.
It is proved in {\smc Rinehart} [80] that
$d$ is an ${U(A,L)}$-linear differential whence
$K(A,L)$ is indeed a 
chain complex.
It is manifest that the Koszul complex is functorial
in $(A,L)$. Given a morphism
$
(\phi,\psi)
\colon
(A,L)
\longrightarrow
(A',L')
$
of Lie-Rinehart algebras,
we shall denote
the induced morphism by
${
K(\phi,\psi)
\colon
K(A,L)
\longrightarrow
K(A',L').
}$
It is a morphism of $U(A,L)$-modules
where ${K(A',L')}$ is viewed as an
$U(A,L)$-module
via 
${U(\phi,\psi)}$.
Furthermore,
when $L$ is projective or free as a left $A$-module,
$K(A,L)$ is a
projective 
or free
resolution of
$A$ in the category of left $U(A,L)$-modules
according as $L$ is a projective or free
left $A$-module;
details may be found in {\smc Rinehart} [80].
\paragraph
For  an $(R,A)$-Lie algebra $L$
and  an $(A,L)$-module $M$,
as in {\smc Rinehart} [80], we shall write
${
\Ho_A^*(L,M) = \MEExt*{U(A,L)}AM,
}$
and we shall refer to this as the $A$-{\it Lie algebra cohomology
of $L$ with coefficients in\/} $M$.
In the guise of standard homological algebra we view
$\MEExt*{U(A,L)}AM$
as the primary object; it is always defined,
whether or not $L$ is projective as an $A$-module.
It is clear that the above establishes the following.

\proclaim{Proposition 1.14}
Let $L$ be an $(R,A)$-Lie algebra, and let
$M$ be an $(A,L)$-module. 
Then the $A$-multilinear functions
$\roman{Alt}_A(L,M)$ form a subcomplex of 
$\roman{Alt}_R(L,M)$,
and, if
the $L$ underlying $A$-module is projective,
$\roman{Alt}_A(L,M)$
computes 
$${\Ho_A^*(L,M)(=\roman {Ext}^*_{U(L,A)}(A,M))}.
$$
Furthermore, for $M = A$, the usual shuffle product
{\rm (1.4)} induces 
a structure of a differential graded commutative algebra on
$\roman{Alt}_A(L,A)$ which induces that of a graded commutative algebra
on ${\Ho_A^*(L,A)}$.
\endproclaim
For example, when $A$ is the algebra of smooth functions on a 
smooth 
manifold $N$ and $L$ the Lie algebra of smooth vector fields
on $N$, then $\Hom_{U(A,L)}(K(L,A),A)$ 
is the usual {\it de Rham complex of\/} $N$
and the de Rham cohomology of $N$ is
$\MEExt*{U(A,L)}AA$ over
  the algebra $U(A,L)$ of differential operators on $N$.
Likewise,
for a Lie algebra $L$ over $R$
acting trivially on $R$
 and an $L$-module $M$,
the object
$
K(R,L)
$
is the usual Koszul complex;  
in particular, when $L$ is projective
as an $R$-module,
${K(R,L)}$ is the usual Koszul resolution
computing Lie algebra homology and cohomology.
\paragraph
A partial converse to (1.14) is given by
(1.15) below, whose proof
is routine and left to the reader:

\proclaim{Theorem 1.15}
Let $A$ be a commutative $R$-algebra, let $L$ be a left $A$-module,
let
\linebreak
${
\omega
\colon
L 
\longrightarrow
\roman{End}_R(A)
}$
be a morphism
of $R$-modules,
and let
${
\lbrack \cdot,\cdot\rbrack
\colon
L \otimes_R L 
\longrightarrow
L
}$
be a skew symmetric pairing as indicated.
For 
$\alpha \in L$ and $a \in A$,
write
$ \alpha(a) = (\omega(\alpha))(a)$,
define an operator
$d$ on
the graded commutative algebra
$\roman{Alt}_R(L,A)$
of $R$-multilinear alternating functions 
by means of 
{\rm (1.3)},
and suppose that
$d$
passes to an operator on the subalgebra
$\roman{Alt}_A(L,A) \subseteq \roman{Alt}_R(L,A)$
which  endows
$\roman{Alt}_A(L,A)$
with a structure of a
differential graded commutative algebra.
Then 
$\omega$
factors through $\roman{Der}(A)$
so that
$$
\align
(a\,\alpha)(b) &= a\,(\alpha(b)),
\quad 
\alpha \in L, \ a,b \in A,
\tag1.15.1
\\
\alpha_1(\alpha_2(a)) -\alpha_2 ((\alpha_1)(a)) 
&= \lbrack \alpha_1,\alpha_2 \rbrack(a),\, 
\alpha_1,\alpha_2 \in L, a \in A.
\tag1.15.2
\endalign
$$
Furthermore, when $L$ is projective as an $A$-module,
$\omega$ and $\lbrack \cdot,\cdot\rbrack$
yield an $(R,A)$-Lie algebra structure on $L$.
\endproclaim
The next aim is
to introduce
what will be called
\lq\lq induced structures\rq\rq;
their significance for the
present paper will emerge in (3.18) below.
I am indebted to K. Mackenzie
for suggesting the following description
of induced structures
which replaces a more clumsy one
given in an earlier version of the paper.
\paragraph
Suppose the following data are given:
\paragraph
-- $R$-algebras $A$, $A'$,
\newline\indent  
-- an $(R,A)$-Lie algebra $L$, with structure maps
$\omega 
\colon
L
\longrightarrow
\roman{Der}(A)
$
and
\newline\indent \phantom{p} 
$[\cdot,\cdot] \colon L \otimes _R L \to L$,
\newline\indent  
-- an action 
$\tilde\omega 
\colon
L
\longrightarrow
\roman{Der}(A')
$
of $L$ on $A'$
(but $L$ is not assumed to admit an
\newline\indent\phantom{p} 
$A'$-module 
structure),
\newline\indent  
-- a morphism
${\phi 
\colon
A
\longrightarrow
A'
}$
of algebras
which is also a morphism of $L$-modules.
\paragraph
\noindent
Write $L' = A' \otimes _A L$;
it is an $A'$-module in an obvious fashion.
Our aim is
to endow
$L'$ with a structure of an $(R,A')$-Lie algebra.
To this end,
we consider the obvious pairings
$$
A' \otimes_R L \otimes _R A' \otimes_R L
\longrightarrow
A' \otimes _A L,
\tag1.16.1
$$
given by $u \otimes \alpha \otimes v \otimes \beta
\longmapsto
uv \otimes \lbrack \alpha,\beta \rbrack
-(v\beta(u))\otimes \alpha
+ (u \alpha(v)) \otimes \beta,$
where $u,v \in A',\,\alpha,\beta \in L,$ and
$$
A' \otimes_R L \otimes_R A' 
\longrightarrow A',
\quad u \otimes \alpha \otimes v \longmapsto u\cdot \alpha(v),
\quad u,v \in A',\alpha \in L.
\tag1.16.2
$$
\proclaim{Proposition 1.16}
Under the above circumstances, suppose that, for every
$a \in A$ and for every $\alpha \in L$,
$$
\tilde\omega (a \alpha) = \phi(a) \tilde\omega(\alpha) .
\tag 1.16.0
$$
Then
{\rm (1.16.1)} induces on $L'$ a structure
$[\cdot,\cdot]' \colon L' \otimes _R L' \to L'$
of a Lie algebra (over $R$),
{\rm (1.16.2)} induces an action
$\omega' \colon L' \to \roman{Der}(A')$,
and 
$[\cdot,\cdot]'$ and
$\omega'$ endow
$L'$ with a structure of an 
${(R,A')}$-Lie algebra
in such a way that
$$
(\phi,\phi \otimes \roman{Id})
\colon (A,L) \longrightarrow (A',L')
\tag1.16.3
$$
is a morphism of Lie-Rinehart algebras.
\endproclaim
\demo{Proof}
This comes down to routine checking and is therefore left to the reader. \qed
\enddemo
In the situation of this Proposition
we shall say that $(L',[\cdot,\cdot]',\omega')$ is
 {\it induced from\/}
${\phi}$;
often we shall then write $L'$ rather than  $(L',[\cdot,\cdot]',\omega')$.
Moreover, 
we then have the two morphisms
$$
(\roman{Id},\phi_*)
\colon
\roman{Alt}_A(L,A) 
\longrightarrow
\roman{Alt}_A(L,A') 
$$
and
$$
((\phi \otimes \roman{Id})^*,\roman{Id})
\colon
\roman{Alt}_{A'}(L',A') 
\longrightarrow
\roman{Alt}_A(L,A')
$$
of differential graded algebras.
However, the latter is a standard
adjointness isomorphism
and hence
$\phi$ induces a morphism
$$
\roman{Alt}_A(L,A)
\longrightarrow
\roman{Alt}_{A'}(A' \otimes _A L,A')  
$$
of differential graded algebras.
\paragraph
\noindent
{\smc Example 1.16.4.}
Let 
$A$ be an algebra over $R$, let
$g$ be a Lie algebra over $R$,
and let $\omega \colon g \to \roman{Der}(A)$ be an action
of $g$ on $A$. 
Then with the obvious change in notation,
(1.16.1) and (1.16.2) endow
$A \otimes g$ with a structure of an $(R,A)$-Lie algebra.
It is referred to as a {\it crossed product\/}
(produit crois\'e) in
{\smc Malliavin} [69].
We shall have to say more about this 
$(R,A)$-Lie algebra
in (3.18) below.
\paragraph
A geometric analogue of the above notion of induced structure
may be found in {\smc Higgins-Mackenzie} [38].

\beginsection 2. Extensions, connection, curvature

In this Section we relate extensions of $(R,A)$-Lie algebras
with formal concepts of connection and curvature.
This extends 
earlier work of {\smc Atiyah} [5] and
others.
Historical comments will be given at the end of this Section.
\paragraph
Let $A$ be an algebra,
and let $L'$, $L$, $L''$ be $(R,A)$-Lie algebras.
Extending common notation, we refer to
a short exact sequence
$$
\Me
\colon
0
\to
L'
\to
L
\to
L''
\to
0
\tag2.1
$$
in the category of $(R,A)$-Lie algebras as an
{\it extension\/} of
$(R,A)$-Lie algebras;
notice in particular that the
Lie algebra $L'$ necessarily acts trivially on $A$.
If also
$
\bar
\Me
\colon
0
\to
L'
\to
\bar L
\to
L''
\to
0
$
is an extension of
$(R,A)$-Lie algebras,
as usual, $\Me$ and $\bar \Me$ are
said to be {\it congruent\/},
if there is a
morphism
$
(\roman {Id},\cdot,\roman {Id})
\colon \Me \longrightarrow\bar \Me
$
of extensions of
$(R,A)$-Lie algebras.
\paragraph
An extension of the kind (2.1) may be represented
by a 2-cocycle, provided
the extension 
(2.1)
splits in the category of $A$-modules, e.~g. if
$L''$ is projective as an $A$-module.
More precisely, let
$
\omega 
\colon
L''
\to
L
$
be a section of $A$-modules
for the projection $L \to L''$.
We shall occasionally refer
to $\omega$ as an $\Me$-{\it connection\/}.
Given an
$\Me$-connection, define
the corresponding
($\Me$-){\it curvature\/}
$
\Omega
\colon
L'' \otimes_A L''
\to
L'
$
as the morphism $\Omega$ of $A$-modules
satisfying
$$
[\omega(\alpha),\omega(\beta)]
=
\omega[\alpha,\beta]
+ \Omega(\alpha,\beta)
\tag2.2
$$
for every $\alpha, \beta  \in L''$;
a little thought reveals that $\Omega$ is indeed well defined
as an alternating $A$-bilinear 2-form on $L''$ with values in $L'$.
Since in $L$ the Jacobi identity holds,
the morphism $\Omega$ must satisfy the following
2-cocycle condition:
The Lie algebra $L$ acts on $L''$ via the adjoint representation
$\roman{ad} \colon L \to \roman{End}(L'')$.
Define an operator
$D^{\omega}$
on the complex
$\roman{Alt}_A(L',L'')$
by means of the formal analogue of (1.3), i.~e. by
$$
\aligned
(D^{\omega}f)(\alpha_1,\dots\alpha_n)
&=
\quad 
(-1)^n
\sum_{i=1}^n (-1)^{(i-1)}
\roman {ad}(\omega(\alpha_i))
(f (\alpha_1, \dots,\widehat{\alpha_i},\dots, \alpha_n))
\\
&\phantom{=}+\quad
(-1)^n
\sum_{j<k} (-1)^{(j+k)}f(\lbrack \alpha_j,\alpha_k \rbrack,
\alpha_1, \dots,\widehat{\alpha_j},\dots,\widehat{\alpha_k},\dots,\alpha_n).
\endaligned
\tag2.3
$$
Then the Jacobi identity in $L$ boils down to
$$
D^{\omega}(\Omega) = 0.
\tag2.4
$$
We refer to the latter as the {\it generalized Bianchi identity\/}.
\paragraph
We mention in passing that these concepts
of $\Me$-connection and
$\Me$-curvature generalize
the notions of
principal
connection and curvature;
details will be explained
in a follow up paper [41 I].
\paragraph
As one would expect,
the above 2-cocycle is unique up to a coboundary.
More precisely, let
$
\omega' 
\colon
L''
\to
L
$
be another section 
of $A$-modules, and 
let
$
\Omega'
\colon
L'' \otimes L''
\to
L'
$
be the corresponding morphism of $A$-modules so that
$$
[\omega'(\alpha),\omega'(\beta)]
=
\omega'[\alpha,\beta]
+ \Omega'(\alpha,\beta).
$$
Then $\omega- \omega'$ 
factors through
an $A$-linear morphism
$u \colon L'' \to L'$ and, if we view $u$ as an element of
$\roman{Alt}_A(L'',L')$,
we obtain $\Omega'-\Omega$ as the 1-coboundary of $u$
in the sense that, for $\alpha,\beta \in L$,
$$
(\Omega'-\Omega) (\alpha,\beta)
=
\roman{ad}(\omega(\alpha))(u(\beta)) + 
\roman{ad}(\omega'(\beta))(u(\alpha))
- u([\alpha,\beta]).
\tag2.5
$$
\paragraph
It is well known that when $L'$ is abelian
which under the present circumstances means that
$L'$ is just an $A$-module
with trivial Lie bracket,
the
adjoint action of $L$ on $L'$ induces an action 
$\rho \colon L'' \to \roman{End}(L')$
of $L''$ on
$L'$.
Inspection shows that this then endows $L'$ in fact with the structure of
an $(A,L'')$-module.
Furthermore, the cocycle condition
then boils down to
$d\Omega = 0 \in 
 \roman{Alt}_A(L'',L')$ 
with respect to the action $\rho$ -- which is now independent of
a section 
of $A$-modules of the kind $L'' \to L$,
and the classical
argument 
in Eilenberg-Mac Lane cohomology,
cf. e.~g. VII.3 in {\smc Hilton-Stammbach} [39],
may easily be extended to a proof of the following.
\proclaim{Theorem 2.6}
Let $L'$ and $L''$ be $(R,A)$-Lie algebras,
assume that $L'$ is abelian,
and let
$\rho \colon L'' \to \roman{End}(L')$
be a structure of an
$(A,L'')$-module
 on
$L'$.
 Then
the assignment of a 2-cocycle $\Omega \in \roman{Alt}_A(L'',L')$ 
to its extension {\rm (2.1)} of $(A,L'')$-Lie algebras
yields  a bijection between the congruence classes 
of extensions of $L'$ by $L''$ 
whose underlying extension of $A$-modules split
and the
classes in  
$\Ho^2(\roman{Alt}_A(L'',L')) \,\,( =\Ho_A^2(L'',L')$ 
if $L''$ is $A$-projective).
\endproclaim
\paragraph
In the case where $L'$ is non-abelian,
$\Omega$ is a non-abelian 2-cocycle 
in a suitable sense and hence it does {\it not\/} lead to a 
cohomology class in 
the above sense.
We shall explain elsewhere
a formal analogue of the classical Chern-Weil construction [41 I];
(2.6) will then be a special case thereof.
\paragraph
We now relate the above material 
to the classical notions of
connection and curvature.
We pursue this here only
as far as needed for the study of Poisson
algebras in the next Section.
We shall elaborate further on these ideas in [41].
\paragraph
Let $L$ be an $(R,A)$-Lie algebra and let
$M$ be a
left $A$-module.
Then an $R$-linear morphism
$\omega
\colon
L
\to \roman{End}(M)
$
will be referred to as a (linear)
$L$-{\it connection\/} on 
 $M$, 
if for 
$\alpha
 \in L,\,a \in A,\,m \in M,$
$$
\align
(a\,\alpha)(m) &= a(\alpha(m)),
\tag2.7.a
\\
\alpha(a\,m) &=   a\,\alpha(m) + \alpha(a)\,m.
\tag2.7.b
\endalign
$$
Here is another way to say this:
Write 
$
D^{\omega} \colon M
\to
\Hom(L,M)
$ 
for the indicated  adjoint of $\omega$,
so that
$$
D^{\omega}_{\alpha}(m) = (\omega(\alpha))(m),
\quad \alpha \in L,\ m \in M;
\tag2.7.c
$$
then
instead of (2.7.a) we may require that
$D^{\omega}$ factors through
$\Hom_A(L,M)$ , i.~e. that it may be displayed as
$$
D^{\omega} \colon M
\to
\Hom_A(L,M),
\tag2.7.a'
$$
and the rule (2.7.b) then assumes the formally well known form 
$$
D^{\omega}_{\alpha}(am) = (\alpha (a))m + a D^{\omega}_{\alpha}(m).
\tag 2.7.b'
$$
Given a linear
$L$-connection $\omega\colon L \to \roman{End}(M)$ on $M$,
as usual,
extend 
$D^{\omega}$
to the operator of {\it covariant derivative \/} on 
$\roman{Alt}_A(L,M)$
again by means of  the formal analogue of (1.3), i.~e. by
$$
\aligned
(D^{\omega}f)(\alpha_1,\dots\alpha_n)
&=
\quad 
(-1)^n
\sum_{i=1}^n (-1)^{(i-1)}
(\omega(\alpha_i))
(f (\alpha_1, \dots,\widehat{\alpha_i},\dots, \alpha_n))
\\
&\phantom{=}+\quad
(-1)^n
\sum_{j<k} (-1)^{(j+k)}f(\lbrack \alpha_j,\alpha_k \rbrack,
\alpha_1, \dots,\widehat{\alpha_j},\dots,\widehat{\alpha_k},\dots,\alpha_n),
\endaligned
\tag2.8
$$
and define the {\it curvature\/} 
$
\Omega
\colon
L \otimes_A L \to \roman{End}(M)
$
of $\omega$
as the adjoint of
$$
D^{\omega} D^{\omega}
\colon
M
\longrightarrow
\roman{Alt}^2_A(L,M).
\tag2.9
$$
Explicitly, with $D = D^{\omega}$,
this is the formally well known formula
$$
\Omega(\alpha,\beta) = 
D_\alpha D_\beta - D_\beta D_\alpha - D_{[\alpha,\beta]}.
$$
The standard argument shows that $\Omega$ is a \lq\lq tensor\rq\rq, i.~e.
that it is a 2-form with values in
$\roman{End}_A(M)$.
Furthermore,
we now literally have the {\it Bianchi identity\/}
$$
D^{\omega}(\Omega) = 0 
\in  \roman{Alt}_A(L,\roman{End}_A(M)).
\tag2.10
$$
An $L$-connection will be said to be {\it flat\/},
if its curvature is zero.
We note that when $A$ is the algebra over the
reals $\bold R$ 
of smooth functions on a smooth
finite dimensional manifold $N$,
when $L$ is the
$(\bold R,A)$-Lie algebra of smooth vector fields
on $N$, and when
$M$ is the $A$-module
of smooth sections of a vector bundle
over $N$,
the above notions boil down to the usual ones.
\paragraph
Let $M$ be an $A$-module, and let $L$ be an $(R,A)$-Lie algebra.
We now introduce 
an $(R,A)$-Lie algebra $\roman{DO}(A,L,M)$ 
that acts on $M$ by 
the analogue of infinitesimal gauge transformations.
I am indebted to A. Weinstein for asking whether
there is such an object
in general since it is well known to exist 
in the special case where $M$ is the module of sections of a vector bundle.
The elements of
$\roman{DO}(A,L,M)$ 
may be viewed as acting as
\lq\lq differential operators\rq\rq\ on $M$, whence
the notation \lq DO\rq. 
A related (geometric) object
for a smooth vector bundle $E$,
denoted $\roman{CDO}(E)$,
is introduced on p. 103 of {\smc Mackenzie} [64].
\paragraph
Consider the direct sum $\roman {End}_R(M) \oplus L$,
equipped with the obvious componentwise Lie algebra structure so that,
for every $\beta,\beta' \in \roman {End}_R(M)$ and
every $\alpha,\alpha' \in L$,
$$
[(\beta,\alpha),(\beta',\alpha')] 
= (\beta \beta' - \beta' \beta,[\alpha,\alpha'])
\in \roman {End}_R(M) \oplus L,
\tag2.11.1
$$
and let
$\roman{DO}(A,L,M) \subseteq \roman {End}_R(M) \oplus L$
be the $R$-submodule consisting of those pairs 
\linebreak
$(\beta,\alpha)\in \roman {End}_R(M) \oplus L$
that satisfy
$$
\beta(am) = (\alpha(a))m + a (\beta(m)),\quad a \in A,\, m \in M.
\tag2.11.2
$$
Further,
for $a \in A,\, m \in M$,
and $\beta \in \roman {End}_R(M)$,
define $a \beta \in \roman {End}_R(M)$
by
$$
(a\beta)(m) = a (\beta(m)).
\tag2.11.3
$$
Then a little thought reveals that,
in view of (1.2.a) and the commutativity of $A$,
 for every
$(\beta,\alpha)\in \roman{DO}(A,L,M)$ and
$a,b \in A,\, m \in M$,
$$
(b\beta)(am) = ((b \alpha)(a))(m) + a (b \beta)(m),
\tag2.11.4
$$
whence the rule
$$
(b, \beta,\alpha) \mapsto (b\beta,b \alpha),\quad 
(\beta,\alpha) \in \roman{DO}(A,L,M),\, b \in A,
\tag2.11.5
$$
endows
$\roman{DO}(A,L,M)$
with a structure of a (left) $A$-module.
Moreover,
an easy calculation shows that,
for every $(\beta,\alpha),(\beta',\alpha') \in \roman{DO}(A,L,M)$,
and
$a \in A,\, m \in M$,
$$
(\beta \beta' - \beta' \beta)(am)
=
([\alpha,\alpha'](a))m
+
a
(\beta \beta' - \beta' \beta)(m) \in M,
\tag2.11.6
$$
whence
$\roman{DO}(A,L,M)$
inherits a structure of a Lie algebra over $R$
from
$\roman {End}_R(M) \oplus L$.
Finally, the obvious morphism
$$
\roman{DO}(A,L,M) \longrightarrow L
\tag2.11.7
$$
is manifestly a morphism of $A$-modules and of
Lie algebras over $R$,
and this morphism,
combined with the given action of $L$ on $A$,
induces an action of
$\roman{DO}(A,L,M)$ on $A$
in such a way that
$\roman{DO}(A,L,M)$
is an $(R,A)$-Lie algebra.
It is clear that the obvious inclusion 
$\roman{End}_A(M) \subseteq\roman{End}_R(M) $
induces 
an injection
$\roman{End}_A(M) \subseteq \roman{DO}(A,L,M) $
of $(R,A)$-Lie algebras.
Now, with the zero morphism
$\roman{End}_A(M) \to \roman{Der}(A)$
and the obvious Lie algebra structure,
the object $\roman{End}_A(M)$ is an $(R,A)$-Lie algebra,
and
by construction, we have an exact sequence
$$
0
\longrightarrow
\roman{End}_A(M) 
\longrightarrow 
\roman{DO}(A,L,M)
\longrightarrow
L
\tag2.11.8
$$
of $(R,A)$-Lie algebras.
Moreover, 
it is clear that by construction
the 
obvious morphism
$$
\roman{DO}(A,L,M) \longrightarrow \roman{End}_R(M)
\tag2.11.9
$$
induces a structure of an
$(A,\roman{DO}(A,L,M))$-module.
We refer to
the $(R,A)$-Lie algebra
$\roman{DO}(A,L,M)$ as the
{\it infinitesimal gauge algebra\/}
of $M$ {\it with respect to\/} $L$,
and to the endomorphisms 
$\beta$
of $M$ 
coming from
$\roman{DO}(A,L,M)$
as
{\it infinitesimal gauge transformations\/}
of $M$ {\it with respect to\/} $L$.
The reason for this terminology
will be given in (2.16)(2) below.
We note that
$\roman{DO}(A,L,M)$ is an invariant of
$M$ and $L$.
Moreover, 
borrowing some terminology from Galois theory,
we shall say that
$M$ is $L$-{\it normal\/},
if the morphism (2.11.7) is surjective,
i.~e. if for every $\alpha \in L$ there is an $R$-linear
endomorphism $\beta$ of $M$ so that
(2.11.2) holds.
Thus,  an $L$-normal $A$-module $M$
has an exact sequence
$$
\Me_M \colon
0
\longrightarrow
\roman{End}_A(M) 
\longrightarrow 
\roman{DO}(A,L,M)
\longrightarrow
L
\longrightarrow
0
\tag2.11.10
$$
of $(R,A)$-Lie algebras
that is natural in terms of the given data.
It is clear that an $A$-module
$M$ admits an $L$-connection if and only if
it is $L$-normal and if the
corresponding extension
(2.11.10)
splits in the category of $A$-modules.
Moreover, if this happens to be the case,
after a connection 
$\omega$ for $M$
has been chosen,
the corresponding curvature $\Omega$
for the connection
is just a corresponding 2-cocycle
for the extension (2.11.10)
defined by (2.2) above.
\paragraph
Once an $L$-connection has been chosen, there is a more direct
(and less invariant) construction of the extension
(2.11.10). 
We shall need it later and therefore give the details now:
\paragraph
Let $\omega \colon L \to \roman{End}(M)$ be an
$L$-connection on $M$ with curvature $\Omega$.
As an $A$-module, let
$$
\roman{DO}(A,\omega)  = \roman{End}_A(M) \oplus L;
\tag2.12.1
$$
furthermore, 
we extend the bracket on
$\roman{End}_A(M)$
to one on
$\roman{DO}(A,\omega)$ by
means of
$$
\aligned
[(0,\alpha),(0,\beta)]
&=
(\Omega(\alpha,\beta),[\alpha,\beta]),
\quad \alpha,\,\beta \in L,
\\
[(\mu,0),(0,\alpha)]
&=
([\mu,\omega(\alpha)],0),
\quad \mu \in \roman{End}_A(M),\,\alpha\in L,
\endaligned
\tag2.12.2
$$
where at first the commutator $[\mu,\omega(\alpha)]$
is taken in 
$\roman{End}(M)$.
A little thought reveals that
(1.2.b) entails indeed
$[\mu,\omega(\alpha)] \in \roman{End}_A(M) \subseteq \roman{End}(M)$,
and that
the Bianchi identity (2.10) says 
that this bracket satisfies the
Jacobi identity, whence we obtain a
Lie algebra.
Furthermore, (1.2.b) implies
 that 
$\roman{End}_A(M)$ is a Lie ideal in
$\roman{DO}(A,\omega)$ -- in fact,
these two properties are equivalent --,
 and with the obvious morphisms, 
$$
0 \to  \roman{End}_A(M) \to \roman{DO}(A,\omega) \to L \to 0
\tag2.12
$$
then constitutes an extension of
Lie algebras.
We next define an action of
$\roman{DO}(A,\omega)$ on $A$ through the projection onto $L$;
then 
it is not hard to see
that (2.12)
is indeed an extension of
$(R,A)$-Lie algebras.
For example, (1.2.a) implies that there is no problem with the $A$-module
structures.
Finally, 
the obvious morphism
$$
\iota + \omega
\colon
\roman{DO}(A,\omega) = \roman{End}_A(M) \oplus L
\to
\roman{End}(M),
\tag2.12.3
$$
where
$\iota \colon \roman{End}_A(M) \to \roman{End}(M)$
refers to the obvious embedding,
endows $M$ with the structure of an
$(A,\roman{DO}(A,\omega))$-module,
and it is clear that if we take the obvious section
$L \to \roman{DO}(A,\omega)$
of $A$-modules
and define the corresponding 2-cocycle
by (2.2),
this 2-cocycle is just the curvature $\Omega$ of our
$L$-connection.
\paragraph
It is clear that a choice of $\omega$
induces a congruence isomorphism
$$
(\roman {Id},\cdot,\roman {Id})
\colon \roman{(2.11.10)} \longrightarrow \roman{(2.12)},
$$
and hence, up to congruence of extensions of $(R,A)$-Lie algebras,
an extension of the kind (2.12)
depends only on $A$ and $M$ and {\it not\/} on a particular choice of $\omega$.
This can also be seen in the following more direct and less invariant way:
Let $\omega' \colon L \to \roman{End}(M)$ be another
$L$-connection on $M$, and let
$$
0 \to  \roman{End}_A(M) \to \roman{DO}(A,\omega') \to L \to 0
\tag2.12'
$$
be the corresponding extension of
$(R,A)$-Lie algebras.
The rule (1.1.b) implies at once
that the difference $u = \omega - \omega'$
factors through 
$\roman{End}_A(M)$,
i.~e. that $u$ may be displayed as
$u \colon L \to \roman{End}_A(M)$.
 It is then straightforward to check that the rule
$(0,\alpha) \mapsto (u(\alpha),\alpha)$
induces a congruence isomorphism
$
(\roman {Id},\cdot,\roman {Id})
\colon \roman{(2.12')} \longrightarrow \roman{(2.12)}.
$
Summarizing we spell out the following.
\proclaim{Proposition 2.13}
For  an $(R,A)$-Lie algebra $L$ and an $A$-module $M$
there is,
up to congruence of extensions, at most one
extension of
the kind {\rm (2.12)}.
\endproclaim
Next we indicate how the classical argument 
(see e.~g. {\smc Koszul} [53])
yields the existence of connections
for projective $A$-modules $M$,
indeed relatively projective ones:
Let $M$ be a relatively free $A$-module, 
and write $M = M^{\sharp} \otimes A$.
Then 
an
$L$-connection 
$\omega \colon L \to \roman{End}(M)$ (with respect to $M^{\sharp}$)
is given by
$$
(\omega(\alpha)) (b) = 0,\quad (\omega(\alpha))(ab) = \alpha(a)b,
\quad a \in A, b \in M^{\sharp}.
\tag2.14.1
$$
Likewise, if $M$ is relatively projective, let $N$ be an $A$-module so that
$M \oplus N$ is relatively free,  let         
$$
D^{\oplus} \colon
M \oplus N
\to
\Hom_A(L,M \oplus N)
$$
be the adjoint (2.8.a')
of a connection 
of the kind (2.14.1),
and let $D$ be the composite
$$
D
\colon
M
@>>>
M \oplus N
@>{D^{\oplus}}>>
\Hom_A(L,M \oplus N)
@>>>
\Hom_A(L,M), 
\tag2.14.2
$$
where the unlabelled arrows are the obvious morphisms.
Then the adjoint
$$
\omega \colon
L
\to
\roman{End}(M)
\tag2.14.3
$$
is an $L$-connection on $M$.
\paragraph
Before we spell out the following, we note
that for a projective rank one  module $M$
the algebras $\roman{End}_A(M)$ and $A$ are canonically isomorphic.
As usual, we denote the Picard group of $A$
by $\roman{Pic}(A)$;
we remind the reader that it consists of classes of
projective rank one   modules,
with addition being induced by the tensor product.
\proclaim{Theorem 2.15}
Let $L$ be an $(R,A)$-Lie algebra.
Then the assignment 
to the class $[M] \in \roman{Pic}(A)$ of a projective rank one   module $M$
of the class $[\Omega_M] \in \Ho^2(\roman{Alt}_A(L,A))$
of the curvature of an $L$-connection on $M$
is a homomorphism
$$
\roman{Pic}(A)
\longrightarrow \Ho^2(\roman{Alt}_A(L,A))
\tag2.15.1
$$
of $R$-modules.
\endproclaim

\noindent
{\smc REMARK after publication.}
The map (2.15.1) is only a homomorphism of abelian groups.
This is enough to validate what follows.

\demo{Proof}
Theorem 2.13 implies that
(2.15.1) is well defined.
To see that it is a homomorphism
of $R$-modules,
let $M_1$ and $M_2$ be projective rank one  modules,
and let\linebreak
$\omega_1 \colon L \to \roman{End}(M_1)$ and
$\omega_2 \colon L \to \roman{End}(M_2)$
be $L$-connections,
with curvatures $\Omega_1$ and $\Omega_2$, respectively.
For $\alpha \in L,\, x_1 \in M_1,\,x_2 \in M_2$, let
$$
(\omega(\alpha))(x_1 \otimes x_2)
=
\left((\omega_1(\alpha))(x_1)\right) \otimes x_2
+
x_1 \otimes \left((\omega(\alpha))(x_2)\right).
$$
This
defines
an $L$-connection on $M_1 \otimes M_2$.
It is easy to see that it has
 curvature $\Omega_1 + \Omega_2$. \qed
\enddemo
For illustration, let $N$ be a smooth real manifold, let
$A$ be the algebra of smooth complex functions on 
$N$, and let $L$ be the Lie algebra of 
smooth complex vector fields
on $N$, i.~e.
that of  smooth sections of the complexified
tangent bundle $T^{\bold C}N $,
where $\bold C$ refers to the complex numbers.
 Then 
$\roman{Pic}(A)$
is the group of classes of complex line bundles,
the assignment 
of its first chern class to a line bundle
yields an isomorphism
$
\roman{Pic}(A) \longrightarrow \Ho^2(N,\MZ),
$
and the morphism
(2.15.1)
is part of the corresponding Chern-Weil map.
\paragraph
\noindent
{\smc Remarks 2.16.} Some historical comments 
about the algebraic approach to connection theory
seem appropriate.
\roster
\item
Let
$R$ be the ring of the reals or that of the complex numbers,
let
$N$ be a  manifold,
either smooth or complex analytic,
with tangent bundle $TN$, let
 $A$ be the algebra of smooth or analytic functions on 
 $N$, let $L$ be the Lie algebra of 
vector fields
on $N$, either smooth or complex analytic,
let
 $\zeta$ be a smooth or complex analytic vector bundle on $N$,
and let $M= \Gamma(\zeta)$ be the 
$A$-module
of sections of $\zeta$,
either smooth or complex analytic.
Moreover, let
$\xi \colon P \to N$ be a principal bundle for $\zeta$,
with structure group $G$ and Lie algebra $g$,
and consider the corresponding extension
$$
0
@>>>
V
@>>>
TP
@>>>
TN \times_N P
@>>>
0
\tag2.16.1
$$
of 
vector bundles over $P$,
where $V$ is the vertical subbundle; the latter is 
well known to be
canonically
isomorphic to the trivial bundle $g \times P$.
A treatment of the notions of connection and curvature
by pure algebra
goes back at least to {\smc Cartan} [13],
and  Cartan's notions of algebraic connections and curvature
are derived from formal properties of (2.16.1)
and its spaces of sections.
\item
The sequence (2.16.1) inherits an obvious
$G$-action;
when we divide it out  we obtain
an extension
$$
0
@>>>
g(\xi)
@>>>
TP/G
@>>>
TN
@>>>
0
\tag2.16.2
$$
of vector bundles over $N$,
where
$g(\xi) = (V/G)$
is the bundle associated to the principal bundle
by the adjoint representation.
Then the spaces
$\Gamma(g(\xi))$ and
$\Gamma(TP/G)$
of sections
inherit obvious structures of
Lie algebras, 
in fact of $(R,A)$-Lie algebras,
and
the corresponding sequence of sections
$$
0
@>>>
\Gamma(g(\xi))
@>>>
\Gamma(TP/G)
@>>>
\Gamma(TN)
@>>>
0
\tag2.16.3
$$
is an extension of  $(R,A)$-Lie algebras;
moreover, 
the obvious action of
$\Gamma(TP/G)$
 on the space $\Gamma(\zeta)$
of sections for $\zeta$
endows the latter with a structure of
 an $(A,\Gamma(TP/G))$-module
in our sense.
Indeed, when $G$ is the full linear group on the fibre of
$\zeta$, the sequence (2.16.3) is just the above sequence
(2.11.10),
and the sections of
$TP/G$
are exactly the infinitesimal gauge transformation
of $\zeta$ 
covering infinitesimal diffeomorphisms of $N$
in the usual sense
(which correspond bijectively to the $G$-invariant vector fields on $P$).
In general, the principal bundle comes with a linear representation of
$G$ 
on the fibre,
and the structure induces in particular a morphism
$$
\CD
0
@>>>
\Gamma(g(\xi))
@>>>
\Gamma(TP/G)
@>>>
\Gamma(TN)
@>>>
0
\\
@.
@VVV
@VVV
@V{\roman{Id}}VV
@.
\\
0
@>>>
\roman{End}_A(M) 
@>>> 
\roman{DO}(A,L,M)
@>>>
L
@>>>
0
\endCD
$$
of extensions of $(R,A)$-Lie algebras.
\paragraph
The sequence (2.16.2) was introduced
by {\smc Atiyah} [5] (Theorem 1)
and is now usually called
the {\it Atiyah sequence \/} of the principal 
bundle $\xi$.
Furthermore,
the classical notions of linear connection and curvature
have been described by {\smc Atiyah} [5] in formally the same way as
in the present Section.
In particular,
when we write
$L=\Gamma(TN)$, we see
that for any vector bundle
$\zeta$ the $A$-module
$\Gamma(\zeta)$ of sections
is $L$-normal
in the above sense
in both the smooth and complex analytic
context;
however, 
while the exact sequence
(2.16.3)
will always split in the smooth case,
it will in general {\it not\/}
split in the complex analytic case;
it will split e.~g. over a Stein manifold,
but a counterexample is given e.~g. in Proposition 22
of {\smc Atiyah} [5].
For a complete account to Atiyah sequences see
App. A in {\smc Mackenzie} [64].
\item
Notions of algebraic connections and curvature 
of the above kind
(corresponding to the linear notions in the classical case)
have been introduced by {\smc Koszul} [53].
However, 
we do not know whether our
extension theory,
the 
infinitesimal gauge algebra,
and the related material 
are already in the literature.
\item
The formal properties 
of the concept of an Atiyah sequence
have 
been incorporated in
the more general concept of
a {\it Lie algebroid\/};
this is the geometric analogue of an 
$(R,A)$-Lie algebra.
{\smc Pradines} [78]
associated to any differentiable groupoid
a Lie algebroid
as a first order infinitesimal invariant.
This 
generalizes the
construction of the Atiyah sequence
of a principal bundle (reproduced above).
For a complete account to
differentiable groupoids,
Lie algebroids,
and connection theory
while staying in the category of
finite rank vector bundles
see {\smc Mackenzie} [64];
his geometric notions are
similar to the above  algebraic
ones of
 connection and curvature.
Further relevant references are
{\smc Almeida-Molino} [2],
{\smc Coste-Dazord-Weinstein} [17],
{\smc Weinstein} [102] -- [106].
\endroster

\bigskip
\noindent
{\bf 3. Poisson algebras}
\medskip
\noindent
As before, let $R$ be a commutative ring with $1$,
and let
$A$ be a commutative $R$-algebra.
Recall that
a {\it Poisson algebra\/} structure
on $A$ 
(over $R$)
is a Lie bracket
$
\{\cdot,\cdot\} \colon A \otimes A 
\longrightarrow
A
$
on $A$ so that, for every
$a,b,c \in A$,
$$
\{ab,c\} = a\,\{b,c\}+\{a,c\}\,b.
\tag3.1             
$$
We note that this implies at once
that, for
$a \in A$ and $r \in R$,
we have
$\{a,r\} = 0$.
For example,
the ring of smooth functions on
 a smooth Poisson manifold,
see e.~g. {\smc Weinstein\/} [102], [103],
is a Poisson algebra.
In the present Section we shall 
characterize
such a structure on an arbitrary $R$-algebra $A$
in terms of a suitable
$(R,A)$-Lie algebra structure 
together with an appropriate closed 2-form
on the $A$-module
$D_A$ of K\"ahler differentials
for $A$.
\paragraph
For convenience we review briefly
the construction of
K\"ahler differentials:
Let
$A$ be a commutative $R$-algebra, and let
$I = \roman{ker}(\mu)\colon A \otimes A \to A$
so that, with the obvious $A \otimes A$-module structures,
$
  0
@>>>
I
@>>>
A \otimes A
@>{\mu}>>
A
@>>>
0
$
is
an extension
of
$A \otimes A$-modules.
As an $A \otimes A$-module,
$I$ is generated by the elements
$a \otimes 1 - 1 \otimes a \in A \otimes A$;
hence as an $R$-module, $I$ 
is generated by elements of the form
$
b(a\otimes 1 - 1 \otimes a)c,\ a,b,c \in A.
$
Let $D_A$ be the
$A$-{\it module of 
formal differentials\/}
or
{\it K\"ahler differentials\/}
for $A$,
cf. e.~g. {\smc Kunz} [55];
as an $A$-module, it is generated by
elements
$
da,\, a \in A,
$
subject to the relations
$$
d(bc) = (db)c + b dc, \quad dr =0,
\tag3.2.1
$$
where $b,c \in A,\  r \in R$.
It is well known that
the rule
$
\quad a \otimes 1 - 1 \otimes a
\mapsto 
da \quad
$
induces an isomorphism
$
\roman {Tor}^{A \otimes A}_1(A,A)
=
I/I^2
\longrightarrow
D_A
$
of $A$-modules.
In the standard way, $D_A$ represents the functor
$\roman{Der}(A,-)$
from the category of $A$-modules to itself.
More precisely,
given an $A$-module $M$ and
an element $h \in \Hom_{A \otimes A}(I,M)$,
the $R$-endomorphism
$d_h \colon A \to A$
defined by
$
d_h(a) = h(a\otimes 1 - 1 \otimes a)
$
is a derivation $A \to M$,
and it is well known that
the rule
$h \mapsto d_h$
induces a natural isomorphism
$$
\Hom_A(D_A,M)
\longrightarrow
\roman{Der}(A,M)
\tag3.2.2
$$
of $A$-modules;
in fact, this property {\it characterizes\/}
the K\"ahler differentials.
In particular, the $A$-module 
$\roman{Hom}_A(\roman{Der}(A),A)$
may canonically be identified with the double dual $D_A^{**}$
of $D_A$ and there is a canonical map
$$
D_A \longrightarrow \roman{Hom}_A(\roman{Der}(A),A)
\tag3.2.3
$$
of $A$-modules.
Moreover,
with the usual differential $d$
given by
$
d(a db_1 \cdots db_k) =da db_1 \cdots db_k,
$
the graded exterior algebra
 $ \Lambda_A[D_A]$
is the standard differential graded commutative algebra of
K\"ahler forms;
it is natural in $A$.
Henceforth we write
$
\Lambda A = (\Lambda_A[D_A],d).
$
It is well known that
(3.2.3)
induces a morphism
$$
\Lambda A 
\longrightarrow
\roman{Alt}_A(\roman{Der}(A),A)
\tag3.2.4
$$
of differential graded algebras.
It is proved in 
{\smc Hochschild-Kostant-Rosenberg} [40] that,
for a regular affine algebra $A$ over a perfect field,
(3.2.3)
and hence (3.2.4)
are isomorphisms.
Hence there is then
{\it no need to distinguish between 
formal differentials and differential forms\/}.
On the other hand, when
$A$ is the algebra of smooth functions on a smooth finite dimensional manifold
$N$,
the algebra $\Lambda A$
of formal differentials
will be much bigger than
the usual algebra
$\roman{Alt}_A(\roman{Der}(A),A)$
of differential forms.
For example, in the case 
$N = \bold R$, the formal differential
$df - f'dt \in D_A$ will be non-zero
when $f$ and $t$ are algebraically independent,
see e.~g. p. 27 of {\smc Krasilsh'chik, Lychagin, and Vinogradov} [54].
We shall say more about this in (3.12) below.
\paragraph
Let $L$ be an $(R,A)$-Lie algebra.
As in the classical case,
by functoriality, 
we can extend the action of $L$ on $A$ to 
one of $L$ on
$D_A$, in fact on $\Lambda A$,
cf. e.~g. {\smc Lang} [56].
Explicitly, given
 $X \in L$ and $\alpha \in D_A$,
for $a,\,b \in A$, let
$$
\lambda_X(bda) = (X(b))da + b d(X(a)).
\tag3.3
$$
This endows $D_A$ with a structure of an $(A,L)$-module
and it is clear that this extends to a structure 
$
\omega \colon L
\longrightarrow
\roman{End}(\Lambda A)
$
of an
$(A,L)$-module
on 
the differential graded algebra
$\Lambda A$
of K\"ahler differentials,
in fact to that of an algebra over $(A,L)$.

\smallskip\noindent
{\smc REMARK after publication.}
This yields only $L$-module, not
$(A,L)$-module structures, on $D_A$ and $\Lambda A$.
This oversight is irrelevant for the rest of the paper.
\smallskip

Given $X \in L$ and $\alpha \in \Lambda A$,
we refer to
$\lambda_X(\alpha)$  as the
{\it Lie derivative\/}
of $\alpha$ with respect to $X$.
It is clear that when $A$ is the ring of smooth functions
on a smooth finite dimensional manifold $N$
and
 when $L$ is the
$(\bold R,A)$-Lie algebra of smooth vector fields on
$N$
so that 
$\roman{Alt}_A(\roman{Der}(A),A)$
coincides with the usual de Rham complex,
the morphism (3.2.4)
is compatible with
Lie derivatives,
where
on the right hand side $\roman{Alt}_A(\roman{Der}(A),A)$
the usual Lie derivative is understood.
It is also clear that,
for an arbitrary algebra $A$ and an arbitrary
$(R,A)$-Lie algebra $L$,
 the usual formula
 yields an
operation of Lie derivative 
$
\omega \colon L
\longrightarrow
\roman{Alt}_A(D_A,A)
$
for multlinear alternating forms 
on $D_A$
with values in $A$;
for example, for 
$X \in L$ and
a bilinear alternating form
$\Omega \colon D_A \otimes_A D_A \to A$ this formula reads
$$
\left(\lambda_X(\Omega)\right)(\alpha,\beta)
=
X(\Omega(\alpha,\beta)) -
\Omega(\lambda_X(\alpha),\beta)
-
\Omega(\alpha,\lambda_X(\beta)).
\tag3.4
$$
\proclaim{Lemma 3.5}
Let
$
\{\cdot,\cdot\} \colon A \otimes A 
\longrightarrow
A
$
be
a Poisson structure
on $A$,
and, for $a,b,u,v \in A$, let
$$
\pi
(adu \otimes bdv)
=
ab\{u,v\} \in A.
\tag3.5.1
$$
Then $\pi$ is an alternating $A$-bilinear form 
$$
\pi
\colon
D_A
\otimes_{A}
D_A
\longrightarrow
A.
\tag3.5.2
$$
\endproclaim
When we wish to emphasize the dependence of
$\pi$ 
on
$\{\cdot,\cdot\}$
we shall write
$\pi_{\{\cdot,\cdot\}}$.
We refer to
$\pi_{\{\cdot,\cdot\}}$
as the {\it Poisson\/} (2-) {\it form\/}
of $(A,\{\cdot,\cdot\})$.
In view of (3.15) below, this form
generalizes the symplectic form of a symplectic manifold.
\demo{Proof of 3.5}
As an $A$-module, 
$D_A$
is generated by
the elements
$
da, \,a \in A,
$
subject to the relations (3.2.1).
On the other hand,
by definition,
for $a,b,c \in A,\ r \in R$, we have
$$
\aligned
\pi
(d(ab) \otimes dc)
&=
\{ab,c\}=
a\{b,c\}+
b\{a,c\} 
=
\pi
((bda + adb)\otimes dc)
\in A,
\\
\pi
(da \otimes dr)
&=
\{a,r\}= 0,
\endaligned
$$
whence $\pi$ is indeed well defined. \qed
\enddemo
Given a Poisson structure
$
\{\cdot,\cdot\} \colon A \otimes A 
\longrightarrow
A
$
on $A$,
let
$
\pi
\colon
D_A \otimes _A D_A
\to A
$
be the corresponding form (3.5.2),
and let
$$
\pi^{\sharp}
\colon
D_A
\longrightarrow
\Hom_{A}(D_A,A) =
\roman{Der}(A)
\tag3.6
$$
be the indicated adjoint of $\pi$;
then $\pi^{\sharp}$  is a morphism of $A$-modules.
The corresponding adjoint $D_A \otimes A \to A$
will be denoted by
$$
\pi^{\flat}\colon D_A \otimes A \longrightarrow A.
\tag3.7
$$
We shall say that (3.5.1) is {\it non-degenerate\/}, if 
$\pi^{\sharp}$ is injective.
For $\alpha \in D_A$, we shall often write
$
\alpha^{\sharp} = \pi^{\sharp}(\alpha) \in 
\roman{Der}(A).
$
Further, given $a \in A$,
we refer to the derivation
$
(da)^{\sharp} = \{a,-\}\colon A \to A
$
as the corresponding
{\it Hamiltonian element\/}.
\proclaim{Theorem 3.8}
Let
$(A,\{\cdot,\cdot\})$ be a Poisson algebra,
let
$
\pi^{\sharp}
\colon
D_A
\to
\roman{Der}(A)
$
be the morphism {\rm (3.6)},
and,
for 
$a,b,u,v \in A$, let
$$
[adu,bdv] =
a\{u,b\}dv
+
b\{a,v\}du
+
abd\{u,v\} \in D_A.
\tag3.8.1
$$
Then $(\pi^{\sharp},[\cdot,\cdot])$
together with the $A$-module structure
endows
$D_A$ with a structure of an
$(R,A)$-Lie algebra 
in such a way that
$\pi^{\sharp}$
is a morphism of
$(R,A)$-Lie algebras.
\endproclaim
\demo{Proof}
We indicate at first that
$[\cdot,\cdot]$
is well defined.
To this end, consider
$$
\aligned
[d(ab),cdv]
&=
\{ab,c\}dv + c d\{ab,v\}
\\
&=
b\{a,c\}dv + 
a\{b,c\}dv +
cbd\{a,v\}
+
c\{a,v\}db +
cad\{b,v\} +
c\{b,v\}da.
\endaligned
$$
On the other hand,
$$
\aligned
[adb,cdv]
&=
a\{b,c\}dv + c \{a,v\} db + acd\{b,v\}
\\
[bda,cdv]
&=
b\{a,c\}dv + c \{b,v\} da + bcd\{a,v\}
\endaligned
$$
whence indeed
$
[d(ab),cdv]
=
[adb +bda,cdv]
$
as desired.
\paragraph
It is clear that the identities (1.1.a) and (1.1.b)
hold.
Furthermore, a calculation shows that
 the Jacobi identity
for $[\cdot,\cdot]$ boils down to
$$
[[du,dv],dw]
+
[[dv,dw],du]
+            
[[dw,du],dv]
=
0,
\tag3.8.2
$$
and 
it is clear that
(3.8.2)
is equivalent to the
Jacobi identity
for $\{\cdot,\cdot\}$.
Hence
$(\pi^{\sharp},[\cdot,\cdot])$
together with the $A$-module structure
is indeed
 a structure of an
$(R,A)$-Lie algebra 
on $D_A$.
\paragraph
Finally, to see that
$\pi^{\sharp}$ is a morphism 
of
$(R,A)$-Lie algebras,
let $a,b,c,u,v \in A$, and consider
$$
[adu,bdv]^{\sharp}(c)
=
a\{u,b\}\{v,c\}
+b\{v,a\}\{u,c\}
+ab\{\{u,v\},c\}.
$$
On the other hand,
$$
\aligned
(adu)^{\sharp}(bdv)^{\sharp}(c)
&=a\{u,b\{v,c\}\}
=a\{u,b\}\{v,c\}+ab \{u,\{v,c\}\},
\\
(bdv)^{\sharp}(adu)^{\sharp}(c)
&=b\{v,a\{u,c\}\}
=b\{v,a\}\{u,c\}+ab \{v,\{u,c\}\}.
\\
\endaligned
$$
In view of the Jacobi identity
$
\{\{u,v\},c\} = \{u,\{v,c\}\} - \{v,\{u,c\}\} 
$
we conclude that
$$
[adu,bdv]^{\sharp}(c)
=
(adu)^{\sharp}(bdv)^{\sharp}(c)
-
(bdv)^{\sharp}(adu)^{\sharp}(c)
$$
whence
$\pi^{\sharp}$ is indeed a morphism
of $(R,A)$-Lie algebras. \qed
\enddemo
\proclaim{Addendum 3.8.3}
In terms of the above notion {\rm (3.3)}
of Lie derivative
the formula {\rm (3.8.1)} may be rewritten
$$
[\alpha,\beta] =
\lambda_{\alpha^{\sharp}}\beta
-
\lambda_{\beta^{\sharp}}\alpha
-
d(\pi_{\{\cdot,\cdot\}}(\alpha,\beta)) \in D_A,\quad \alpha, \beta \in D_A.
\tag3.8.1'
$$
\endproclaim
This is straightforward and left to the reader. \qed
\paragraph
\proclaim{Addendum 3.8.4}
Given two Poisson algebras $(A, \{\cdot,\cdot\})$ and
$(A', \{\cdot,\cdot\}')$
and a morphism
$$
\phi \colon
(A, \{\cdot,\cdot\})
\longrightarrow
(A', \{\cdot,\cdot\}')
\tag3.8.5
$$
of Poisson algebras,
let
$
\phi_* \colon
D_A 
\longrightarrow
D_{A'}
$
be the induced $A$-linear morphism
between the K\"ahler differentials as indicated
where $A$ acts on 
$D_{A'}$ through $\phi$.
Then 
$(\phi,\phi_*)$
is a morphism
$$
(\phi,\phi_*) \colon
(A,D_{\{\cdot,\cdot\}})
\longrightarrow
(A',D_{\{\cdot,\cdot\}'})
\tag3.8.6
$$
of Lie-Rinehart algebras
(cf. Section 1).
More precisely, 
$\phi_*$
is a morphism
$
\phi_* \colon
D_{\{\cdot,\cdot\}}
\longrightarrow
D_{\{\cdot,\cdot\}'}
$
of Lie algebras over $R$
(and one of $A$-modules),
and the diagram
$$
\CD
D_A \otimes A @>{\pi^{\flat}_{\{\cdot,\cdot\}}}>> A
\\
@V{\phi_* \otimes \phi}VV
@V{\phi}VV
\\
D_{A'} \otimes A' @>{\pi^{\flat}_{\{\cdot,\cdot\}'}}>> A'
\endCD
$$
is commutative.
\endproclaim
This is again straightforward and left to the reader.
\paragraph
Henceforth we write
$D_{\{\cdot,\cdot\}}$
for $D_A$ together with the $(R,A)$-Lie algebra structure
given in (3.8) above.
It is clear that 
the machinery of Section 1 applies
to $D_{\{\cdot,\cdot\}}$.
In particular,
we have 
the differential graded commutative algebra
$
\roman{Alt}_A(D_{\{\cdot,\cdot\}},A);
$
we refer to its cohomology as
the {\it Poisson cohomology of\/} $(A,\{\cdot,\cdot\})$,
written
$\Ho^*_{\roman{Poisson}}(A,\{\cdot,\cdot\};A)$.
More generally, for any $(A,D_{\{\cdot,\cdot\}})$-module $M$, we have the
differential graded
$\roman{Alt}_A(D_{\{\cdot,\cdot\}},A)$-module
$\roman{Alt}_A(D_{\{\cdot,\cdot\}},M)$;
we refer to its cohomology as
the {\it Poisson cohomology of\/} $(A,\{\cdot,\cdot\})$
{\it with values in\/} $M$,
written
$\Ho^*_{\roman{Poisson}}(A,\{\cdot,\cdot\};M)$.
In view of (1.14),
when $D_A$ is projective as an $A$-module, we have
${
\Ho^*_{\roman{Poisson}}(A,\{\cdot,\cdot\};M)
=
\roman{Ext}^*_{U(A,D_{\{\cdot,\cdot\}})}(A,M),
}$
where 
$U(A,D_{\{\cdot,\cdot\}})$
refers to the corresponding universal algebra of differential operators 
introduced in Section 1.
A 
cochain, cocycle etc. in
$\roman{Alt}_A(D_{\{\cdot,\cdot\}},M)$
will be referred to as a
{\it Poisson\/} cochain or cocycle etc. as appropriate. 
Notice that Poisson cohomology depends only on the
$(R,A)$-Lie algebra structure on
the $A$-module $D_A$ of K\"ahler differentials
which is derived from $\{\cdot,\cdot\}$ by means of (3.6) and (3.8.1);
see also (3.11.3) below.
Likewise, 
for any right $(A,D_{\{\cdot,\cdot\}})$-module $M$, we have the
chain complex
\linebreak
$M \otimes _{U(A,D_{\{\cdot,\cdot\}})} K(A,D_{\{\cdot,\cdot\}})$
where
$K(A,D_{\{\cdot,\cdot\}})$
refers to the corresponding Koszul complex (1.13).
We refer to its homology as
the {\it Poisson homology of\/} $(A,\{\cdot,\cdot\})$
{\it with values in\/} $M$,
written
$\Ho_*^{\roman{Poisson}}(A,\{\cdot,\cdot\};M)$.
In particular, 
 with respect to the right 
$U(A,D_{\{\cdot,\cdot\}})$-module structure
(1.8.4) on $A$,
we have the  Poisson homology of $(A,\{\cdot,\cdot\})$
 with values in $A$.

\smallskip\noindent
{\smc Remark after publication.}
As pointed out after (1.8.4), the construction given there
does not work. Under the present circumstances, the requisite right
$U(A,D_{\{\cdot,\cdot\}})$-module structure
is given by the association
$$
a \otimes (b du) \mapsto \{ab,u\}
$$
and this is perfectly o.k.
This construction has been elaborated upon in follow up papers:
\smallskip\noindent
Lie-Rinehart algebras, Gerstenhaber algebras, and Batalin-Vilkovisky
algebras, Annales de l'institut Fourier 48 (1998) 425-440, dg-ga/9704005
\smallskip\noindent
Duality for Lie-Rinehart algebras and the modular class,
J. f\"ur die reine und angew. Mathematik 510 (1999) 103--159,
dg-ga/9702008
\smallskip\noindent
Differential Batalin-Vilkovisky algebras arising from twilled
Lie-Rinehart algebras, Banach center publications 51 (2000) 87--102
\smallskip\noindent
Twilled Lie-Rinehart algebras and differential Batalin-Vilkovisky
algebras, 
\linebreak
math.DG/9811069

\smallskip
For an arbitrary 
right $(A,D_{\{\cdot,\cdot\}})$-module $M$, when
$D_A$ is projective as an $A$-module, we have
${
\Ho_*^{\roman{Poisson}}(A,\{\cdot,\cdot\};M)
=
\roman{Tor}_*^{U(A,D_{\{\cdot,\cdot\}})}(M,A).
}$
\paragraph
In view of the naturality explained in (3.8.4),
these notions of Poisson cohomology and homology
are natural in the following sense:
Let
$\phi \colon(A, \{\cdot,\cdot\})\longrightarrow(A', \{\cdot,\cdot\}')$
be a morphism
of Poisson algebras.
By (3.8.4) it induces
a morphism
$$
(\phi,\phi_*) \colon
(A,D_{\{\cdot,\cdot\}})
\longrightarrow
(A',D_{\{\cdot,\cdot\}'})
$$
of Lie-Rinehart algebras, and hence morphisms
$$
U(\phi,\phi_*) \colon
U(A,D_{\{\cdot,\cdot\}})
\longrightarrow
U(A',D_{\{\cdot,\cdot\}'})
\tag3.8.7
$$
of algebras
and
$$
K(\phi,\phi_*) \colon
K(A,D_{\{\cdot,\cdot\}})
\longrightarrow
K(A',D_{\{\cdot,\cdot\}'})
\tag3.8.8
$$
of chain complexes, the second one being
$U(A,D_{\{\cdot,\cdot\}})$-linear in the obvious sense.
It is clear that, given right
$(A,D_{\{\cdot,\cdot\}})$- and
$(A',D_{\{\cdot,\cdot\}'})$-modules 
$M$ and $M'$ and a morphism
$\psi \colon M \to M'$
of right
$(A,D_{\{\cdot,\cdot\}})$-modules
(where the pair $(A,D_{\{\cdot,\cdot\}})$
acts on $M'$ through
$(\phi,\phi_*)$\,),
these combine to a morphism
$$
M \otimes _{U(A,D_{\{\cdot,\cdot\}})} K(A,D_{\{\cdot,\cdot\}})
\longrightarrow
M'\otimes _{U(A',D_{\{\cdot,\cdot\}'})} K(A',D_{\{\cdot,\cdot\}'})
\tag3.8.9
$$
and hence induce a morphism
$$
(\phi,\psi)_*
\colon
\Ho_*^{\roman{Poisson}}(A,\{\cdot,\cdot\};M)
\longrightarrow
\Ho_*^{\roman{Poisson}}(A',\{\cdot,\cdot\};M')
\tag3.8.10
$$
in Poisson homology.
Notice the special case where $M=A$ and $M' = A'$.
Likewise, given left
$(A,D_{\{\cdot,\cdot\}})$- and
$(A',D_{\{\cdot,\cdot\}'})$-modules 
$M$ and $M'$ and a morphism
$\psi \colon M' \to M$
(backwards)
of left
$(A,D_{\{\cdot,\cdot\}})$-modules
(where again the pair $(A,D_{\{\cdot,\cdot\}})$
acts on $M'$ through
$(\phi,\phi_*)$),
these combine to a morphism
$$
\roman{Alt}_{A'}(D_{\{\cdot,\cdot\}'},M')
\longrightarrow
\roman{Alt}_A(D_{\{\cdot,\cdot\}},M)
\tag3.8.11
$$
of chain complexes and hence induce a morphism
$$
(\phi^*,\psi_*)
\colon
\Ho^*_{\roman{Poisson}}(A',\{\cdot,\cdot\};M')
\longrightarrow
\Ho^*_{\roman{Poisson}}(A,\{\cdot,\cdot\};M)
\tag3.8.12
$$
in Poisson cohomology.
\paragraph

\paragraph

We have seen that a Poisson structure 
$\{\cdot,\cdot\}$
on $A$
determines a structure
of an $(R,A)$-Lie algebra on $D_A$.
We now examine the question to what extent
the latter determines the former.
\paragraph
\proclaim{Lemma 3.9}
The Poisson 2-form
$
\pi_{\{\cdot,\cdot\}}
\colon
D_A
\otimes_{A}
D_A
\longrightarrow
A
$
of
 a
 Poisson algebra $(A,\{\cdot,\cdot\})$
given by {\rm (3.5.1)} is a Poisson 2-cocycle,
i.~e. a 2-cocycle in $\roman{Alt}_A(D_{\{\cdot,\cdot\}},A)$.
Moreover, 
$\pi_{\{\cdot,\cdot\}}$  is natural in the obvious sense.
\endproclaim
\demo{Proof}
This is
an easy consequence of the Jacobi identity  for $\{\cdot,\cdot\}$. \qed
\enddemo
Hence for any Poisson algebra $(A,\{\cdot,\cdot\})$,
its Poisson 2-form
$\pi_{\{\cdot,\cdot\}}$
determines a class
$$
[\pi_{\{\cdot,\cdot\}}] \in \Ho^2_{\roman{Poisson}}(A,\{\cdot,\cdot\};A)
\tag3.10.1
$$
which is natural in an obvious sense;
we refer to it as the
{\it Poisson class\/} of $(A,\{\cdot,\cdot\})$.
This class generalizes the class
$[\sigma] \in \Ho^2(N,\bold R)$
of a symplectic structure $\sigma$ on a smooth real manifold $N$
to arbitrary Poisson algebras,
see (3.15) below.
Furthermore, we write $L_{\{\cdot,\cdot\}}$
for the corresponding $(R,A)$-Lie algebra
whose underlying $A$-module looks like
$A \oplus D_A$ and whose Lie structure is 
given by
$$
\aligned
[(0,\alpha),(0,\beta)]
&
=
(\pi_{\{\cdot,\cdot\}}(\alpha,\beta),[\alpha,\beta]),
\\
[(a,0),(0,\alpha)]
&
=
(-\alpha(a),0),
\endaligned
\tag3.10.2
$$
where $ a \in A,\,\alpha,\beta \in D_A$,
so that 
$L_{\{\cdot,\cdot\}}$
fits into an extension
$$
0
@>>>
A
@>>>
L_{\{\cdot,\cdot\}}
@>>>
D_{\{\cdot,\cdot\}}
@>>>
0
\tag3.10.3
$$
of $(R,A)$-Lie algebras
which, in view of (2.6),
is classified by 
$[\pi_{\{\cdot,\cdot\}}]$.
The extension (3.10.3) is natural in Poisson structures.
When the class $[\pi_{\{\cdot,\cdot\}}]$ is zero, a 
Poisson
1-form
$\vartheta
\colon
D_{\{\cdot,\cdot\}}
\to
A
$
so that
$
d\vartheta =
\pi_{\{\cdot,\cdot\}}
$
will be referred
as a {\it Poisson potential\/}. 
In view of (3.15) below, this generalizes the notion of
a symplectic potential, see e.~g. p. 9 of {\smc Woodhouse} [108].
It admits yet another interpretation:
Let $\vartheta \colon D_A \to A$ be 
a 1-form on $D_A$,
viewed as
a derivation
$X = X_{\vartheta} \colon A \to A$,
cf. (3.2.2).
Then the 
Lie derivative (3.4) of
the Poisson 2-form
 $\pi_{\{\cdot,\cdot\}}$
with respect to $X$, which is an element of the $(R,A)$-Lie algebra
$\roman{Der}(A)$,
boils down to
$$
\lambda_X(\pi_{\{\cdot,\cdot\}}) = - d{\vartheta} \in 
\roman{Alt}_A(D_{\{\cdot,\cdot\}},A).
\tag3.10.4
$$
In fact, 
$$
\align
\left(\lambda_X(\pi_{\{\cdot,\cdot\}})\right)(da,db)
&=
X(\{a,b\}) -
\{X(a),b\}
-
\{a,X(b)\}
\\
&=
\vartheta(d\{a,b\}) -
\{\vartheta(da),b\}
-
\{a,\vartheta(db)\}
\\
&= - (d\vartheta)(da,db).
\endalign
$$
In particular,
when $\vartheta$ is a Poisson potential
so that
$
d\vartheta =
\pi_{\{\cdot,\cdot\}}$,
we obtain
$$
\lambda_{-\vartheta}(\pi_{\{\cdot,\cdot\}}) = \pi_{\{\cdot,\cdot\}} \in 
\roman{Alt}_A(D_{\{\cdot,\cdot\}},A),
\tag3.10.5
$$
i.~e. up to sign the Lie derivative
of $\pi_{\{\cdot,\cdot\}}$
with respect to
its Poisson potential
is just $\pi_{\{\cdot,\cdot\}}$.
I am indebted to A. Weinstein for informing me that
this holds for a Poisson manifold and hence for prompting
me to examine the general case.
In the manifold case a vector field $X$ 
corresponding to a Poisson potential
is called a {\it conformal Poisson\/}
or {\it Liouville\/} vector field.
\paragraph
\noindent
{\smc Remark 3.10.6.}
In view of (3.8.4),
a morphism
$\phi \colon (A,\{\cdot,\cdot\})
\to
(A',\{\cdot,\cdot\}')
$
of Poisson algebras
induces the two morphisms
$$
(\phi^*,\roman{Id}_*)
\colon
\Ho^*_{\roman{Poisson}}(A',\{\cdot,\cdot\};A')
\longrightarrow
\Ho^*_{\roman{Poisson}}(A,\{\cdot,\cdot\};A')
\tag3.10.7
$$
and
$$
(\roman{Id},\phi_*)
\colon
\Ho^*_{\roman{Poisson}}(A,\{\cdot,\cdot\};A)
\longrightarrow
\Ho^*_{\roman{Poisson}}(A,\{\cdot,\cdot\};A')
\tag3.10.8
$$
in Poisson cohomology.
It is clear that
the classes 
$[\pi_{\{\cdot,\cdot\}}]\in \Ho^2_{\roman{Poisson}}(A,\{\cdot,\cdot\};A)$
and
$[\pi_{\{\cdot,\cdot\}'}]\in \Ho^2_{\roman{Poisson}}(A',\{\cdot,\cdot\}';A')$
go to the same class
$$
(\phi^*,\roman{Id}_*)[\pi_{\{\cdot,\cdot\}'}]=
(\roman{Id},\phi_*)[\pi_{\{\cdot,\cdot\}}]
\in
\Ho^*_{\roman{Poisson}}(A,\{\cdot,\cdot\};A').
$$
In this sense the Poisson class
of a Poisson algebra is natural;
notice, however, we cannot in general
directly relate
$[\pi_{\{\cdot,\cdot\}}]$
and
$[\pi_{\{\cdot,\cdot\}'}]$
by means of a map sending one class to the other.
\proclaim{Theorem 3.11}
Let
$
\pi
\colon
D_A
\otimes_{A}
D_A
\longrightarrow
A
$
be an
alternating $A$-bilinear form 
on $D_A$,
let
$\omega \colon D_A \to \Hom_A(D_A,A) = \roman{Der}(A)$
be its adjoint,
for $a,b \in A$, let
$
\{a,b\} = \pi(da \otimes db),
$
for $a,b,u,v \in A$, let
$$
[adu,bdv] =
a\{u,b\}dv
+
b\{a,v\}du
+
abd\{u,v\} \in D_A,
\tag3.11.0
$$
and assume that
$(\omega,[\cdot,\cdot])$
together with the $A$-module structure
endows
$D_A$ with a structure of an
$(R,A)$-Lie algebra.
Then
$\{\cdot,\cdot\}$
endows
$A$ with that of a Poisson algebra
if and only if      
$\pi$ is a 2-cocycle in
$(\roman{Alt}_A(D_A,A),d)$.
\endproclaim
\demo{Proof}
We have already seen that the condition is necessary.
To see that it is also sufficient,
assume that $\pi$ is a 2-cocycle, and
let
$$
0 
@>>> 
A
@>>> 
A \oplus_{-\pi} D_A
@>>> 
D_A
@>>> 
0
\tag3.11.1
$$
be the corresponding extension
of 
$(R,A)$-Lie algebras, cf. Section 2 above,
in particular (2.6).
Inspection shows that 
the morphism
$$
\iota_{\pi} 
=
(\roman{Id},d)
\colon A \longrightarrow
A \oplus_{-\pi} D_A
\tag3.11.2
$$
is compatible with
the bracket operations, i.~e., given $a,b \in A$,
$
\iota_{\pi}\{a,b\}=
[\iota_{\pi}a,
\iota_{\pi}b].
$
Since 
$A \oplus_{-\pi} D_A$
is a Lie algebra, 
and
since $\iota_{\pi}$ is manifestly injective,
we conclude that
$\{\cdot,\cdot\}$ satisfies the Jacobi identity. 
Moreover,
the defining relations (3.2.1)
for the K\"ahler differentials
 and the
$A$-bilinearity of $\pi$
imply at once that $A$ satisfies (3.1).
Since $\{\cdot,\cdot\}$ is obviously skew symmetric, we are done.
\qed
\enddemo
\paragraph
\noindent
{\smc Remark 3.11.3.}
Under the circumstances of (3.11),
when $\pi$ is not a 2-cocycle,
we can only conclude that
$\{\cdot,\cdot\}$
induces on
$dA \subseteq D_A$
a structure of a Lie algebra.
On the other hand, the differential graded commutative algebra
$\roman{Alt}_A(D_A,A)$ 
over $R$ and hence its cohomology are still defined.
\paragraph
\noindent
{\smc 3.12. The geometric version}
\paragraph
\noindent
Our notions of
Poisson homology and cohomology are entirely algebraic.
In the special case where
$A$  is the ring of smooth
real functions 
on a smooth finite dimensional Poisson manifold 
suitable geometric notions of
Poisson homology and cohomology
may be found in the literature
{\smc Brylinski} [12],
{\smc Koszul} [52],
{\smc Lichnerowicz} [58].
To relate the various approaches we need some preparations.
\paragraph
Recall that for an $R$-algebra $A$ the $A$-module 
 $D_A$ of K\"ahler differentials represents the functor
$\roman{Der}(A,-)$ 
on the category of $A$-modules.
However, in the smooth case the appropriate category to work with is
that of \lq\lq geometric modules\rq\rq:
Let $N$ be a smooth finite dimensional manifold and let $A$ be
its ring of smooth functions.
Using a terminology introduced on
p. 27 of {\smc Krasilsh'chik, Lychagin, and Vinogradov} [54]
we shall say that
an $A$-module $M$ is {\it geometric\/},
if $\bigcap_{x \in N} \mu_x M = 0$, where
$\mu_x \subseteq A$ denotes the maximal ideal of functions
on $N$ that vanish on  $x \in N$.
The $A$-module 
 $D_A$ of K\"ahler differentials
is apparently {\it not\/} geometric;
however, its geometric analogue is the $A$-module
of smooth 1-forms on $N$ which, in view of the isomorphism
(3.2.2), may be identified with the double dual $D_A^{**}$;
henceforth we denote 
this $A$-module by $D_A^{\roman{geo}}$ -- 
the superscript \lq ${}^{\roman{geo}}$\rq\ stands for
\lq geometric\rq.
The obvious morphism
$$
q \colon D_A \longrightarrow D_A^{\roman{geo}}
\tag3.12.1
$$
of $A$-modules
which sends a formal differential to the corresponding differential
form
is manifestly surjective.
For example, in the case 
$N = \bold R$, the formal differential
$df - f'dt \in D_A$ will be non-zero
when $f$ and $t$ are algebraically independent
but   goes to zero in
$D_A^{\roman{geo}}$ under (3.12.1).
\paragraph
Let $N$ be
a smooth finite dimensional
manifold, let
$A$ be its ring of smooth
functions,
and let
$L$ be the $(\bold R, A)$-Lie algebra of smooth vector fields on
$N$.
Roughly speaking,
the geometric approach
to Poisson structures and Poisson cohomology
proceeds in formally the same way as above
with the K\"ahler differentials replaced by the
smooth 1-forms on $N$.
To explain it
we recall at first
{\smc Lichnerowicz'} [58]  characterization of
a Poisson structure 
on $N$
in terms of 
a certain 2-tensor:
Let
$
[\cdot,\cdot] \colon\Lambda_A[L]\otimes\Lambda_A[L]\longrightarrow\Lambda_A[L]
$
be the
\lq\lq Schouten product\rq\rq\  or \lq\lq Lagrangian concomitant\rq\rq\ 
 {\smc Schouten} [81], [82],
{\smc Nijenhuis} [76],
{\smc Kirillov} [46],
{\smc Koszul} [52].
Given a 2-tensor
 $G \in \Lambda^2_A[L]$,
let
 $\{\cdot,\cdot\} \colon A \otimes A \to A$
be defined by
$$
\{f,g\} = \pi_G(df,dg),\quad f ,g \in A,
\tag3.12.2
$$
where 
$\pi_G$ denotes the image of $G$
under the obvious isomorphism
$$
\phi \colon \Lambda_A[L]
\longrightarrow
\roman{Alt}_A(D_A^{\roman{geo}},A)
\tag3.12.3
$$
of graded $A$-algebras;
conversely,
given
$\{\cdot,\cdot\} \colon A \otimes A \to A$,
define $G$ by (3.12.2).
Then 
$\{\cdot,\cdot\}$
is a Poisson structure on $N$ if and only if
 $[G,G] = 0$, see
{\smc Lichnerowicz} [58].
Given a
Poisson structure
$\{\cdot,\cdot\}$,
we shall refer to
$\pi_G$ as the
{\it geometric\/} Poisson (2-) form
on $N$, and we denote by
$
\pi_G^{\sharp}
\colon
D_A^{\roman{geo}}
\to
L
$
the indicated adjoint.
Furthermore,
for $\alpha \in D_A^{\roman{geo}}$, we write
$$
\alpha^{\sharp} =\pi_G^{\sharp}(\alpha).
$$
The following is known (where
as before $\lambda$ refers to the usual Lie derivative):
\proclaim{Proposition 3.12.4}
Let
$N$ be a smooth Poisson manifold, let
$G \in \Lambda^2_A[L]$ be its 2-tensor,
and,
for $\alpha,\,\beta \in D_A^{\roman{geo}}$, let
$$
[\alpha,\beta] =
\lambda_{\alpha^{\sharp}}\beta
-
\lambda_{\beta^{\sharp}}\alpha
-
d(\pi_G(\alpha,\beta)) 
\in D^{\roman{geo}}_A,\quad \alpha, \beta \in D_A^{\roman{geo}}.
\tag3.12.5
$$
Then $(\pi_G^{\sharp},[\cdot,\cdot])$
together with the $A$-module structure
endows
$D_A^{\roman{geo}}$ with a structure of an
$(\bold R,A)$-Lie algebra 
in such a way that
$\pi_G^{\sharp}$
is a morphism of
$(\bold R,A)$-Lie algebras.
\endproclaim
\paragraph
\noindent
{\smc Remark 3.12.6.}
For
an arbitrary commutative algebra $A$ over an arbitrary ground ring $R$,
the obvious morphism
$$
\Lambda_A[\roman{Der}(A)]
\longrightarrow
\roman{Alt}_A(D_A,A)
\tag3.12.7
$$
is  always defined, but it will
in general {\it not\/} be an isomorphism;
furthermore, given a Poisson structure $\{\cdot,\cdot\}$ on $A$,
its (algebraic) closed 2-form
$\pi_{\{\cdot,\cdot\}} \in \roman{Alt}^2_A(D_A,A)$
introduced in (3.5)
is always defined but need {\it not\/}
come from 
$\Lambda_A[\roman{Der}(A)]$.
In the situation of (3.12.4)
the morphism (3.12.7) is just the composite
of (3.12.3)
with the morphism
$$
q^* \colon \roman{Alt}_A(D_A^{\roman{geo}},A)
\longrightarrow
\roman{Alt}_A(D_A,A),
\tag3.12.8
$$
induced by (3.12.1),
and under (3.12.7)
the 2-tensor $G$ is mapped to
$\pi_{\{\cdot,\cdot\}}$.
Moreover, we shall see in (3.12.13) below that
(3.12.8) is indeed an isomorphism;
N.~B. that such a remark makes sense 
since we are in the geometric case.
This explains 
in which sense 
the algebraic 2-form
$\pi_{\{\cdot,\cdot\}}$
generalizes Lichnerowicz' 
geometric 2-tensor for smooth finite dimensional Poisson manifolds.
\demo{Proof of 3.12.4}
It is shown in
(3.1) of {\smc Weinstein} [105] 
and
(III.2.1) of {\smc Coste-Dazord-Weinstein} [17] 
that 
$G$ yields a morphism $T^*N \to TN$
of vector bundles inducing
$\pi^{\sharp}$
and that this structure together with
the bracket yields
in fact
a structure of a {\it Lie algebroid\/}
on the cotangent bundle
$T^*N$.
This is just another way to spell out the assertion. \qed
\enddemo
Henceforth we write
$D_{\{\cdot,\cdot\}}^{\roman{geo}}$
for
$D_A^{\roman{geo}}$
together with the
$(\bold R,A)$-Lie algebra structure
given in (3.12.5) above.
It is clear that 
the machinery of Section 1 applies to
$D_{\{\cdot,\cdot\}}^{\roman{geo}}$ as well
and (3.8) -- (3.11) 
have precise geometric analogues.
In particular, for a smooth Poisson manifold $N$
we refer to the cohomology of
the differential graded commutative algebra
$
\roman{Alt}_A(D_{\{\cdot,\cdot\}}^{\roman{geo}},A)
$
as
the {\it geometric Poisson cohomology of\/} $(N,\{\cdot,\cdot\})$,
written
$\Ho^*_{\roman{Poisson}}(N,\{\cdot,\cdot\};\bold R)$.
More generally, for any 
vector bundle $\zeta$ over $N$ 
with a flat 
$D_{\{\cdot,\cdot\}}^{\roman{geo}}$-connection
we have the
differential graded
$\roman{Alt}_A(D_{\{\cdot,\cdot\}}^{\roman{geo}},A)$-module
$\roman{Alt}_A(D_{\{\cdot,\cdot\}}^{\roman{geo}},\Gamma(\zeta))$;
we refer to its cohomology as
the {\it geometric Poisson cohomology of\/} $(N,\{\cdot,\cdot\})$
{\it with values in\/} $\zeta$,
written
$\Ho^*_{\roman{Poisson}}(N,\{\cdot,\cdot\};\zeta)$.
In view of (1.14),
since $D_A^{\roman{geo}}$ is projective as an $A$-module, we have
${
\Ho^*_{\roman{Poisson}}(N,\{\cdot,\cdot\};\zeta)
=
\roman{Ext}^*_{U(A,D_{\{\cdot,\cdot\}}^{\roman{geo}})}(A,\Gamma(\zeta)),
}$
where 
$U(A,D_{\{\cdot,\cdot\}}^{\roman{geo}})$
refers to the corresponding universal algebra of differential operators 
introduced in Section 1.
Virtually the same argument as in (3.9) shows that
the geometric 2-form
$\pi_G$ 
of $N$
is a 2-cocycle and represents a class 
${
[\pi_G] \in \Ho^2_{\roman{Poisson}}(N,\{\cdot,\cdot\};\bold R).
}$
Likewise, 
with respect to the right 
$U(A,D_{\{\cdot,\cdot\}}^{\roman{geo}})$-module structure
(1.8.4) on $A$,
we have the chain complex
$A \otimes _{U(A,D_{\{\cdot,\cdot\}}^{\roman{geo}})} 
K(A,D_{\{\cdot,\cdot\}}^{\roman{geo}})$
where
$K(A,D_{\{\cdot,\cdot\}}^{\roman{geo}})$
refers to the corresponding Koszul complex (1.13);
this chain complex 
computes what we call the
 {\it geometric Poisson homology of\/} $(N,\{\cdot,\cdot\})$
{\it with values in\/} $\bold R$.
We write
$\Ho_*^{\roman{Poisson}}(N,\{\cdot,\cdot\};\bold R)$
for the latter.
Since $D_A^{\roman{geo}}$ is projective as an $A$-module, we have
${
\Ho_*^{\roman{Poisson}}(N,\{\cdot,\cdot\};\bold R)
=
\roman{Tor}_*^{U(A,D_{\{\cdot,\cdot\}}^{\roman{geo}})}(A,A).
}$
We note that there is a more general notion of
geometric Poisson homology with values in a suitable vector bundle,
exactly analogous
to that for the algebraic case explained earlier;
we refrain from giving the details.
\paragraph
It has been noted by several people that,
for the ring $A$  of smooth
functions on a smooth 
finite dimensional
Poisson
manifold $N$,
the formula
(3.12.5) yields a
 Lie bracket 
on the space of 1-forms of $N$.
To our knowledge
the first ones to spell out such a bracket were
{\smc Gelfand-Dorfman} [29] (p. 243),
but it is not shown, however, that 
this bracket satisfies
the Jacobi identity.
The bracket occurs also
in {\smc Magri and Morosi} [67], 
see also (2.2) in {\smc Magri-Morosi-Ragnisco} [68],
on p. 266 of {\smc Koszul} [52] 
(denoted by $[\cdot,\cdot]_w =[\cdot,\cdot]_{\Delta}$), 
in (3.22) of {\smc Karasev} [44]
(for a special class of Poisson manifolds),
in (2.1) of {\smc Bhaskara and Viswanath} [8],
in
{\smc Coste-Dazord-Weinstein} [17],
in {\smc Weinstein} [105],
and perhaps work of others;
however, it seems that only
in
(3.1) of {\smc Weinstein} [105] 
and
(III.2.1) of {\smc Coste-Dazord-Weinstein} [17] 
is it pointed out that the bracket yields
in fact
a structure of a {\it Lie algebroid\/}
over $N$ which is 
the geometric analogue of an $(\bold R,A)$-Lie algebra.
\paragraph
Next we relate the algebra with the geometry by means of the
following.
\proclaim{Lemma 3.12.9}
Let
$N$ be a smooth Poisson manifold, and let
$(A,\{\cdot,\cdot\})$
be its real Poisson algebra.
Then   
the surjection
$q \colon D_A \to D_A^{\roman{geo}}$
is 
a morphism  $q \colon D_{\{\cdot,\cdot\}} \to 
D_{\{\cdot,\cdot\}}^{\roman{geo}}$ 
of $(\bold R,A)$-Lie algebras.
In particular,
the adjoint
$\pi_{\{\cdot,\cdot\}}^{\sharp}
\colon D_A \longrightarrow \roman{Der}(A)$
factors through $q$. 
Moreover,
the 
corresponding algebraic and geometric Poisson 2-forms
$\pi_{\{\cdot,\cdot\}} \colon D_A \otimes _A D_A \to A$
and
$\pi_G \colon D_A^{\roman{geo}} \otimes _A D_A^{\roman{geo}} \to A$
are related by
$$
\pi_G (q \otimes q) = \pi_{\{\cdot,\cdot\}}
\colon D_A \otimes _A D_A \longrightarrow A.
$$
\endproclaim
\demo{Proof}
Inspection of the definitions of
$\pi_{\{\cdot,\cdot\}}$ in (3.5)
and
$\pi_G$ in (3.12.2)
reveals that
\linebreak
$\pi_G (q \otimes q) = \pi_{\{\cdot,\cdot\}}$.
Moreover this implies that
the adjoint
$\pi_{\{\cdot,\cdot\}}^{\sharp}
\colon D_A \longrightarrow \roman{Der}(A)$
factors through 
$q \colon D_A \to D_A^{\roman{geo}}$.
Finally, 
since the bracket on $D_A$ may as well be defined
by (3.8.1') and  
since the bracket on $D_A^{\roman{geo}}$ is defined
by (3.12.5),
the morphism
$q \colon D_A \to D_A^{\roman{geo}}$
is in fact a morphism
$
q \colon D_{\{\cdot,\cdot\}} \to D_{\{\cdot,\cdot\}}^{\roman{geo}}
$
of $(\bold R,A)$-Lie algebras. \qed
\enddemo
\proclaim{Corollary 3.12.10}
Let
$N$ be a smooth Poisson manifold, 
let
$(A,\{\cdot,\cdot\})$
be its real Poisson algebra, and
let $\zeta$ be a smooth vector bundle on $N$ with
a flat
$D_{\{\cdot,\cdot\}}^{\roman{geo}}$-connection.
Then the morphism
$q \colon D_A \to D_A^{\roman{geo}}$
induces morphisms
$$
\Ho^*_{\roman{Poisson}}(N,\{\cdot,\cdot\};\zeta)
\longrightarrow
\Ho^*_{\roman{Poisson}}(A,\{\cdot,\cdot\};\Gamma(\zeta))
\tag3.12.11
$$
and
$$
\Ho_*^{\roman{Poisson}}(A,\{\cdot,\cdot\};A)
\longrightarrow
\Ho_*^{\roman{Poisson}}(N,\{\cdot,\cdot\};\bold R)
\tag3.12.12
$$
of graded real vector spaces.
Furthermore,
the class
$[\pi_G] \in
\Ho^2_{\roman{Poisson}}(N,\{\cdot,\cdot\};\bold R)
$
goes to
the class
$[\pi_{\{\cdot,\cdot\}}] \in
\Ho^2_{\roman{Poisson}}(A,\{\cdot,\cdot\};\bold A).
$
\endproclaim
Somewhat surprisingly, we have the following.
\proclaim{Theorem 3.12.13}
Let
$N$ be a smooth Poisson manifold, 
let
$(A,\{\cdot,\cdot\})$
be its real Poisson algebra, and
let $\zeta$ be a smooth vector bundle on $N$ with
a flat
$D_{\{\cdot,\cdot\}}^{\roman{geo}}$-connection. 
Then
the morphism 
$$
q^* \colon \roman{Alt}_A(D_{\{\cdot,\cdot\}}^{\roman{geo}},\Gamma(\zeta))
\longrightarrow
\roman{Alt}_A(D_{\{\cdot,\cdot\}},\Gamma(\zeta))
\tag3.12.14
$$
induced by
$q \colon D_A \to D_A^{\roman{geo}}$
is an isomorphism
of real chain complexes.
Consequently
{\rm (3.12.11)}
is an isomorphism of real graded vector spaces.
\endproclaim
Hence 
in the smooth case there is no need
to distinguish
between
geometric and algebraic
Poisson cohomology
with coefficients in a smooth vector bundle.
\demo{Proof}
Since the morphism
$q \colon D_A \longrightarrow D_A^{\roman{geo}}$
of $A$-modules is surjective,
the induced morphism
(3.12.14)
is an injective morphism of graded $A$-modules
and in view (3.12.9)
 is compatible with the differentials.
We show it is also surjective.
\paragraph
To  see this we observe first that it suffices to
consider the case where $\zeta$ is a trivial line bundle
and hence, as an $A$-module,
$\Gamma(\zeta) \cong A$.
In fact, 
since
$\Gamma(\zeta)$
is a finitely generated projective $A$-module,
we can easily reduce to the case of a finitely
generated free $A$-module,
and from this we reduce further to the case of
a free $A$-module of dimension 1.
\paragraph
To see that
$
\roman{Alt}_A(D_A^{\roman{geo}},A)
\longrightarrow
\roman{Alt}_A(D_A,A)
$
is surjective,
we observe first that
we already know that
$$
\roman{\Hom}_A(D_A^{\roman{geo}},A)
\longrightarrow
\roman{\Hom}_A(D_A,A)
$$
is surjective, in fact an isomorphism.
Indeed, cf. (3.2.2),
$\roman{\Hom}_A(D_A,A)$
is canonically isomorphic to the
$A$-module $\roman{Der}(A)$ of derivations of $A$,
i.~e. to the smooth vector fields on $N$,
while
$\roman{\Hom}_A(D_A^{\roman{geo}},A)$
is the dual of the (projective)
$A$-module of smooth 1-forms which is again
canonically isomorphic 
to the smooth vector fields on $N$.
In other words,
in degree one the isomorphism amounts to the
classical fact that
on a smooth finite dimensional manifold
a derivation is always induced from a smooth vector field.
The general case is only slightly more complicated than this.
In fact, let
$$
\phi\colon (D_A)^{\otimes _A k}
\longrightarrow
A
$$
be an $A$-multilinear $k$-form on $D_A$ with values in $A$,
and let
$$
\phi^{\sharp}
\colon 
D_A 
\longrightarrow
\Hom_A((D_A)^{\otimes _A (k-1)},A)
$$
be its adjoint.
By induction we may assume that
the induced morphism
$$
\Hom_A((D_A^{\roman{geo}})^{\otimes _A (k-1)},A)
\longrightarrow
\Hom_A((D_A)^{\otimes _A (k-1)},A)
$$
is an isomorphism, and hence we can identify
$\Hom_A((D_A)^{\otimes _A (k-1)},A)$
with the space of sections of
the $(k-1)$'th tensor power 
$TN^{\otimes (k-1)}$
of the tangent bundle
$TN$.
In other words,
$\phi^{\sharp}$ may be viewed
as a differential operator
between the trivial line bundle on $N$ and
$TN^{\otimes (k-1)}$.
However, such a differential operator
is a section in the
smooth vector bundle
$\Hom(T^*N,TN^{\otimes (k-1)})$ over $N$,
and hence
$\phi^{\sharp}$
factors through
$D_A 
\longrightarrow
D_A^{\roman{geo}}$.
Consequently $\phi$
factors through
$(D_A^{\roman{geo}})^{\otimes _A k}$,
and hence the induced
morphism
$$
\Hom_A((D_A^{\roman{geo}})^{\otimes _A k},A)
\longrightarrow
\Hom_A((D_A)^{\otimes _A k},A)
$$
is an isomorphism.
This completes the proof. \qed
\enddemo
Next we relate our notions
of Poisson homology and cohomology
 with what is in the literature.
\proclaim{Theorem 3.13}
Let $N$ be
a smooth 
finite dimensional
Poisson
manifold,
and let
$A$ be its ring of smooth
real functions 
with Poisson structure $\{\cdot,\cdot\}$.
Then the geometric Poisson cohomology 
of $(N,\{\cdot,\cdot\})$ with values in $\bold R$
coincides with the Poisson cohomology of $N$
introduced by {\rm {\smc Lichnerowicz} [58]}.
In other words,
Lichnerowicz'
Poisson cohomology
is naturally isomorphic to
$\roman{Ext}^*_{U(A,D_{\{\cdot,\cdot\}}^{\roman{geo}})}(A,A)$
with respect to the corresponding universal algebra
$U(A,D_{\{\cdot,\cdot\}}^{\roman{geo}})$ of differential operators.
Furthermore,
this cohomology also coincides with the
(algebraic)
Poisson cohomology 
of $(A,\{\cdot,\cdot\})$ with values in $A$.
\endproclaim
A brief comparison of the definitions
shows that
the first statement
holds.
In fact, 
in the present special case,
the chain complex 
$\Hom_{U(A,D_{\{\cdot,\cdot\}}^{\roman{geo}})}
(K(A,D_{\{\cdot,\cdot\}}^{\roman{geo}}),A)$
which computes
$\Ho^*_{\roman{Poisson}}(N,{\{\cdot,\cdot\}};\bold R)
=\roman{Ext}^*_{U(A,D_{\{\cdot,\cdot\}}^{\roman{geo}})}(A,A)$
is 
the same as that introduced by
Lichnerowicz.
Thus our approach yields a description of
this cohomology in terms of standard homological algebra,
i.~e. as an Ext over a suitable ring.
The \lq furthermore\rq\ statement
is an immediate consequence of Theorem 3.12.13.
It shows that our notion of Poisson cohomology
generalizes 
Lichnerowicz' notion
for smooth Poisson manifolds
to arbitrary Poisson algebras.
\proclaim{Theorem 3.14}
Let $N$ be
a smooth 
finite dimensional
Poisson
manifold,
and let
$A$ be the ring of smooth
functions on $N$.
Then the geometric Poisson homology 
of $(N,\{\cdot,\cdot\})$ with values in $\bold R$
coincides with the canonical homology of $N$
introduced by {\rm {\smc Brylinski} [12]}.
In other words,
 Brylinski's
canonical homology
is naturally isomorphic to
$\roman{Tor}_*^{U(A,D_{\{\cdot,\cdot\}}^{\roman{geo}})}(A,A)$
with respect to the corresponding universal algebra
$U(A,D_{\{\cdot,\cdot\}}^{\roman{geo}})$ of differential operators.
\endproclaim
This is again seen by
a brief comparison of the definitions.
In the present special case,
the chain complex 
$A \otimes _{U(A,D_{\{\cdot,\cdot\}}^{\roman{geo}})}
K(A,D_{\{\cdot,\cdot\}}^{\roman{geo}})$
coincides with 
the chain complex introduced by {\smc Koszul} [52]
and used in
{\smc Brylinski} [12] for the definition of
canonical homology.
Thus our approach 
yields a description of
Brylinski's canonical homology in terms of standard homological algebra,
i.~e. as a Tor over a suitable ring.
However, unlike the Poisson cohomology case,
there is no reason why 
the two notions of Poisson
homology should coincide.
More precisely,
algebraic Poisson homology
is computed by $\Lambda_A[D_A]$, equipped with a suitable differential,
while
geometric Poisson homology
is computed by $\Lambda_A[D_A^{\roman{geo}}]$   
equipped with a suitable differential,
and the surjection
${
\Lambda_A[D_A]
\longrightarrow
\Lambda_A[D_A^{\roman{geo}}]   
}$
induced by (3.12.1)
is compatible with the differentials.
Notice that
as a graded $A$-algebra
$\Lambda_A[D_A^{\roman{geo}}]$
is just the de Rham complex of $N$.
Further, the obvious morphism
$
\iota
\colon
R 
\longrightarrow
A
$
gives rise to the morphism
$
(\iota,\roman {Id})
\colon
(R,D_{\{\cdot,\cdot\}}^{\roman{geo}})
\longrightarrow
(A,D_{\{\cdot,\cdot\}}^{\roman{geo}})
$
which, cf. Section 1,
induces a morphism
$
K(\iota,\roman {Id})
\colon
K(R,D_{\{\cdot,\cdot\}}^{\roman{geo}}) 
\longrightarrow
K(A,D_{\{\cdot,\cdot\}}^{\roman{geo}})
$
of the corresponding Koszul complexes.
This yields an obvious morphism from
ordinary 
{\smc Chevalley-Eilenberg} [16] Lie algebra homology
to
geometric Poisson homology
whose existence 
was observed by {\smc Brylinski} [12].
\paragraph
Let $N$ be a smooth manifold,
let $A$ be the algebra of smooth real valued functions 
on $N$, and
let $L$ be the $(\bold R,A)$-Lie algebra
of smooth vector fields.
As in the classical case, call
an $A$-bilinear alternating 2-form
$
\sigma
\colon L \otimes_A L \longrightarrow A
$
a {\it symplectic structure\/}
if its adjoint
$
\sigma^{\sharp}
\colon
L 
\longrightarrow
\Hom_A(L,A)
$
is an isomorphism of $A$-modules,
cf. e.~g. p. 28 of {\smc Marsden\/} [70].
If this is the case, and if 
$N$ is finite dimensional
so that (i) 
the $A$-module $D_{\{\cdot,\cdot\}}^{\roman{geo}}=\Gamma(T^*N)$ of
smooth 1-forms
is a finitely generated
projective $A$-module
and
(ii) the obvious map
$D_{\{\cdot,\cdot\}}^{\roman{geo}} \to \Hom_A(L,A)$
is a natural isomorphism of $A$-modules,
the inverse of $\sigma^{\sharp}$ induces a Poisson structure
on $A$.
Hence:
\proclaim{Proposition 3.15}
Let $(N,\sigma)$ be a 
smooth
finite dimensional symplectic manifold,
let $A$ be its algebra  of smooth real functions,
let $L$ be the
$(\bold R,A)$-Lie algebra
of smooth vector fields on $N$, and
let $\{\cdot,\cdot\}$ the the corresponding Poisson structure
on $A$.
Then the adjoint 
$
\sigma^{\flat} \colon L \to D_{\{\cdot,\cdot\}}^{\roman{geo}}
$
of $\sigma$ is an isomorphism of $(\bold R,A)$-Lie algebras.
\endproclaim
This entails at once the following.
\proclaim{Addendum to 3.13 (Lichnerowicz [58])}
For a smooth finite dimensional symplectic manifold
the de Rham and Poisson cohomologies coincide.
\endproclaim
Furthermore,
(3.15) explains the ${}^*$-operation
introduced in (2.2.2) of {\smc Brylinski} [12]:
In fact,
this operation is just the composite of
the Poincar\'e duality map with
the isomorphism
$
\roman{Tor}_*^{U(A,D_{\{\cdot,\cdot\}}^{\roman{geo}})}(A,A)
\to
\roman{Tor}_*^{U(A,L)}(A,A) = \Ho_*(N,\bold R)
$
induced 
by the inverse of $\sigma^{\flat}$ --
where we have to keep in mind that a symplectic manifold is orientable.
We note that,
for an arbitrary algebra $A$ over a commutative ring $R$, 
the concept of a symplectic structure
does not seem to extend properly 
whereas 
our algebraic approach to
Poisson structures works in complete generality.
\paragraph
We conclude this Section with a number of examples:
\paragraph
\noindent
{\smc Example 3.16.} 
Let $A$ be a commutative algebra, let
$L$ be an $(R,A)$-Lie algebra,
and let
$S_A [L]$ be the symmetric algebra on $L$ over $A$.
 Then the 
$L$-action on $A$ and the
bracket operation 
on
$L$ induce
an obvious Poisson structure
$$
\{\cdot,\cdot\}
\colon
S_A [L]
\otimes
S_A [L]
\longrightarrow
S_A [L]
\tag3.16.1
$$
on $S_A [L]$;
explicitly, this structure is determined by
$$
\alignat 2
\{\alpha,\beta\} &= [\alpha,\beta],\quad &&\alpha,\beta \in L,
\tag3.16.2.1
\\
\{\alpha,a\} &= \alpha(a) \in A,\quad &&a \in A,\,\alpha \in L,
\tag3.16.2.2
\\
\{u,vw\}&= \{u,v\}w + v\{u,w\},\quad &&u,v,w \in S_A[L].
\tag3.16.2.3
\endalignat
$$
We now study this example in the light of our notion of Poisson cohomology:
Write $S =  S_A [L]$, and let
$\pi \colon D_S \otimes _S D_S \to S$
be the 2-form (3.5.1);
as before, we write $D_{\{\cdot,\cdot\}}$ for $ D_S $
together with the $(R,S)$-Lie algebra structure.
We assert that $\pi$ is a Poisson coboundary,
i.~e. $\{\cdot,\cdot\}$ admits a Poisson potential.
As an algebra,
$S$ is generated by the elements of $A$ and
those of $L$.
Hence 
as an $S$-module, $D_S$ is generated by
the formal differentials $da,\, a \in A$, and
$d\alpha,\,\alpha \in L$, and
there is an obvious surjection 
$S \otimes_A D_A \oplus S \otimes_A L \to D_S$
given by
$ 1 \otimes da \mapsto da,\, a \in A,$ and
$ 1 \otimes \alpha \mapsto d\alpha,\,\alpha \in L,$
with a slight abuse of notation.
We assert that the 1-form
$\vartheta \colon D_S \to S$ 
given by
$$
\vartheta(da)=0,\, a \in A,
\quad
\vartheta(d\alpha)= 
\alpha,\, \alpha  \in L,
\tag3.16.3
$$
is a Poisson potential for $\{\cdot,\cdot\}$.
Indeed, from
the formula
$$
\align
(d\vartheta)(du,dv) &=  (du)^{\sharp}(\vartheta(dv))
                      - (dv)^{\sharp}(\vartheta(du))
                      - \vartheta[du,dv]
\\ 
                    &=  \{u,\vartheta(dv)\}
                      - \{v,\vartheta(du)\}
                      - \vartheta(d\{u,v\})
\endalign
$$
we conclude at once
$$
\pi =d\vartheta \in \roman{Alt}^2_S(D_{\{\cdot,\cdot\}},S).
\tag3.16.4
$$
\paragraph
When
$L$ is projective as an $A$-module,
so that, in view of Rinehart's result reproduced as (1.9) above,
the obvious map
$S_A[L] \to E^0(U(A,L))$
is an isomorphism,
the non-commutative Poisson algebra $U(A,L)$ of differential operators 
is
a {\it deformation\/} in the sense of {\smc Gerstenhaber} [29]
of the 
commutative Poisson algebra $(S_A[L],\{\cdot,\cdot\})$.
In the even more special case where
$N$ is
a smooth 
finite dimensional
manifold, $A$  the ring of smooth
functions on $N$,
and $L$ the Lie algebra of smooth vector fields on $N$,
$S_A[L]$ is the algebra of smooth functions on the cotangent bundle
of $N$ which are polynomial on each fibre;
these are called {\it polynomial observables\/},
see e.~g. p. 84 of {\smc Woodhouse\/} [108].
Furthermore, the
Poisson structure (3.16.1) then comes from the classical one
on the smooth functions on the cotangent bundle of $N$;
this Poisson structure 
is the \lq\lq symmetric concomitant\rq\rq\  of {\smc Schouten} [81], [82].
For these matters, see also p. 180 of {\smc Vinogradov and Krasil'shchik} [100]
and (3.3) of {\smc Braconnier} [10], [11].
When $N$ is finite dimensional, the relation
(3.16.4) is the algebraic analogue of the
usual local formula expressing the symplectic structure
of the cotangent bundle $T^*N$
as a coboundary 
$\sigma =d \bold {p dq}$
of the 1-form $\bold {p dq}$
written out in local coordinates $(\bold {q},\bold {p})$.
The relationship between
$S_A[L]$ and $U(A,L)$
is the heart of the geometric quantization program
initiated by 
{\smc I. Segal} [84], {\smc Kostant} [48], and {\smc Souriau} [93].
\paragraph
\noindent
{\smc Example 3.17.}
We now modify the above example: Under the circumstances of
(3.16),
let $\psi \colon D_{\{\cdot,\cdot\}}
\to
S \otimes_A L$
be the morphism of $S$-modules determined by
$\psi(da) = 0,\,
a \in A$, and
$\psi(d\alpha) = \alpha,\,\alpha \in L$,
and let
$\psi^* \colon
\roman{Alt}_S(S \otimes _A L,S)
\to
\roman{Alt}_S(D_{\{\cdot,\cdot\}},S)$
be the induced morphism of graded algebras.
Inspection shows its composite
$$
\roman{Alt}_A(L,A) 
\to
\roman{Alt}_S(D_{\{\cdot,\cdot\}},S)
\tag3.17.1
$$
with the obvious map
$
\roman{Alt}_A(L,A) 
\to
\roman{Alt}_A(L,S) 
\cong\roman{Alt}_S(S \otimes _A L,S)
$
is a morphism
of differential graded algebras.
We note that 
(1.16.1) and (1.16.2)
endow
$S \otimes _A L$
in fact
with a structure of an
$(R,S)$-Lie algebra, but $\psi$ is {\it not\/} a morphism
of
$(R,S)$-Lie algebras.
In particular, let
$\chi \colon L \otimes_A L \to A$
be a closed alternating 2-form and
let $\chi^{\flat}\colon D_S \otimes_S D_S \to S$
be its image in
$\roman{Alt}^2_S(D_{\{\cdot,\cdot\}},S)$.
Then
$\pi_{\{\cdot,\cdot\}}+ \chi^{\flat}$
is a closed alternating $S$-bilinear 2-form on $D_S$. Moreover,
let $\omega \colon D_S \to \Hom_S(D_S,S)$ 
be its adjoint, and
define 
a new bracket
$[\cdot,\cdot] \colon D_S \otimes_S D_S \to D_S$
by (3.11.0).
Inspection shows that this yields a structure of an
$(R,S)$-Lie algebra on $D_S$ and hence,
by virtue of (3.11), a new Poisson algebra structure 
$\{\cdot,\cdot\}$
on $S = S_A[L]$.
Explicitly,
this structure is determined by
$$
\{\alpha,\beta\} = [\alpha,\beta] + \chi(\alpha,\beta),
\quad \alpha,\beta \in L,
\tag3.17.2
$$
together with
(3.16.2.2) and (3.16.2.3).
This yields examples of Poisson algebras with non-trivial Poisson class.
\paragraph
To explain the geometric analogue of this
class of examples,
let $N$ be a smooth finite dimensional manifold,
let $\tau 
\colon T^*N \to N$
be its  cotangent bundle,
let
$\sigma =d \bold {p dq}$
be the standard symplectic form on $T^*N$,
and let $\chi$ an arbitrary closed 2-form on $N$.
Then 
$\sigma + \tau^*(\chi)$
is again a symplectic structure on 
$T^*N$.
The above 2-form
$\pi_{\{\cdot,\cdot\}}+ \chi^{\flat}$
is the formal analogue of the image
of
$\sigma + \tau^*(\chi)$
under the isomorphism
$\roman{Alt}_C(L_C,C)\to
\roman{Alt}_C(D_C^{\roman{geo}},C)$
in the complex 
$\roman{Alt}_C(D_C^{\roman{geo}},C)$
computing
the Poisson cohomology of
the algebra 
$C =C^{\infty}(T^*N)$
of smooth functions
on $T^*N$;
here $L_C$ refers to the smooth vector fields on $T^*N$,
and $\roman{Alt}_C(L_C,C)\to
\roman{Alt}_C(D_C^{\roman{geo}},C)$
to the isomorphism induced by the adjoint of
$\sigma + \tau^*(\chi)$, cf. (3.15).
In particular, the cohomology class
$[\sigma + \tau^*(\chi)] \in \Ho^2(T^*N)$
is non-zero if 
$[\chi] \in \Ho^2(N)$
is non-zero.
It is clear that
this yields examples of Poisson algebras with non-trivial Poisson class.
\paragraph
\noindent
{\smc Example 3.18.}
Another variant of (3.16) yields an example which 
in a special case
goes back to
{\smc Lie} [63]:
Let $g$ be 
a Lie algebra over $R$,
let $S =S[g]$
be the symmetric algebra
on $g$, and let 
$\{\cdot,\cdot\}
\colon
S\otimes S\to S$
be the corresponding Poisson structure
(3.16.1) (with the obvious change in notation).
We assume that $g$ is projective as an $R$-module.
Then the morphism
$$
S \otimes g \to D_S,
\quad s \otimes x \mapsto s dx,\ s \in S,\, x \in g,
\tag3.18.1
$$
is an isomorphism of $S$-modules.
Furthermore,
we assert that,
when $S \otimes g$
is endowed with the induced
$(R,S)$-Lie algebra structure
illustrated in
(1.16.4),
the morphism (3.18.1)
is an isomorphism
 $S \otimes g\to D_{\{\cdot,\cdot\}}$
of
$(R,S)$-Lie algebras,
where as before
$D_{\{\cdot,\cdot\}}$
refers to
$D_S$ 
with the $(R,S)$-Lie algebra structure
given by (3.8.1) above.
In fact, for $x,y \in g \subseteq S$, the Poisson bracket
$\{x,y\} \in S$ is defined by
$\{x,y\}= [x,y] \in g \subseteq S$.
Now, on the one hand,
the bracket (3.8.1) on $D_S$ is given by
$$
[adx,bdy] =
a\{x,b\}dy
+
b\{a,y\}dx
+
abd\{x,y\},
\quad a,b \in S,\, x,y \in g,
$$
while on the other hand the 
bracket which is part of the
induced $(R,S)$-Lie algebra structure
(1.16.4) on $S \otimes g$
is given by
$$
[ax,by] =
ax(b)y
-
by(a)x
+
ab[x,y],
\quad a,b \in S,\, x,y \in g,
$$
where we have discarded the
tensor product symbol and written $ax = a \otimes x$ etc.
Hence (3.18.1) is indeed an isomorphism
of $(R,S)$-Lie algebras.
In this way,
for an arbitrary Poisson algebra $A$,
the $(R,A)$-Lie algebra structure
(3.8.1)
on the $A$-module $D_A$
of formal differentials appears as a generalization
of 
the concept of
induced structure
in
(1.16.4).
\paragraph
Under the present circumstances we can identify 
Poisson cohomology with a well known object:
It is clear that
the obvious injection map
$g \to S \otimes g$
induces
an isomorphism
$$
(\roman{Alt}_S(S \otimes g,S),d)
\longrightarrow
(\roman{Alt}_R(g,S),d)
\tag3.18.2
$$
of chain complexes where $d$ refers to the
corresponding
{\smc Cartan-Chevalley-
\linebreak
Eilenberg} differentials (1.3).
In view of the isomorphism
(3.18.2) we conclude:
\proclaim {3.18.3}
For any 
$g$-module $M$,
the Poisson cohomology of
$S[g]$ with values in $M$
(with the obvious 
$(S[g],D_{\{\cdot,\cdot\}})$-module structure on $M$)
is isomorphic to the usual Lie algebra cohomology of $g$
with values in $M$.
\endproclaim
We note that
in view of what was said in (3.16) above,
the Poisson class of 
$(S,\{\cdot,\cdot\})$ in
$\Ho^2_{\roman{Poisson}}(S,\{\cdot,\cdot\};S)$
is trivial.
In fact, in view of (3.16.3) a Poisson potential 
$\vartheta \colon D_S \to S$ is given by
$\vartheta (dx) = x,\, x \in g$.
\paragraph   
An algebra closely related to
$S[g]$
is classical and has been studied extensively in the literature:
Assume the ground ring is that of the reals
$\bold R$, let $g$ be a finite dimensional
real Lie algebra,
and let $g^*$ be its dual.
Then the 
obvious modification of the above construction
yields a 
Poisson algebra structure on the algebra of smooth,
analytic, or some other class
 functions
on
$g^*$ (the prior construction yields the polynomial functions).
In this way $g^*$ becomes a Poisson manifold 
(in the appropriate category),
even though it is not in general a symplectic manifold.
Such structures were studied by {\smc Lie} [63]
and others
{\smc Berezin} [6],
{\smc Kostant} [50],
{\smc Kirillov} [46],
{\smc Souriau} [93], [94],
{\smc Lichnerowicz} [61],
{\smc Weinstein} [102] -- [106],
{\smc Koszul} [52],
{\smc Conn} [18], [19].
Furthermore,
the obvious modifications of
(3.18.1) and (3.18.2)
yield at once the following:
\proclaim {3.18.4}
The Poisson cohomology of
$C^{\infty}(g^*)$
with values in   $C^{\infty}(g^*)$
is isomorphic to the usual Lie algebra cohomology of $g$
with values in $C^{\infty}(g^*)$
with respect to the coadjoint representation.
\endproclaim
This cohomology plays a significant role in the
\lq\lq linearization problem\rq\rq \ 
for real Poisson manifolds
{\smc Weinstein} [102], [103],
{\smc Conn} [18], [19].
\paragraph
\noindent
{\smc Example 3.19.}
Let $A = R[u_1,u_2]$ be the polynomial algebra
in $u_1$ and $u_2$ as indicated, let $p \in A$ be an arbitrary
polynomial, and define a Poisson structure 
$\{\cdot,\cdot\}$ on $A$ by
$\{u_1,u_2\} = p$.
Then a little thought reveals that
if this Poisson structure admits a Poisson potential
$\vartheta$, for degree reasons we must have
$\vartheta(du_1) = \alpha u_1$ and
$\vartheta(du_2) = \beta u_2$ 
for suitable constants $\alpha, \beta \in R$.
Now
$$
\aligned
d\vartheta(du_1,du_2) &=
\{u_1,\vartheta(du_2)\}
+
\{\vartheta(du_1),u_2\}
-
\vartheta(d  \{u_1,u_2\} )
\\
&=
(\alpha + \beta) p 
- \alpha \frac{\partial p} {\partial u_1} u_1
- \beta \frac{\partial p} {\partial u_2} u_2,
\endaligned
\tag3.19.1
$$
and from this it is straightforward to decide whether
$(A,\{\cdot,\cdot\})$ admits a Poisson potential.
For example, when 
$p$ is a non-trivial homogeneous quadratic polynomial
$p = a u_1^2 + b u_1 u_2 + cu_2^2$,
(3.19.1)
entails
$$
d\vartheta(du_1,du_2) = 
(\alpha - \beta)(cu_2^2 - a u_1^2),
$$
and hence the Poisson class
$[\pi_{\{\cdot,\cdot\}}] \in 
\Ho^2_{\roman{Poisson}}(A,\{\cdot,\cdot\};A)$
will be non-trivial
unless 
\roster
\item
$b=0,\,a =0,\,\alpha -\beta= 1, $ or
\item
$b=0,\,c =0,\,\beta-\alpha=1$.
\endroster
It was pointed out to me by A. Weinstein that
the question whether or not 
$(A,\{\cdot,\cdot\})$ admits a Poisson potential
is closely related with the question whether or not
the vector field $X = (u_1,u_2)$
generating dilations of the plane
leaves invariant the Poisson structure, i.~e.
whether or not
the Lie derivative
$\lambda_X(p)$ equals $p$.

\bigskip
\noindent
{\bf 4. Linear representations of the underlying Lie algebra}
\medskip
\noindent
Let $(A,\{\cdot,\cdot\})$ be Poisson algebra,
and let $D_{\{\cdot,\cdot\}}$
be the corresponding $(R,A)$-Lie algebra
introduced in (3.8).
In this Section we 
relate the Poisson class (3.10.1)
of $(A,\{\cdot,\cdot\})$
with the representation theory
of the 
underlying Lie algebra. 
This is motivated by the Dirac quantization problem
{\smc Dirac} [22], [23],
see also 
{\smc Kostant} [48], 
{\smc Kostant-Sternberg} [51], 
{\smc Simms-Woodhouse} [87], {\smc Woodhouse} [108].
The problem is to construct an $R$-linear representation
so that {\it the elements of the ground ring $R$ act by scalar 
multiplication\/};
this is non-trivial since under the adjoint representation
the elements of the ground ring act trivially.
\paragraph
Let
$\pi =\pi_{\{\cdot,\cdot\}}$ be the Poisson 2-form 
(3.5.1)
of
$(A,\{\cdot,\cdot\})$, and let
$\bar L_{\{\cdot,\cdot\}}= A \oplus_{-\pi} D_{\{\cdot,\cdot\}}$
be the indicated $(R,A)$-Lie algebra whose structure
is given by (3.10.2) except that 
$-\pi_{\{\cdot,\cdot\}}$ comes into play instead of
$\pi_{\{\cdot,\cdot\}}$, so that the corresponding extension
$$
0 
@>>> 
A
@>>> 
\bar L_{\{\cdot,\cdot\}}
@>>> 
D_{\{\cdot,\cdot\}}
@>>> 
0
\tag4.1
$$
of $(R,A)$-Lie algebras is classified by
$-[\pi_{\{\cdot,\cdot\}}] 
\in \Ho^2_{\roman{Poisson}}(A,{\{\cdot,\cdot\}};A)$,
cf. (2.6).
Moreover, let
$$
\iota_{\{\cdot,\cdot\}} 
=
(\roman{Id},d)
\colon A \longrightarrow
\bar L_{\{\cdot,\cdot\}}.
\tag4.2.1
$$
It is clear that
$dA \subseteq D_A$
inherits a Lie algebra structure from that
on $A$.
Moreover, in view of (3.8.1), for $u,v \in A$,
$$
d\{u,v\} = [du,dv],
$$
whence $dA \subseteq D_{\{\cdot,\cdot\}}$
is a sub Lie algebra over $R$;
we denote $dA$
together with this Lie algebra structure
by $\roman{Ham}_{\{\cdot,\cdot\}}$.
We note that
the corresponding morphism
(3.6)
maps the elements of
$\roman{Ham}_{\{\cdot,\cdot\}}$
to the Hamiltonian elements in $\roman{Der}(A)$,
cf. what was said just before (3.8).
Direct inspection proves the following:
\proclaim{Proposition 4.2}
Let $(A,\{\cdot,\cdot\})$ 
be a Poisson algebra.
Then the morphism {\rm (4.2.1)}
is one of Lie algebras over $R$ and makes commutative
the diagram
$$
\CD
0
@>>>
R
@>>>
A
@>>>
\roman{Ham}_{\{\cdot,\cdot\}}
@>>>
0
\\
@.
@VVV
@V{\iota_{\{\cdot,\cdot\}}}VV
@VVV
@.
\\
0 
@>>> 
A
@>>> 
\bar L_{\{\cdot,\cdot\}}
@>>> 
D_{\{\cdot,\cdot\}}
@>>> 
0
\endCD
\tag4.2.2
$$
in the category of 
 Lie algebras over $R$.
Furthermore, 
$\iota = \iota_{\{\cdot,\cdot\}}$ satisfies the formula
$$
\iota (ab) =
a\iota (b) +
b\iota (a) -(ab,0),
\tag4.2.3
$$
where $\iota = \iota_{\{\cdot,\cdot\}}$.
\endproclaim
\paragraph
We shall say that
$(A,\{\cdot,\cdot\})$
is {\it representable\/}
if
there is a projective
rank one   $A$-module
$M$
and a
$(D_{\{\cdot,\cdot\}})$-connection
$
\nabla
\colon
D_{\{\cdot,\cdot\}}
\longrightarrow
\roman{End}(M)
$
with curvature 
$-\pi_{\{\cdot,\cdot\}}$.
In view of (2.15),
$(A,\{\cdot,\cdot\})$
is representable
if and only if 
its Poisson class
$
[\pi_{\{\cdot,\cdot\}}] \in \Ho^2_{\roman{Poisson}}(A,{\{\cdot,\cdot\}};A)
$
(cf. (3.10))
lies in the image of the corresponding morphism
(2.15.1).
We note that $\roman{End}_A(M) \cong A$ canonically.

\proclaim{Theorem 4.3}
Let
$(A,\{\cdot,\cdot\})$ be a
representable Poisson algebra,
and let
$M$ be a projective
rank one   $A$-module
together with
a 
$(D_{\{\cdot,\cdot\}})$-connection
$
\nabla
\colon
D_{\{\cdot,\cdot\}}
\longrightarrow
\roman{End}(M)
$
with curvature 
$-\pi_{\{\cdot,\cdot\}}$.
Then the composite of
$\iota_{\{\cdot,\cdot\}}
\colon
A
\to
\bar L_{\{\cdot,\cdot\}}
$
with the flat 
\linebreak
$(\bar L_{\{\cdot,\cdot\}})$-connection
$
\bar L_{\{\cdot,\cdot\}} 
\longrightarrow
\roman{End}(M)
$
on $M$ is a representation of
the $R$-Lie algebra underlying
$A$ on $M$ viewed as an $R$-module
having the property that
the \lq\lq constants\rq\rq, i. e. the elements of $R$,
act by multiplication.
Explicitly, the representation is given by
$a \mapsto \hat a$ where, for $a \in  A,$
$\hat a$ refers to the operator given by
$$
\hat a x =  \nabla_{da} x + ax,
\quad 
x \in M.
$$
When the form
{\rm (3.5.1)}
is non-degenerate,
the representation is faithful.
\endproclaim
\demo{Proof}
This follows at once from (4.2). \qed
\enddemo
In particular, let 
$(A,\{\cdot,\cdot\})$ be a Poisson algebra with zero Poisson class
$[\pi_{\{\cdot,\cdot\}}]$,
and let $\vartheta \colon D_{\{\cdot,\cdot\}} \to A$ be
a Poisson potential 
for $(A,\{\cdot,\cdot\})$
so that
$d\vartheta = \pi_{\{\cdot,\cdot\}}$.
Define a 
$D_{\{\cdot,\cdot\}}$-connection 
$\nabla \colon D_{\{\cdot,\cdot\}} \to \roman{End}(A)$
on $A$, viewed as a free $A$-module with basis 1,
by
$$
\nabla_{\alpha}(a) = \alpha^{\sharp}(a) - a \vartheta(\alpha).
\tag4.3.1
$$
Then $\nabla$ has curvature $-\pi_{\{\cdot,\cdot\}}$,
whence
$(A,\{\cdot,\cdot\})$ is representable.
Notice this applies in particular to the 
example
in (3.18).
\paragraph
We now take
the ground ring $R$ to be  that of the reals $\bold R$ 
and introduce a concept
of \lq\lq prequantization\rq\rq\ for
real Poisson algebras
by means of a variant of the above.
The problem of quantizing Poisson algebras
which are not associated with a symplectic manifold
arises in the physics of singular constrained systems,
see e. g.
{\smc Gotay} [31], {\smc \'Sniatycki-Weinstein} [92],
and Section 5 below.
\paragraph
We write $\bold C$ for the complex numbers.
Consider the complexified algebra $A\otimes \bold C$, with the
obvious Poisson structure, still denoted by
$\{\cdot,\cdot\}$, and let
$D_{\{\cdot,\cdot\}}^{\bold C}$
and  
$\pi_{\{\cdot,\cdot\}}^{\bold C}$
be the 
corresponding 
$(\bold C,A\otimes \bold C)$-
Lie algebra
and Poisson 2-form, respectively. 
We note that these objects arise from the former ones
by merely an extension of \lq\lq scalars\rq\rq.
We shall say that the {\it real\/} Poisson algebra
$(A,{\{\cdot,\cdot\}})$
is {\it quantizable\/}
if 
there is a 
projective rank one  
$(A\otimes \bold C)$-module
$M$
with a 
$(D_{\{\cdot,\cdot\}}^{\bold C})$-connection
$\nabla$
having curvature
$$
K_{\nabla}=
-
i\, \pi_{\{\cdot,\cdot\}}^{\bold C}\in 
\roman{Alt}^2_{A\otimes \bold C}
(D_{\{\cdot,\cdot\}}^{\bold C},A\otimes \bold C).
$$
In view of (3.12.13) and (3.15),
when $(A,{\{\cdot,\cdot\}})$
is the algebra of smooth real valued functions on a
finite dimensional symplectic manifold $(N,\sigma)$,
then $(A,{\{\cdot,\cdot\}})$ is quantizable
in our sense
if and only if there is
a 
{\it prequantum bundle\/}
for $(N,\sigma)$, i.~e.
a complex line bundle
$\lambda \colon E \to N$
with a 
connection
having  curvature 
$-
i\, \sigma$. 
The reader will note that 
some of the numerical constants used here 
differ from those in the literature but this is of course of no account.
\paragraph
Let $M$ be
an $(A \otimes \bold C)$-module.
An $A$-linear pairing
$
(\cdot,\cdot) \colon M \otimes_A  M
\longrightarrow A\otimes \bold C
$
on  $M$
will be said to be
a {\it Hermitian structure\/}
on $M$ if it
has the usual Hermitian properties.
Further,
we shall say 
that a $(D_{\{\cdot,\cdot\}}^{\bold C})$-connection $\nabla$ on $M$
{\it preserves\/}
the  Hermitian structure $(\cdot,\cdot)$
if for $\alpha \in D_{\{\cdot,\cdot\}}^{\bold C}$ and $s,s' \in M$,
$
\alpha(s,s') = (\nabla_\alpha s,s') + (s,\nabla_{\bar \alpha}s')
$
where as usual the symbol \lq$\ \bar {}\ $\rq\  refers to complex conjugation.
Extending common notation, we shall say that
an $\bold R$-linear operator 
$\vartheta \colon M \to M$
is {\it symmetric\/}, if
for $s,s' \in M$,
$
(\vartheta s,s') = (s,\vartheta s').
$
The symmetric operators on $M$ constitute 
a Lie algebra over $\bold R$
with Lie bracket $[\cdot,\cdot]_s$ given by
$
[\vartheta,\xi]_s = 
i\, [\vartheta,\xi],
$
where
$\vartheta$ and $\xi$ are arbitrary symmetric operators and
where $[\vartheta,\xi] = \vartheta \xi - \xi \vartheta$
as usual.
\proclaim{Theorem 4.4}
Let $(A,{\{\cdot,\cdot\}})$
be a quantizable Poisson algebra,
 let $M$ be a
projective rank one  
$(A\otimes \bold C)$-module
with a 
$(D_{\{\cdot,\cdot\}}^{\bold C})$-connection
$\nabla$
having curvature
$
K_{\nabla}=
-
i\, \pi_{\{\cdot,\cdot\}}^{\bold C},
$
and, for $\alpha \in D_{\{\cdot,\cdot\}}$ 
and
$\,a \in A$, let
$$
\omega(\alpha) 
= -i\,\nabla_{\alpha},
\quad
\omega(a) 
= \mu_a,
\tag4.4.1
$$
where $\mu_a$ refers to multiplication by $a \in A $.
This defines
a representation 
\linebreak
$
\omega
\colon
\bar L_{\{\cdot,\cdot\}}
\longrightarrow
\roman{End}_{\bold R}(M)
$
of
the $(\bold R,A)$-Lie algebra $\bar L_{\{\cdot,\cdot\}}$
by 
$\bold R$-linear
operators on $M$ so that,
for $a \in A$ and $u,v \in \bar L_{\{\cdot,\cdot\}}$,
$$
\omega(a) 
= \mu_a,\quad 
i[\omega(u),\omega(v)] 
= 
\omega([u,v]),\quad 
\omega(a u) 
= (\mu_a\omega)(u),
$$
where $\mu_a\omega$ refers to composition of operators.
When 
$M$ has a Hermitian structure in such a way that
$\nabla$ is compatible with the Hermitian structure,
the representation is 
a homomorphism into the real Lie algebra
of symmetric operators on $M$.
\endproclaim
\demo{Proof} This is left to the reader. \qed
\enddemo
We note that, in contrast to (4.3),
the morphism $\omega$ in (4.4) does {\it not\/}
endow
$M$ with a structure of a
$(A,\bar L_{\{\cdot,\cdot\}})$-module.
We also note that, as usual, the factor $i$
in (4.4.1)
has the effect that,
when
$M$ has a Hermitian structure and
$\nabla$ is compatible with the Hermitian structure,
 a differential operator is represented by
a symmetric operator rather than an antisymmetric one.
\paragraph
When we write out
the composite of
$\iota_{\{\cdot,\cdot\}}
\colon
A
\longrightarrow
\bar L_{\{\cdot,\cdot\}}
$
with $\nabla$,
we arrive at a proof of the following:
\proclaim{Corollary 4.5}
Let
$(A,\{\cdot,\cdot\})$
be a quantizable Poisson algebra,
and let
$M$ be a 
projective rank one  
$(A\otimes \bold C)$-module
with a 
$(D_{\{\cdot,\cdot\}}^{\bold C})$-connection
$\nabla$
having curvature
$
K_{\nabla}=
-
i\, \pi_{\{\cdot,\cdot\}}^{\bold C}.
$
Then, for $a \in A$, the formula
$$
\hat a s = -
i\, \nabla_{da} s + as,
\ s \in M,
\tag4.5.1
$$
yields a representation of the $\bold R$-Lie algebra underlying
the Poisson algebra $A$
by 
$\bold R$-linear
operators on $M$ so that
\roster
\item
the constants act by
multiplication,
and
\item
for $a,b \in A$,
$
\widehat {\{a,b\}} = 
i\,[\hat a,\hat b].
$
\endroster
When 
$M$ has a Hermitian structure in such a way that
$\nabla$ is compatible with the Hermitian structure,
the representation is 
a homomorphism into the real Lie algebra
of symmetric operators on $M$.
\endproclaim

In {\smc Dirac's} terms [22], [23],
the Poisson bracket  $\{\cdot,\cdot\}$
is thus the classical counterpart of
the quantum commutator 
$
i\,[\cdot,\cdot]$.
When $(A,{\{\cdot,\cdot\}})$
is the algebra of smooth real valued functions on a
finite dimensional symplectic manifold $(N,\sigma)$,
up to constants, the Corollary 
yields precisely 
{\smc Kostant's} prequantization construction [48],
see also
{\smc Simms-Woodhouse} [87],
{\smc Woodhouse} [108].
In this case, the significance of the corresponding
Atiyah sequence analogous to (4.1) (cf. 2.16.1)
for quantization
was noticed already by {\smc Almeida-Molino} [2].  
Corollary 4.5 goes beyond traditional prequantization since it
yields a prequantization for
e.~g. a Poisson algebra of the kind in (3.18)
(involving a Lie algebra $g$ over the ground ring and the symmetric algebra
$S[g]$ on $g$).
\paragraph
The usual completion of
the geometric quantization procedure
as explained e.~g. in {\smc Woodhouse} [108]
points the way to the completion of 
the present kind of quantization.
Here we confine ourselves with some remarks:
The usual notions of an isotropic 
and Lagrangian
subspace
make perfectly good sense for
$D_{\{\cdot,\cdot\}}$ and the 2-form
$\pi_{\{\cdot,\cdot\}}$, and hence 
so does 
the concept of a {\it polarization\/}.
Indeed, given a real Poisson algebra
$(A,\{\cdot,\cdot\})$,
a {\it complex  polarization\/} $P$ may be defined as a
$(\bold C,A\otimes \bold C)$-Lie sub algebra 
$P \subseteq D_{\{\cdot,\cdot\}}^{\bold C}$
which is maximally isotropic
with respect to
the complexified 2-form
$\pi^{\bold C}_{\{\cdot,\cdot\}}$.
However, since
$D_{\{\cdot,\cdot\}}^{\bold C}$
will in general {\it not\/}
act faithfully on $A$,
it remains to be seen what the appropriate generalization of
a polarization should be, cf. what will be said in the next Section.
Whatever choice of polarization then has been made,
one can then introduce analogues of the
cohomology spaces introduced e.~g. on p. 216 
of {\smc Woodhouse} [108] and the usual inner product problem
 arises,
cf. {\smc Rawnsley-Schmid-Wolf} [79].
In the next Section
we shall 
employ the above
to quantize
a system described in terms of
a Poisson algebra
that is {\it not\/} associated with a
symplectic manifold.

\beginsection 5. Poisson reduction and a non-classical example

In this Section we extend a standard construction
to the present setting.
This will yield examples of
Poisson algebras of a kind different from 
those in (3.16) -- (3.19).
\paragraph
Let $(A,\{\cdot,\cdot\})$ be a Poisson algebra over $R$,
let $g$ be a Lie algebra over $R$, 
 and let
$\delta \colon g \to A$
be a morphism of Lie algebras (over $R$).
Let $I \subseteq A$ be the ideal in $A$ generated by 
the image $\delta(g) \subseteq A$.
Since $\delta$ is compatible with the Lie structures,
$I$ is closed under the Poisson bracket;
note, however, that $I$ is not necessarily a Lie ideal.
Further, with the obvious structure
explained in (1.18), $A \otimes g$ is an
$(R,A)$-Lie algebra, 
$A$ is an $(A,A \otimes g)$-module,
and the obvious map
$\delta^{\sharp}
\colon A \otimes g \to A$ given by 
$\delta^{\sharp}(a \otimes y) = a \delta(y) \subseteq A$
is a morphism 
of $(A,A \otimes g)$-modules whence
the quotient $A/I$ inherits a $g$-action.
Let $A_{\roman{red}} = (A/I)^g$
be the sub algebra of $g$-invariants.
To describe it we recall
that for an arbitrary Lie algebra $k$ and a sub Lie algebra 
$h \subseteq k$, the {\it normalizer\/} 
$h^k \subseteq k$ 
of $h$ in $k$
consists of all
$\alpha \in k$ having the property that
$[\alpha,x] \in h$ for every $x \in h$.
Now
$A_{\roman{red}} = (A/I)^g$
it consists of all classes of elements $a \in A$ for which
$\{a,I\} \subseteq I$,
i.~e.
$A_{\roman{red}}$ is the image $I^A/I$ in $A/I$ of the normalizer 
$I^A\subseteq A$ of $I$ in $A$
in the sense of Lie algebras.
Inspection shows that $I^A$ inherits a Poisson algebra structure, and hence so 
does
$A_{\roman{red}} = I^A/I$.
We write
$\{\cdot,\cdot\}_{\roman{red}}$
for the latter structure
and 
call
$(A_{\roman{red}},\{\cdot,\cdot\}_{\roman{red}})$  
the {\it reduced\/} Poisson algebra
of $(A,\{\cdot,\cdot\})$ (with respect to $\delta$).
Special cases of this 
kind of reduction are
the reduction procedures of
{\smc Marsden and Weinstein} [72], and of
{\smc \'Sniatycki and Weinstein} [92], see
also {\smc Weinstein} [103],
p. 51 of {\smc Kostant and Sternberg} [51],
and {\smc Stasheff} [95], [96].
In fact,
let $A$ be the ring of smooth functions on a 
Poisson manifold $N$, let
$\delta \colon g \to A$
be a morphism of real Lie algebras,
let $I \subseteq A$ be the ideal in $A$ generated by 
$\delta(g) \subseteq A$,
and define the corresponding {\it moment mapping\/}
$J \colon N \to g^*$
as usual by
$$
(J(x))(\xi) = (\delta(\xi))(x),\quad x \in N,\, \xi \in g.
$$
When $0$ is a regular value of $J$,
$J^{-1}(0)$ is a smooth manifold,
$I$ coincides with
the ideal of smooth functions 
on $N$ that vanish on $J^{-1}(0)$,
 and
$A/I$ 
is canonically isomorphic to the algebra of smooth functions
on $J^{-1}(0)$;
moreover, the $g$-action on $A/I$ then induces a foliation
of $J^{-1}(0)$.
Furthermore,
the normalizer $I^A$ consists of the smooth functions $f$ on $N$
having the following property:
\proclaim {5.1}
For $X \in g$,
$\{\delta(X),f \}$
vanishes on $J^{-1}(0)$.
\endproclaim
\noindent
Geometrically this means that $f$ is constant
along the 
restriction to
$J^{-1}(0)$ of
any integral curve in $N$
of the vector field 
$\{\delta(X),-\}$.
Hence
the reduced Poisson algebra $A_{\roman{red}} = (A/I)^g$
then appears as the algebra of 
classes of
smooth functions
on $N$ having the property (5.1), where two such functions are identified
if they 
assume the same values
on $J^{-1}(0)$.
In the situation of symplectic reduction
the assumptions are made that (i)
$N$ is a finite dimensional symplectic manifold
with a Hamiltonian $G$-action, where $G$ is a Lie group
with Lie algebra $g$,
that (ii) the map $J$ is the corresponding moment mapping,
that (iii) 
$0$ is a regular value thereof, 
and that (iv) the foliation on $J^{-1}(0)$ 
comes from a principal $G$-bundle so that the space of leaves is the base $B$
of the bundle. Then $B$
inherits a symplectic structure
$\sigma_B$, the algebra of smooth functions on $B$
may be identified with 
$A_{\roman{red}}$,
and the reduced Poisson structure
$\{\cdot,\cdot\}_{\roman{red}}$ 
is induced from $\sigma_B$.
When $0$ is not a regular value of $J$,
the {\smc \'Sniatycki-Weinstein}-reduction  [92] 
is formally the same as the above algebraic procedure
except that in the description in [92] 
still
a Lie group $G$ with Lie algebra $g$
comes into play.
However, the question whether
the ideal $I$ 
contains all smooth functions on $N$ that vanish
on $J^{-1}(0)$
then becomes a delicate problem.
When this is so,
the algebra $A/I$ 
may be identified with the algebra of smooth functions
in the sense of {\smc Whitney} on $J^{-1}(0)$,
see e.~g. {\smc Malgrange}~[68],
and
the normalizer $I^A$ still consists of the smooth functions $f$ on $N$
having the property (5.1).
\paragraph
We now illustrate the
algebraic reduction procedure
 with an example which also occurs in  {\smc Gotay} [31].
This example has quadratic singularities.
We note that {\smc Arms, Marsden and Moncrief} [3]
have shown that singular momentum mappings
typically have quadratic singularities.
\paragraph
Let $Q$ be four dimensional
Minkowski space-time,
let $T^*Q$ be its cotangent bundle,
with the usual coordinates $(x_0,x_1,x_2,x_3,p_0,p_1,p_2,p_3)$
and symplectic form
$\sigma = \sum dp_j \wedge dx_j$,
let $A$ be its ring of smooth functions,
and let
$J \colon T^*Q \to \bold R$ be the moment mapping
$J(x,p) = p_0^2 - p_1^2 -p_2^2 - p_3^2$.
Then $J = m^2$ 
describes the 
{\it constrained system\/}, see e.~g. p. 256 of
{\smc Woodhouse} [108],
for a 
spinless relativistic particle
with rest mass $m$.
Henceforth we write $C(c)$ 
for the zero locus of
$J-c=0$.
Notice that, geometrically, for $c>0$, 
$C(c)$  is a product of an $\bold R^4$ with a
2-sheeted hyperboloid 
$H(c)$ 
(each copy of which is topologically a
cone on
$S^2$),
 for
$c<0$, $C(c)$ is a 
product of an $\bold R^4$ with a 1-sheeted hyperboloid $H(c)$
(which is topologically an $S^2 \times \bold R$),
and for $c=0$ these degenerate
to a product of an $\bold R^4$ with a cone $C$
(= the zero locus of $p_0^2 - p_1^2 -p_2^2 - p_3^2= 0$).
\proclaim{Lemma 5.2. [Gotay]}
The ideal 
$I$ in $A$ of
smooth functions
that vanish on $C(0)$
coincides with
the ideal in $A$ generated by $J$.
\endproclaim
\demo{Proof} See
(4.1) of {\smc Gotay} [31]. \qed 
\enddemo
Hence
the algebra $A/I$ 
may be identified with the algebra of smooth functions
in the sense of {\smc Whitney} on $J^{-1}(0)$.
We reproduced the statement of the Lemma
since it will enable us to
illustrate geometrically the quantization result
obtained below.
The quantization procedure itself will {\it not\/}
make use of the Lemma.
\paragraph
We view
$\bold R$ as an abelian Lie algebra,
with the single basis element $\bold 1$.
By construction, for every $c \in \bold R$,
the adjoint $\delta$ of the corresponding moment mapping is given by
$\delta(\bold 1) = J - c \in A$,
and, whatever value $c$ assumes,
the vector field $\{\delta (\bold 1), -\}= \{J, -\}$ 
may be described by
$$
\{J, -\} =
-2p_0 \frac {\partial}{\partial x_0}
+2p_1 \frac {\partial}{\partial x_1}
+2p_2 \frac {\partial}{\partial x_2}
+2p_3 \frac {\partial}{\partial x_3}.
\tag5.2.1
$$
In view of what was said 
above, 
whatever value $c$ assumes,
the reduced Poisson algebra
$A_{\roman{red}}$
consists of classes of smooth functions
$f$ on $T^*Q$ 
having the property that
$\{J,f\}$
vanishes on $C(c)$,
where two functions are identified whenever they coincide
on $C(c)$.
These are precisely classes of functions that
are constant along the flow lines
$$
(t,(x,p)) \mapsto (x_0+tp_0,x_1-tp_1,x_2-tp_2,x_3-tp_3,p_0,p_1,p_2,p_3),
\quad t \in \bold R,
\tag5.2.2
$$
with $p_0^2 -p_1^2 -p_2^2 -p_3^2 = c$.
In {\smc Dirac's} terms [23]
the elements of
$A_{\roman{red}}$
are the corresponding true (classical)
observables. 
\paragraph
The rule (5.2.2) defines  a smooth action 
$$
G \times T^*Q \longrightarrow T^*Q
\tag5.2.3
$$
of the
additive group $G$ of the real numbers on
$T^*Q$, which is free
on
$T^*Q-\bold R^4 \times \{0\}$
and fixes the singular set
$\bold R^4 \times \{0\}$
pointwise.
Since for $c \ne 0$,
$(p_0,p_1,p_2,p_3) = (0,0,0,0) \not \in C(c)$,
(5.2.3) yields a structure of a principal
$G$-bundle on $C(c)$;
hence,
by symplectic reduction
{\smc Marsden-Weinstein} [72],
the quotient 
$B(c) =C(c)/G$
then inherits
a symplectic structure.
Likewise,
when $c = 0$,
(5.2.3) yields a structure of a principal
$G$-bundle on 
$C_0 = C(0)-\bold R^4 \times \{0\}$ and
again by symplectic reduction
the quotient 
$B_0=C_0/G$
inherits
a symplectic structure.
In both cases, the reduced Poisson algebra
in our sense then manifestly coincides with the standard Poisson algebra of 
smooth functions
on $B(c)$ or $B_0$ (as appropriate).
However, such a statement makes no sense for $C(0)$.
\paragraph
The problem we wish to study next is the
quantization of the reduced Poisson algebra
$(A_{\roman{red}},\{\cdot,\cdot\}_{\roman{red}})$.
Since for $c \ne 0$ 
the quotient $B(c)$ inherits a smooth symplectic structure this case
can be handled by standard geometric
quantization theory, and
 we assume henceforth that $c=0$.
In {\smc Gotay} [31] 
the corresponding particle is referred to
 as a
{\it photon\/}
and
treated by a method different from that to be given below.
We shall not use the name photon
since in the physics
literature a photon
means something else
(i.~e. an irreducible representation of the Poincar\'e group
with mass zero and helicity $\pm 1$).
\paragraph
We denote 
by $D_{\roman{red}}$
the
$A_{\roman{red}}$-module
of K\"ahler differentials
for $A_{\roman{red}}$ 
with the corresponding
$(\bold R,A_{\roman{red}})$-Lie algebra
structure
spelled out in (3.8), and we write
$\pi_{\roman{red}} \colon D_{\roman{red}} 
\otimes_{A_{\roman{red}}} D_{\roman{red}}
\to 
A_{\roman{red}}      
$
for its Poisson 2-form.
Since the Poisson structure on
$T^*Q$
is that induced from the cotangent bundle symplectic structure,
the morphism
$\vartheta^{\roman{geo}}
\colon D_{\{\cdot,\cdot\}}^{\roman{geo}} \to A$
of $A$-modules
given by
$$
\vartheta^{\roman{geo}}(dp_j) = p_j,\ 
\vartheta^{\roman{geo}}(dx_j) = 0, \quad 0 \leq j \leq 3,
\tag5.3
$$
is a geometric Poisson potential for the geometric Poisson 2-form
$\pi_{\{\cdot,\cdot\}}^{\roman{geo}}$ 
associated with the real Poisson algebra $A$ of real smooth functions
on $T^*Q$.
We note that, 
with $M = T^*Q$,
under the isomorphism
$T^*M \to TM$
of vector bundles
induced by the symplectic structure,
this geometric Poisson potential
corresponds to the usual symplectic potential
$\Theta = \sum p_i dx_i$,
viewed as a 1-form on $M$. 
In fact 
the isomorphism
$T^*M \to TM$
identifies $dp_i$ with 
the vector field $\frac{\partial}{\partial x_i}$ etc.;
cf. what was said in (3.16).
\proclaim{Theorem 5.4}
The geometric Poisson potential
$\vartheta^{\roman{geo}}$
induces an algebraic
Poisson potential $ \vartheta_{\roman{red}}
\colon  D_{\roman{red}} \to  A_{\roman{red}}$
for
the Poisson algebra $ A_{\roman{red}}$,
i.~e. a 1-form
$ \vartheta_{\roman{red}}$ so that
$$
d (\vartheta_{\roman{red}}) = 
 \pi_{\roman{red}} 
\in
\roman{Alt}^2_{A_{\roman{red}}}
( D_{\roman{red}}, A_{\roman{red}}).
$$
\endproclaim
\demo{Proof}
Let
$\pi_{\{\cdot,\cdot\}}
\colon 
D_{\{\cdot,\cdot\}} \otimes _A
D_{\{\cdot,\cdot\}} \to A
$ 
be the corresponding algebraic Poisson 2-form
for $(A,{\{\cdot,\cdot\}})$.
It is the composite of
the obvious surjection
$$
D_{\{\cdot,\cdot\}} \otimes _A
D_{\{\cdot,\cdot\}} 
\longrightarrow
D_{\{\cdot,\cdot\}}^{\roman{geo}} \otimes _A
D_{\{\cdot,\cdot\}}^{\roman{geo}}
$$
with the geometric Poisson 2-form
$\pi_{\{\cdot,\cdot\}}^{\roman{geo}}$.
Hence the composite
of
the obvious surjection
$$
D_{\{\cdot,\cdot\}} 
\longrightarrow
D_{\{\cdot,\cdot\}}^{\roman{geo}} 
$$
with
the geometric Poisson potential
$\vartheta^{\roman{geo}}$
is an algebraic Poisson potential
$$
\vartheta
\colon D_{\{\cdot,\cdot\}} \to A
$$
for
$(A,{\{\cdot,\cdot\}})$.
We assert that $\vartheta$ passes to a Poisson potential for $I^A$.
To see this we observe first that,
by construction, 
for
$f \in A$,
$$
\vartheta(df)
=
\left(
 \frac{\partial f}{\partial p_0} p_0
+ \frac{\partial f}{\partial p_1} p_1
+ \frac{\partial f}{\partial p_2} p_2
+ \frac{\partial f}{\partial p_3} p_3 \right) \in A.
\tag5.4.1
$$
We note, for clarity, that here $df$ refers to the formal differential
$df \in D_A$.
From this and the description
(5.2.1) of the vector field
$\{J, -\}$
we deduce that, for
$f \in A$,
$$
\{J, \vartheta (df)\} = \vartheta d\{J, f\}-\{J, f\}.
\tag5.4.2
$$
In fact, a straightforward calculation shows that, for $0 \leq j \leq 3$,
$$
\left\{J, p_j \frac{\partial f}{\partial p_j} \right\} = 
p_j \frac{\partial}{\partial p_j}\{J, f\}  - 
p_j \left\{\frac{\partial J}{\partial p_j} , f\right\}.
$$
From the description (5.4.1) of
$\vartheta (df)$ we conclude at once that
(5.4.2) holds.
We now recall that the normalizer
$I^A$ of $I$ in $A$ in the sense of Lie algebras
consists of all $f \in A$ 
having the property that
$\{J, f\}$  lies in $I$.
However, this implies that,
for $f \in I^A$,
$\{J, \vartheta (df)\}$
lies in $I$ as well,
i.~e.
that 
for $f \in I^A$
the value
$\vartheta (df)$ lies in $I^A$.
In fact, given
$f \in I^A$, there is an $h \in A$ so that
$$
\{J, f\} = h\,J.
$$
However, 
$$
\aligned
\vartheta (dJ) =\vartheta (d(p_0^2 - p_1^2- p_2^2- p_3^2 )) &= 
  p_0 \vartheta (d p_0)
- p_1 \vartheta (d p_1)
- p_2 \vartheta (d p_2)
- p_3 \vartheta (d p_3)
\\
&=
2(p_0^2 - p_1^2- p_2^2- p_3^2 ) = 2J \in  A.
\endaligned
$$
Moreover
$$
d(h\,J) = (dh) J + h\,dJ.
$$
Consequently
$$
\vartheta (d\{J, f\}) = 
(\vartheta (dh))\,J + 2h\,J \in I,
$$
whence, in view of (5.4.2),
$\{J, \vartheta (df)\} \in I$ and hence 
$\vartheta (df) \in I^A$ as asserted.
Therefore the algebraic Poisson potential
$\vartheta$ induces 
a morphism 
$\vartheta_I$
of $I^A$-modules
so that the diagram
$$
\CD
D_{I^A}    
@>>>
D_{A}
\\
@V{\vartheta_I}VV
@V{\vartheta}VV
\\
{I^A}
@>>>
{A}
\endCD
$$
is commutative,
and $\vartheta_I$ is
an algebraic Poisson potential for the Poisson algebra $I^A$,
with corresponding Poisson 2-form
$$
\pi_I 
\colon D_{I^A} \otimes_{I^A} D_{I^A}
\longrightarrow
I^A.
$$
\paragraph
By the general theory 
of K\"ahler differentials,
see e.~g. {\smc Kunz} [55],
we know that there is an exact sequence
$$
 I/ I^2
\longrightarrow
 A_{\roman{red}} \otimes _{I^A}  D_{I^A}
\longrightarrow
 D_{\roman{red}}
\longrightarrow
0
$$
of $ A_{\roman{red}}$-modules
where the first arrow is obtained
by sending a class $f \mod  I^2$
to the class of its differential $df$
and where the second one is the obvious morphism of
$ A_{\roman{red}}$-modules.
We now assert that the morphism
$
 A_{\roman{red}} \otimes _{I^A}  D_{I^A}
\longrightarrow  A_{\roman{red}} 
$
induced by
$\vartheta_I$
vanishes on
the image of
$ I/ I^2$
and hence induces an algebraic Poisson potential
$$
 \vartheta
\colon  D_{\roman{red}}
\longrightarrow
 A_{\roman{red}}
$$
for $\{\cdot,\cdot\}_{\roman{red}}$.
To see this, let $f\in  I$, so that
$f = hJ$ for some $h \in A$.
Then the image of
the class
$f \mod  I^2$
in
$ A_{\roman{red}} \otimes _{I^A}  D_{I^A}$
is $[h] \otimes dJ$ where
$[h] \in  A_{\roman{red}} = I^A/ I$
denotes the class of 
$h \in I^A$.
However,  cf. what was said above,
$
\vartheta_I(dJ) = 2J = 0 \in  A_{\roman{red}},
$
and hence
$\vartheta$
induces an algebraic Poisson potential
$$
  \vartheta_{\roman{red}}
\colon  D_{\roman{red}}
\longrightarrow
 A_{\roman{red}}
$$
as asserted. \qed
\enddemo
\paragraph
Since 
$[ \pi_{\roman{red}}] = 0 \in 
\Ho^2_{\roman{Poisson}}
(A_{\roman{red}},{\{\cdot,\cdot\}}_{\roman{red}};
A_{\roman{red}})$,
the Poisson algebra
$(A_{\roman{red}},\{\cdot,\cdot\}_{\roman{red}})$
is representable:
Let 
$M= A_{\roman{red}}\langle b\rangle$
be the free $A_{\roman{red}}$-module with a single basis
element $b$, and define
a $D_{\roman{red}}$-connection
$D$ on $M$ by
$$
D_{\alpha}(b) = -( \vartheta_{\roman{red}} (\alpha))\,b,\quad
\alpha \in  D_{\roman{red}}.
$$
Then Theorem 4.3 yields an $\bold R$-linear representation of
$(A_{\roman{red}},\{\cdot,\cdot\})$,
viewed as a real Lie algebra,
on $M$. 
\paragraph
Likewise,
the Poisson algebra
$(A_{\roman{red}},\{\cdot,\cdot\}_{\roman{red}})$
is quantizable.
In fact, 
let
\linebreak
$ A_{\roman{red}}^{\bold C} =   A_{\roman{red}} \otimes {\bold C}$,
$ D_{\roman{red}}^{\bold C} =   D_{\roman{red}} \otimes {\bold C}$,
etc.,
and define a complex
$ D_{\roman{red}}^{\bold C}$-connection
$\nabla$ on $M^{\bold C}=M\otimes \bold C$
by
$$
\nabla_{\alpha}(b) = -i ( \vartheta_{\roman{red}}(\alpha)) \,b,
\quad \alpha \in  D_{\roman{red}}^{\bold C}.
\tag5.5
$$
Its curvature $K_{\nabla}$ is manifestly given by
$$
K_{\nabla} = -i\, \pi^{\bold C}_{\roman{red}},
\tag5.6
$$
since
$d\,(i \vartheta_{\roman{red}}) = i\, \pi_{\roman{red}}$.
In view of (4.5),
the corresponding representation of the $ A_{\roman{red}}$
underlying real Lie algebra on 
${
M^{\bold C} =  A_{\roman{red}}\otimes \bold C \langle b \rangle
}$ 
is given 
by the formula
$$
\hat a \, s = -i\,\nabla_{da} s + as,
\quad a \in  A_{\roman{red}},\, s \in 
M^{\bold C}.
\tag5.7
$$
\paragraph
Our prequantization construction for
$(A_{\roman{red}},\{\cdot,\cdot\}_{\roman{red}}$
is now complete. To construct the quantum state space
we must introduce a
\lq\lq polarization\rq\rq\  and a Hilbert space structure.
Because of the 
singularity in the classical picture
the situation is 
rather subtle and we
shall give all requisite details.
\paragraph
First we introduce an analogue of the \lq\lq horizontal polarization\rq\rq\ 
in $T^*Q$.
Since 
$ D_{\roman{red}}$
does {\it not\/} act faithfully on
$ A_{\roman{red}}$
some care is needed here, and we proceed as follows:
\paragraph
We write
$(A_0,\{\cdot,\cdot\}_0)$,
$D_0$, and
$\pi_0 \colon D_0 \otimes_{A_0} D_0
\to A_0      
$
for the corresponding structures arising
from the reduced manifold
$B_0= C_0/G$.
As a smooth manifold, $B_0$ is 
a disjoint union of two copies of
the product
$(\bold R^3-0) \times \bold R^3$.
We can in fact
identify $B_0$
with the cotangent bundle on
the two component space
$P = P_{+} \cup P_{-}$ where
$P_{+} (\cong\bold R^3-0)$
and $P_{-} (\cong\bold R^3-0)$
are the subspaces of the cone $C$
given by
$p_0 >0$ and
$p_0 <0$, respectively.
To see this we introduce coordinates as follows:
For $1 \leq j \leq 3$, let
$$
\xi_j = 
\frac {x_0}{p_0} p_j + x_j,
$$
where $p_0 = \pm\sqrt{p_1^2 +p_2^2 +p_3^2}$
according as we are on
$P_{+}$ or $P_{-}$.
In these coordinates the Poisson brackets
in $A_0$ assume the form
$$
\{\xi_i,p_j\} = \delta_{i,j}, \quad
\{\xi_i,\xi_j\} = 0, \quad
\{p_i,p_j\} = 0.
$$
This implies at once that,
as a Poisson manifold, $B_0$
is the cotangent bundle $T^*P$ on $P$.
In particular, $B_0$ is a symplectic manifold.
We denote the symplectic structure by $\sigma$.
Let 
$D_0^{\roman{geo}}$ be the corresponding
$(\bold R,A_0)$-Lie algebra                                   
of smooth 1-forms on $B_0$.
Then, cf. (5.3) above,
the formulas
$$
\vartheta_0^{\roman{geo}}(dp_j) = p_j,\ 
\vartheta_0^{\roman{geo}}(d\xi_j) = 0, \quad 1 \leq j \leq 3,
\tag5.3.0
$$
yield a geometric Poisson potential 
$\vartheta_0^{\roman{geo}}
\colon D_{0}^{\roman{geo}} \to A_0$
for the induced geometric Poisson 2-form
on $B_0$.
Let $L_0$ be the $(\bold R, A_0)$-Lie algebra of smooth
vector fields on $B_0$,
and let
$\sigma^{\flat}\colon D_0^{\roman{geo}} \to L_0$
be the isomorphism of
$(\bold R,A_0)$-Lie algebras
induced by the symplectic structure $\sigma$ on $B_0$,
cf. (3.15).
Under $\sigma^{\flat}$
the geometric Poisson potential 
$\vartheta_0^{\roman{geo}}$
corresponds to
the symplectic potential
$\sum_{1 \leq j \leq 3} p_j d\xi_j$
on $B_0$.
We note that, since 
in the present description of $B_0$ the coordinates
$p_j$ are 
\lq positions\rq\  and the $\xi_j$ \lq momenta\rq,
this is {\it not\/} the standard symplectic potential
$\bold {p dq}$ of a cotangent bundle.
\paragraph
To relate the Poisson algebras $A_{\roman{red}}$ and
$A_0$ we
observe that the inclusion $C_0 \to C(0)$ induces
a morphism 
$A_{\roman{red}} \to A_0$
of Poisson algebras and hence,
as explained in (3.8.4), a morphism
$(A_{\roman{red}},D_{\roman{red}}) \to (A_0,D_0)$
of Lie-Rinehart algebras.
Furthermore, the  
composite of 
$\vartheta_0^{\roman{geo}}$ with the obvious surjection
$
D_{0} \to \colon D_{0}^{\roman{geo}} 
$
is an algebraic Poisson potential
$\vartheta_0 \colon D_0 \longrightarrow A_0$.
By construction, 
the diagram
$$
\CD
 D_{\roman{red}}
@>>>
D_0
\\
@V{ \vartheta}VV
@V{\vartheta_0}VV
\\
{ A_{\roman{red}}}
@>>>
{A_0}
\endCD
$$
is commutative.
\paragraph
Let $H_0^{\roman{geo}} \subseteq D_0^{\roman{geo}}$
be the $(\bold R,A_0)$-Lie subalgebra
generated by the 1-forms
\linebreak
 $dp_1,dp_2,dp_3$.
Under the isomorphism
$\sigma^{\flat}\colon D_0^{\roman{geo}} \to L_0$
of
$(\bold R,A_0)$-Lie algebras
$H_0^{\roman{geo}}$
is mapped isomorphically to the
{\it vertical\/} polarization
of 
$B_0$,
viewed as the cotangent bundle on
the space
$P$.
Notice, however, that 
on the cotangent bundle $T^*Q$ of Minkowski space
$Q$ the 1-forms
$dp_0,dp_1,dp_2,dp_3$
 generate an object corresponding to the
{\it horizontal\/} polarization
of $T^*Q$.
Let $H_0 \subseteq D_0$
be the 
pre-image of
$H_0^{\roman{geo}}$
under the obvious surjection
$D_0 \longrightarrow D_0^{\roman{geo}}$
of
$(\bold R,A_0)$-Lie algebras,
and let
$H_{\roman{red}} \subseteq  D_{\roman{red}}$
be the
pre-image of
$H_0$
under the above morphism
$(A_{\roman{red}}, D_{\roman{red}})
 \longrightarrow (A_0,D_0)$
of Lie-Rinehart algebras.
A little thought reveals that
$ H_{\roman{red}}$
is indeed an $(\bold R, A_{\roman{red}})$-Lie algebra.
It is clear that
the formal differentials $d[p_0],d[p_1],d[p_2],d[p_3]$
lie in $ H_{\roman{red}}$
where $[p_i] \in A_{\roman{red}}$ denotes the class
of $p_i \in I^A$.
Our philosophy is that
an object of the kind $ H_{\roman{red}}$
is the {\it correct\/}
generalization of the concept of a
polarization.
\paragraph
The object $ H_{\roman{red}}$
acts on $M^{\bold C}$
in the usual way through the connection (5.5),
and
the $H_{\roman{red}}$-invariants
$(M^{\bold C})^{H_{\roman{red}}} \subseteq M^{\bold C}$  
consist of elements $ab \in M^{\bold C},\ a \in A_{\roman{red}}^{\bold C}$,
satisfying
$$
\nabla_{d[p_j]} (ab) =0,\quad 0 \leq j \leq 3.
$$
To evaluate $\nabla_{d[p_j]} (ab)$,
we pick a representative 
$\phi \in A\otimes C$ of $a \in A_{\roman{red}}^{\bold C}$ and compute
$$
\nabla_{d[p_j]} (ab) = [(dp_j)^{\sharp}\phi - i \vartheta(dp_j)\phi]b
=[\{dp_j,\phi\}] - i p_j\phi]b
=-\left[\frac{\partial \phi}{\partial x_j} + i p_j\phi\right] b.
$$
Hence
$(M^{\bold C})^{H_{\roman{red}}}$ may be identified
with classes of smooth complex functions $\phi$
on $T^*Q$  
having the property
$\{J,\phi\} \in I \otimes \bold C$
and which, for
$0 \leq j \leq 3,$
satisfy 
$$
\frac {\partial \phi}{\partial x_j}+ i\,p_j \phi \in I \otimes \bold C;
\tag5.9.1
$$
in view of {\smc Gotay's} Lemma 5.2
this means precisely that
these $\phi$ satisfy the equations
$$
\frac {\partial \phi}{\partial x_j}+ i\,p_j \phi =0,
\quad \quad 0 \leq j \leq 3,
\tag5.9.2
$$
on  $C(0)$.
These classes of functions $\phi$
restrict to honest \lq\lq wave functions\rq\rq\ on $B_0$
in the polarization $F=\roman{im} (H_0^{\roman{geo}}) \subseteq L_0$.
In fact, for a function $\phi$ on $T^*Q$
constant along the flow lines
(5.2.2)
with $p_0^2 -p_1^2 -p_2^2 -p_3^2 = 0$,
let  $\phi_0$ be the function on $B_0$
given by
$$
\phi_0(\xi_1,\xi_2,\xi_3,p_1,p_2,p_3) =
\phi(0,\xi_1,\xi_2,\xi_3,p_0,p_1,p_2,p_3),
$$
where 
$p_0 = \pm\sqrt{p_1^2 +p_2^2 +p_3^2}$
according as
$(p_1,p_2,p_3)\in P_{+}$ or $(p_1,p_2,p_3) \in P_{-}$.
If $\phi$ satisfies (5.9.1),
$\phi_0$ satisfies the equations
$$
\frac {\partial \phi_0}{\partial \xi_j}+ i\,p_j \phi_0 =0, 
\quad 1 \leq j \leq 3.
\tag5.9.3
$$
\paragraph
To see that
$(M^{\bold C})^{H_{\roman{red}}} \subseteq M^{\bold C}$ 
is non-empty
we recall that each \lq wave function\rq\ 
on $T^*Q$
in the horizontal polarization
is a solution $\phi$ of the equations
$$
\frac {\partial \phi}{\partial x_j}+ i\,p_j \phi =0,
\quad 0 \leq j \leq 3.
$$
Hence such a $\phi$ satisfies
$$
\{J,\phi\} = 2i\, J\,\phi \in I\otimes\bold C
$$
and therefore represents an element of $A_{\roman{red}}^{\bold C}$,
in fact of
$(M^{\bold C})^{H_{\roman{red}}} \subseteq M^{\bold C}$,
whence the latter
is non-empty.
\paragraph
The next task is to construct the 
requisite Hilbert space.
To this end
we play the same game 
with $B_0$,\ $\pi_0$,\ $\vartheta_0$, etc. 
Virtually the same formula as in (5.5) above then yields a
$D^{\bold C}_0$-connection
$\nabla_0$ on $M_0 = (A_0 \otimes {\bold C}) \langle b \rangle$
having curvature
$-i \pi_0$.
Now
$B_0$ is the cotangent bundle on the space
$P = P_{+} \cup P_{-}$,
and $M_0$
is the space of smooth sections of the trivial complex line bundle
$\Lambda = B_0 \times \bold C$ over $B_0$, with connection
$\nabla_0$.
This situation can be handled by standard geometric quantization
theory, cf.
{\smc \'Sniatycki} [88] and {\smc Woodhouse} [108]:
Pick a metalinear structure on the (real) polarization
$F=\roman{im} (H_0^{\roman{geo}}) \subseteq L_0$.
By the general theory, a metalinear structure
(on the vertical polarization of a cotangent bundle $T^*P$)
exists
if and only if $0 = w_1^2(P) \in \Ho^2(P,\MZ/2)$, where
$w_1(P) \in \Ho^1(P,\MZ/2)$ is the first Stiefel-Whitney class;
the set of all
metalinear structures, if non-empty, is parametrized up to equivalence
by $\Ho^1(P,\MZ/2)$.
Since in our case $P$ is homotopy equivalent to
a disjoint union of two 2-spheres,
a metalinear structure exists
and is unique up to equivalence.
Let $\sqrt {\Lambda^3 F}$ be 
the corresponding half form bundle
over $P$; it carries a canonical
partial flat connection covering $F$.
Under the present circumstances 
the bundle $\sqrt {\Lambda^3 F}$
 is trivial.
Let $\nu$ be a section thereof, suitably normalized
(with respect to its values on appropriately chosen metaframes).
The corresponding space of smooth wave functions
consists of the
polarized sections
$\psi =\phi_0 b \otimes \nu$
 of $\Lambda \otimes \sqrt {\Lambda^3 F}$
where $\phi_0$ is a function on
$B_0$.
We note that \lq polarized\rq \ amounts to
the requirement  that
$\phi_0$ satisfies the equation
(5.9.3).
Moreover,
for a smooth complex function $\alpha$ on $P$,
let
$\phi^{\alpha}$ be the function on $B_0$
given by
$$
\phi^{\alpha}(\xi,p) =
\alpha(p)\roman e^{i\langle \xi,p\rangle}.
$$
Then
the association
$
\alpha \mapsto \phi^{\alpha} b \otimes \nu
$
defines an isomorphism
between
the space of smooth complex functions on
$P$
and
the 
space of smooth wave functions.
Moreover,
perhaps up to constants, the inner product
between two such wave functions
$\psi$ and $\psi'$
is given by
$$
\langle \psi,\psi' \rangle = \int _P \alpha \overline{\alpha'} \varepsilon,
$$
where 
$\varepsilon = dp_1dp_2dp_3$
refers to the natural volume element on $P$.
We write $Y_0$ for the Hilbert space
arising from this construction.
Now given two
functions
$\phi$ and $\phi'$
on $T^*Q$
representing elements of
$(M^{\bold C})^{H_{\roman{red}}} \subseteq M^{\bold C}$,
we define their inner product by
$$
\langle [\phi],[\phi'] \rangle 
= \langle \phi_0 b \otimes \nu,\phi'_0 b \otimes \nu \rangle,
$$
where $\phi_0$ and $\phi'_0$ are the restrictions to $B_0$
of the classes
$[\phi]$ and  $[\phi']$ in $A_{\roman{red}}$.
This yields a pre-Hilbert space;
its completion $Y$
will be our 
 quantum state space
for
$(A_{\roman{red}},\{\cdot,\cdot\}_{\roman{red}})$.
We note that $Y$ may be viewed as a sub Hilbert space of $Y_0$.
\paragraph
Our final task is to construct the
quantum operators
on $Y$ corresponding to the classical observables
in $A_{\roman{red}}$.
As in standard geometric quantization theory,
the polarization restricts the observables
in $A_{\roman{red}}$ that eventually become \lq\lq quantized\rq\rq,
and the formula (5.7) describing the quantum operators must
be modified according to the Lie derivative of
the chosen half forms.
Now in general,
if $g$ is just a Lie algebra,
$h \subseteq g$ a sub Lie algebra, 
and if
$U$ is a $g$-module,
with structure map
$\phi \colon g \to \roman{End}(U)$,
the action of $g$ on $U$ passes to an action 
of the 
pre-image $\phi^{-1}(N_{\overline g}(\overline h)) \subseteq g$ of the
 normalizer
$N_{\overline g}(\overline h)$
of $ \overline h= \phi(h)$ in $\overline g= \phi(g)$
on the invariants $U^h$.
In our situation, 
we recall the corresponding extension (4.1) 
of $(\bold R, A_{\roman{red}})$-Lie algebras
which now looks like
$$
0 
@>>> 
 A_{\roman{red}}
@>>> 
\overline L_{\{\cdot,\cdot\}_{\roman{red}}}
@>>> 
 D_{\roman{red}}
@>>> 
0
$$
and take
$g = \overline L_{\{\cdot,\cdot\}_{\roman{red}}}$,
$h =  H_{\roman {red}}$,
viewed as a sub Lie algebra 
of $\overline L_{\{\cdot,\cdot\}_{\roman{red}}}$
through the obvious
$ D_{\roman{red}}$-connection
$$
 D_{\roman{red}}
\longrightarrow
 A_{\roman{red}} \oplus  D_{\roman{red}} =
\overline L_{\{\cdot,\cdot\}_{\roman{red}}}.
$$
Moreover, we take
$U = M^{\bold C}$.
Our construction will directly quantize the pre-image
in $ A_{\roman{red}}$
of the corresponding normalizer 
$N_{\overline g}(\overline h)$
with respect to the injection
$$
\iota_{\roman{red}} \colon
 A_{\roman{red}}
\longrightarrow
\overline L_{\{\cdot,\cdot\}_{\roman{red}}},
$$
cf. (4.2) above.
This pre-image is
the
sub Lie algebra of
$ A_{\roman{red}}$
consisting of those
$a=[f] \in  A_{\roman{red}}$
which satisfy
the conditions
$$
[df,dp_k] \in H_0^{\roman{geo}},
\quad 0 \leq k \leq 3,
 \tag5.10
$$
where $f$ is a function representing $a \in  A_{\roman{red}}$,
and where
$df$ and $dp_k$ refer to smooth 1-forms on
$B_0$.
To understand what this really means,
for $0 \leq k \leq 3$, we compute
$$
\align
[df,dp_k] &= d\{f,p_k\}
=
d\left(\frac {\partial f}{\partial x_k}\right)
\\
&=
\frac {\partial^2 f}{\partial x_0 \partial x_k} dx_0
+
\dots
+
\frac {\partial^2 f}{\partial x_3 \partial x_k} dx_3
+
\frac {\partial^2 f}{\partial p_0 \partial x_k} dp_0
+
\dots
+
\frac {\partial^2 f}{\partial p_3 \partial x_k} dp_3 \in D_0^{\roman{geo}},
\endalign
$$
where 
$d$\ still refers to the usual exterior derivative of
smooth functions.
We conclude from this that
(5.10) is equivalent to the requirement that,
for every 
$a \in  A_{\roman{red}}$,
each representative $f \in A $ 
of $a$ 
has the property that,
for
$0 \leq j,k\leq 3$, the second partial derivatives
$$
\frac {\partial^2 f}{\partial x_j \partial x_k} 
\tag5.11
$$
lie in the ideal $I$, i.~e.,
in view of {\smc Gotay's} Lemma 5.2,
vanish on the constraint $C(0)$.
Hence 
(5.7)
induces a representation of
the sub Lie algebra 
$ A_{\roman{red}}^{ H_{\roman{red}}}$
of
$( A_{\roman{red}},\{\cdot,\cdot\}_{\roman{red}})$
consisting of classes of smooth functions $f$ on
$T^*Q$ that 
(i) 
satisfy
$\{J,f\} \in I$
{\it and\/} (ii) 
have the additional property that
the second partial derivatives 
$\frac {\partial^2 f}{\partial x_j \partial x_k} $
lie in $I$.
To complete the construction we recall that
on the 
space of smooth wave functions
over $B_0$,
for a function
$k \in A_0$ whose Hamiltonian vector field
$(dk)^{\sharp}$ preserves the polarization $F$
(i.~e. $[(dk)^{\sharp},F] \subseteq F]$),
according to standard geometric quantization theory,
the corresponding operator $\hat k$ is given by
$$
\hat k s = -i\,(\nabla_0)_{dk} s + ks 
- \frac 12 i\,\roman {div}_{\varepsilon}((dk)^{\sharp})s,
\quad k \in  A_0,\, s \in M_0^{\bold C}.
\tag5.12
$$
Here 
for a vector field $Z$ on $B_0$ preserving the polarization $F$
the expression
$\roman {div}_{\varepsilon}(Z)$
is defined by
$$
\lambda_Z \varepsilon =\roman {div}_{\varepsilon}(Z) \varepsilon,
$$
where as before
$\lambda_Z$ refers to the operation of Lie derivative.
Now a vector field
$$
Z = \alpha_1 \frac {\partial}{\partial \xi_1}
  + \alpha_2 \frac {\partial}{\partial \xi_2}
  + \alpha_3 \frac {\partial}{\partial \xi_3}
  + \beta_1 \frac {\partial}{\partial p_1}
  + \beta_2 \frac {\partial}{\partial p_2}
  + \beta_3 \frac {\partial}{\partial p_3}
$$
preserves the polarization $F$ if and only if
$$
\frac {\partial \beta_1}{\partial \xi_1}=0,
\quad
\frac {\partial \beta_2}{\partial \xi_2}=0,
\quad
\frac {\partial \beta_3}{\partial \xi_3}=0,
$$
and a straightforward calculation yields
$$
\lambda_Z \varepsilon =
\left(
\frac {\partial \beta_1}{\partial p_1}+
\frac {\partial \beta_2}{\partial p_2}+
\frac {\partial \beta_3}{\partial p_3}
\right)\varepsilon .
$$
For
$f \in A$ we now write
$$
\roman {div}_{\varepsilon}(f) =
\left(
\frac {\partial^2 f}{\partial x_1 \partial p_1}+
\frac {\partial^2 f}{\partial x_2 \partial p_2}+
\frac {\partial^2 f}{\partial x_3 \partial p_3}
\right) \in A.
$$
A calculation shows that,
whenever
the second partial derivatives
$\frac {\partial^2 f}{\partial x_j \partial x_k}$
of $f$
lie in $I^2$,
$$
\{J,\roman {div}_{\varepsilon}(f)\} \in I,
$$
whence then $\roman {div}_{\varepsilon}(f) \in I^A$.
For
$a \in A_{\roman{red}}$ having the property that
the second partial derivatives 
$\frac {\partial^2 f}{\partial x_j \partial x_k}$ of each representative
$f$ lies in $I^2$
and for
$s \in M^{\bold C}$ we now define
$$
\hat a s = -i\,\nabla_{da} s + as 
- \frac 12 i \,[\roman {div}_{\varepsilon}(f)] s.
\tag5.13
$$
It is manifest that for such an $a \in A_{\roman{red}}$
the restriction $a_0 \in A_0$ to $B_0$
acts on $Y_0$  via (5.12) and hence
the injection $Y \to Y_0$
is compatible with the actions.
Consequently $\hat a$ is a symmetric operator on $Y$;
it is self-adjoint whenever
the Hamiltonian vector field on $B_0$  induced
by $a$ is complete.
For example,
for
$a \in A_{\roman{red}}$ 
represented by a function $f$ of the kind
$$
f(x,p) = u(p) + x_0 v_0(p)
+ x_1 v_1(p)
+ x_2 v_2(p)
+ x_3 v_3(p),
$$
and for a state represented by a \lq wave function\rq\ 
$\phi$ on $T^*Q$ of the kind
$$
\phi(x,p) = \alpha(p) \roman e^{-i\langle x,p\rangle}
$$
where $\alpha$ is a function depending on $p_0,\dots,p_3$ only
as indicated,
this amounts to
$$
\hat a[\phi] = [\hat{\alpha} \roman e^{-i\langle -,-\rangle}],
\tag5.14
$$
where
$$
\hat{\alpha}= u\,\alpha 
-i \,\sum _0^3 v_j \frac{\partial \alpha}{\partial p_j}
-\frac 12 i\, \sum _1^3 \alpha \frac{\partial v_j}{\partial p_j}.
$$
\paragraph
For illustration,
we mention that any function $f$ 
independent of $x_0,\cdots, x_3$
trivially satisfies $\{J,f\} \in I$
and trivially has 
second partial derivatives 
$\frac {\partial^2 f}{\partial x_j \partial x_k}$
in $I^2$,
and 
these requirements are met by
the 
{\it boosts\/}
$\beta_j,\ 1 \leq j \leq 3$, defined by
$$
\align
\beta_j(x_0,x_1,x_2,x_3,p_0,p_1,p_2,p_3)&=
x_0p_j + x_j p_0,\quad 1 \leq j \leq 3,
\\
\intertext{and by the components $\alpha_{k,j},\, 1 \leq k< j \leq 3,$
of angular momentum given by}
\alpha_{k,j}(x_0,x_1,x_2,x_3,p_0,p_1,p_2,p_3)&=
x_kp_j - x_j p_k,\quad 1 \leq k< j \leq 3,
\endalign
$$
too.
The corresponding operators look like
$$
\widehat {[\beta_j]} [\phi] = -i\,[\{\beta_j,\phi\}],\quad
\widehat {[\alpha_{k,j}]} [\phi]= -i\,[\{\alpha_{k,j},\phi\}],
$$
where $\phi \in A$ represents $[\phi] \in A_{\roman{red}}^{\bold C}$.
Notice the Poisson brackets
${
\{\beta_k,\beta_j\} = \alpha_{k,j}.
}$
This shows that indeed 
some observables are quantized by our construction.
In particular, it yields among others a precise description
of those quantum observables on $Y_0$ that arise from
quantum observables of the system
classically described by 
$( A_{\roman{red}},\{\cdot,\cdot\}_{\roman{red}})$.
\paragraph
We leave it to the experts 
to 
decide whether ours is a 
physically meaningful quantization
of a spinless relativistic particle with zero rest mass.
\paragraph
\noindent
{\smc Remark 5.15.}
To get more results a
machinery is needed which enables one to handle
geometric versions of objects 
of the kind $D_{A/I}$ etc.
For example, let $A$ be the algebra of smooth functions
on a smooth finite dimensional manifold $N$, let $K \subseteq N$
be compact, and let $I\subseteq A$ be the ideal 
of functions that vanish on $K$.
When we view e.~g. $A/I$ as the smooth functions 
$C^{\infty}(K)$ in the sense of {\smc Whitney} on $K$,
see e.~g. {\smc Malgrange}~[68],
a candidate for the corresponding geometric object would be
$$
D_{C^{\infty}(K) } /\left(\cap_{x \in K}\mu_x D_{C^{\infty}(K) } \right) ,
$$
where $\mu_x \subseteq C^{\infty}(K) $ refers to the maximal ideal
of functions that vanish on $x \in K$.
An alternate approach would be to 
introduce Poisson varieties and/or schemes
with singularities, thereby
staying entirely within
algebraic geometry.
We intend to address both issues elsewhere.

\bigskip
\centerline{\smc References}
\widestnumber\key{9999}
\medskip
\noindent
\ref \no  1
\by R. Abraham and J. E. Marsden
\book Foundations of Mechanics
\publ Benjamin/Cummings Publishing Company
\yr 1978
\endref
\ref \no  2
\by R. Almeida and P. Molino
\paper Suites d'Atiyah, feuilletages,
et quantification g\'eom\'etrique
\jour S\'em. Geom. Diff. U. Sci. Tech. Languedoc, 1984/85
\yr 1985
\pages 39--59
\endref
\ref \no  3
\by J. M. Arms, J. E. Marsden, and V. Moncrief
\paper  Symmetry and bifurcation of moment mappings
\jour Comm. Math. Phys.
\vol 78
\yr 1981
\pages  455--478
\endref
\ref \no  4
\by V. I. Arnold
\book Mathematical Methods of Classical Mechanics
\bookinfo Graduate Texts in Mathematics, No.~60
\publ Springer
\publaddr Berlin-Heidelberg-New York
\yr 1978
\endref
\ref \no  5
\by M. F. Atiyah
\paper Complex analytic connections in fibre bundles
\jour Trans. Amer. Math. Soc.
\vol 85
\yr 1957
\pages 181--207
\endref
\ref \no  6
\by F. A. Berezin
\paper Some remarks about the associated
envelope of a Lie algebra
\jour Funct. Anal. Appl.
\vol 1
\yr 1967
\pages  91--102
\endref 
\ref \no  7
\by R. Berger
\paper G\'eom\'etrie alg\'ebrique de Poisson
\jour C. R. Acad. Sci. Paris A
\vol 289
\yr 1979
\pages 583--585
\endref
\ref \no  8
\by K. H. Bhaskara and K. Viswanath
\paper Calculus on Poisson manifolds
\jour Bull. London Math. Soc.
\vol 20
\yr 1988
\pages 68--72
\endref
\ref \no  9
\by N. Bourbaki
\book Vari\'et\'es diff\'erentiables et analytiques
\bookinfo Fascicule des resultats
\publ Herman
\publaddr Paris
\yr 1969
\endref
\ref \no 10
\by J. Braconnier
\paper Alg\`ebres de Poisson
\jour C. R. Acad. Sci. Paris S\'erie A
\vol 284
\yr 1977
\pages  1345--1348
\endref
\ref \no 11
\by J. Braconnier
\paper Applications des crochets 
g\'en\'eralis\'es aux alg\`ebres de Lie et de Poisson
\jour C. R. Acad. Sci. Paris S\'erie A
\vol 286
\yr 1978
\pages  877--880
\endref
\ref \no 12
\by J. L. Brylinski
\paper A differential complex for Poisson manifolds
\jour J. of Diff. Geom.
\vol 28
\yr 1988
\pages 93--114
\endref
\ref \no 13
\by H. Cartan
\paper Notions d'alg\`ebre diff\'erentielle;
applications aux groupes de Lie et aux vari\'ete\'es
o\`u op\`ere un groupe
de Lie
\jour Coll. Topologie Alg\'ebrique
\vol 
\paperinfo Brussels
\yr 1950
\pages  15--28
\endref
\ref \no 14
\by H. Cartan and S. Eilenberg
\book Homological Algebra
\publ Princeton University Press
\publaddr Princeton
\yr 1956
\endref
\ref \no 15
\by P. R. Chernoff and J. E. Marsden
 \book Properties of 
infinite dimensional Hamiltonian systems
\bookinfo Lecture Notes in Mathematics, No.~425
\publ Springer 
\publaddr Berlin-Heidelberg-New York
\yr 1974
\endref
\ref \no 16
\by C. Chevalley and S. Eilenberg
\paper Cohomology theory of Lie groups and Lie algebras
\jour  Trans. Amer. Math. Soc.
\vol 63
\yr 1948
\pages 85--124
\endref
\ref \no 17
\by A. Coste, P. Dazord, et A. Weinstein
\paper 
Groupoides symplectiques
\jour Publications du d\'epartement de Math\'ematiques,
Universit\'e Claude Bernard-Lyon
\vol 2A
\yr 1987 
\pages 1--62
\endref
\ref \no 18
\by J. Conn
 \paper Normal forms for analytic Poisson structures
\jour Ann. of Math.
\vol 119
\yr 1984
\pages  577--601
\endref
\ref \no 19
\by J. Conn
 \paper Normal forms for smooth Poisson structures
\jour Ann. of Math.
\vol 121
\yr 1985
\pages  565--593
\endref
\ref \no 20
\by T. J. Courant
\book 
Dirac manifolds
\bookinfo Ph. d. thesis, University of California at Berkeley
\yr 1987
\endref
\ref \no 21
\by M. De Wilde et P. B. A. Le Compte
\paper Crochet de Nijenhuis-Richardson et cohomologie de l'alg\`ebre
de Poisson
\jour C. R. Acad. Sci. Paris I 
\vol 300
\yr 1985
\pages 9--13
\endref
\ref \no 22
\by P. A. M. Dirac
\paper Generalized Hamiltonian systems
\jour Can. J. of Math.
\vol 12
\yr 1950
\pages 129--148
\endref
\ref \no 23
\by P. A. M. Dirac
\book Lectures on Quantum Mechanics
\publ Belfer Graduate School of Science
\publaddr Yeshive University, New York
\yr 1964
\endref
\ref \no 24
\by I. Ya. Dorfman
\paper 
Dirac structures of integrable evolution equations
\jour  Phys. Lett. A
\vol 125
\yr 1987
\pages 240--246
\endref
\ref \no 25
\by M. G. Eastwood, R. Penrose, and R. O. Wells, JR.
\paper  Cohomology and massless fields
\jour Comm. Math. Phys.
\vol 78
\yr 1981
\pages  305--351
\endref
\ref \no 26
\by G. L. Feldman
\paper Global dimensions of rings of differential operators
\jour Trans. Moscow Math. Soc. No. 1
\yr 1982
\pages 123--147
\endref
\ref \no 27
\by I. M. Gelfand and I. Ya. Dorfman
 \paper Hamiltonian operators and algebraic structures related to them
\jour Funct. Anal. Appl.
\vol 13
\yr 1979
\pages  248--262
\endref
\ref \no 28
\by I. M. Gelfand and I. Ya. Dorfman
 \paper Hamiltonian operators and infinite dimensional Lie algebras
\jour Funct. Anal. Appl.
\vol 15
\yr 1981
\pages  173--187
\endref
\ref \no 29
\by I. M. Gelfand and I. Ya. Dorfman
 \paper Hamiltonian operators and the classical Yang-Baxter equations
\jour Funct. Anal. Appl.
\vol 16
\yr 1983
\pages  241--248
\endref
\ref \no 30
\by M. Gerstenhaber
\paper On the deformation of rings and algebras. I.
\jour Ann. of Math.
\vol 79
\yr 1964
\pages  59--103
\moreref
\paper II.
\jour Ann. of Math.
\vol 84
\yr 1966
\pages  1--19
\moreref
\paper III.
\jour Ann. of Math.
\vol 88
\yr 1968
\pages  1--34
\moreref
\paper IV.
\jour Ann. of Math.
\vol 99
\yr 1974
\pages  257--276
\endref
\ref \no 31
\by M. J. Gotay
\paper Poisson reduction and quantization for the $n+1$ photon
\jour J. of Math. Phys.
\vol 25
\yr 1984
\pages 2154--2159
\endref
\ref \no 32
\by M. J. Gotay
\paper Kostant-Souriau quantization of Robertson-Walker cosmologies
with a scalar field
\jour Lecture Notes in Physics
\vol 94
\yr 1979
\pages 293--295
\paperinfo in: Group theoretical methods in Physics, Austin, 1978,
ed. W. Beiglb\"ock, A. B\"ohm, and E. Takasugi
\publ Springer 
\publaddr Berlin-Heidelberg-New York
\endref
\ref \no 33
\by M. J. Gotay and J. \'Sniatycki
\paper Quantization of presymplectic dynamical systems
\jour Comm. Math. Phys.
\vol 82
\yr 1981
\pages 377--389
\endref
\ref \no 34
\by M. J. Gotay, J. A. Isenberg,
J. Marsden, R. Montgomery, J. \'Sniatycki,
and P. B. Yasskin
\book Constraints and moment mappings in relativistic field theory
\publ Math. Sci. Research Institute
\publaddr Berkeley California 94720
\finalinfo To appear, approx. 1989 
\endref
\ref \no 35
\by V. Guillemin and S. Sternberg
\book Symplectic techniques in Physics
\publ Cambridge University Press
\publaddr London/New York
\yr 1984
\endref
\ref \no 36
\by V. W. Guillemin and S. Sternberg
\paper Geometric quantization and multiplicities of group representations
\jour Invent. Math.
\vol 67
\yr 1982
\pages 515--538
\endref
\ref \no 37
\by J. Herz
\paper Pseudo-alg\`ebres de Lie
\jour C. R. Acad. Sci. Paris 
\vol 236
\yr 1953
\pages 1935--1937
\endref
\ref \no 38
\by P. J. Higgins and K. Mackenzie
\paper Algebraic constructions in the category of Lie algebroids
\paperinfo to appear
\jour J. of Algebra
\endref
\ref \no 39
\by P. J. Hilton and U. Stammbach
\book A Course In Homological Algebra
\bookinfo Graduate texts in mathematics, vol. 4
\publ Springer
\publaddr Berlin--Heidelberg--New York
\yr 1971
\endref
\ref \no 40
\by G. Hochschild, B. Kostant, and A. Rosenberg
\paper Differential forms on regular affine algebras
\jour  Trans. Amer. Math. Soc.
\vol 102
\yr 1962
\pages 383--408
\endref
\ref \no 41
\by J. Huebschmann
\paper Extensions of Lie-Rinehart algebras. I. The Chern-Weil construction
\paperinfo Preprint, Universit\"at Heidelberg, 1989
\moreref
\paper II. The spectral sequence and the Weil algebra
\paperinfo Preprint, Universit\"at Heidelberg, 1989
\moreref
\paper III. Principal extensions
\paperinfo in preparation
\endref
\ref \no 42
\by N. Jacobson
\book Lie algebras
\bookinfo Interscience tracts in pure and applied mathematics, vol. 10
\publ Interscience, J. Wiley and Sons
\publaddr New York and London
\yr 1962
\endref
\ref \no 43
\by F. W. Kamber and Ph. Tondeur
\book Invariant differential operators and the cohomology
of Lie algebra sheaves
\bookinfo Memoirs of the Amer. Math. Soc. vol. 113
\yr 1971
\publ Amer. Math. Soc.
\publaddr Providence, R. I.
\endref
\ref \no 44
\by M. V. Karasev
\paper Analogues of objects of Lie group theory for non-linear
Poisson brackets
\paperinfo (Russian)
\jour Izvestiya
\vol 50, 1-3
\yr 1986
\moreref
\jour engl. translation: Math. of the USSR Izvestiya
\vol 28
\yr 1987
\pages 497--527
\endref
\ref \no 45
\by C. Kassel
\paper L'homologie cyclique des alg\'ebres enveloppantes
\jour Invent. Math.
\vol 91
\yr 1988
\pages  221--251
\endref
\ref \no 46
\by A. A. Kirillov
\paper Infinite dimensional groups
\jour Proc. Intern. Congr. Math. Helsinki
\vol 
\yr 1978
\pages 705--708
\endref
\ref \no 47
\by A. A. Kirillov
\book \'El\'ements de la th\'eorie des repr\'esentations
\publ Mir
\publaddr Moscow
\yr 1974
\endref
\ref \no 48
\by B. Kostant
\paper Quantization and unitary representations
\jour Lecture Notes in Mathematics
\vol 170
\yr 1970
\pages 87--207
\paperinfo
Lectures in Modern Analysis and Applications, III, ed. C. T. Taam
\publ Springer
\publaddr Berlin-Heidelberg-New York
\endref
\ref \no 49
\by B. Kostant
\paper Graded manifolds, graded Lie theory, and prequantization
\jour Lecture Notes in Mathematics
\vol 570
\yr 1977
\pages 177--306
\paperinfo in Diff. Geom. Meth. in Math. Physics, Bonn, 1975
\publ Springer 
\publaddr Berlin-Heidelberg-New York
\endref
\ref \no 50
\by B. Kostant
\paper Orbits, symplectic structures, and representation theory
\jour Proc. U. S. Japan seminar in Diff. Geom. Kyoto 1965
\publ Nippon Hyoronisha
\publaddr Tokyo
\yr 1965
\endref
\ref \no 51
\by B. Kostant and S. Sternberg
\paper Symplectic reduction, BRS-cohomology, and infinite dimensional
Clifford algebras
\jour Ann. of Physics (N. Y.)
\vol 176
\yr 1987
\pages 49--113
\endref
\ref \no 52
\by J. L. Koszul
\paper Crochet de Schouten-Nijenhuis et cohomologie
\jour Asterisque,
\vol hors-serie,
\yr 1985
\pages 251--271
\paperinfo in E. Cartan et les Math\'ematiciens d'aujourd hui, 
Lyon, 25--29 Juin, 1984
\endref
\ref \no 53
\by J. L. Koszul
\book Lectures on fibre bundles and differential geometry
\publ Tata Institute of Fundamental research
\publaddr Bombay
\yr 1960
\finalinfo Reprinted: Springer 1986
\endref
\ref \no 54
\by I. S. Krasil'shchik, V. V. Lychagin, and A. M. Vinogradov
\book Geometry of Jet Spaces and Nonlinear Partial Differential Equations
\bookinfo Advanced Studies in Contemporary Mathematics, vol. 1
\publ Gordon and Breach Science Publishers
\publaddr New York, London, Paris, Montreux, Tokyo
\yr 1986
\endref
\ref \no 55
\by E. Kunz
\book K\"ahler differentials
\bookinfo Vieweg Advanced Lectures in Mathematics
\publ Friedrich Vieweg \& Sohn
\publaddr Braunschweig/Wiesbaden
\yr 1986
\endref
\ref \no 56
\by S. Lang
\book Differential manifolds
\publ Springer
\publaddr Berlin-Heidelberg-New York
\yr 1985
\endref
\ref \no 57
\by A. Lichnerowicz
\paper Vari\'et\'e symplectique et dynamique associ\'ee \`a une
sous-vari\'et\'e
\jour C. R. Acad. Sci. Paris A
\vol 280
\yr 1975
\pages 523--527
\endref
\ref \no 58
\by A. Lichnerowicz
\paper Les vari\'et\'es de Poisson et leurs alg\`ebres de Lie
associ\'ees
\jour J. Diff. Geo.
\vol 12
\yr 1977
\pages 253--300
\endref
\ref \no 59
\by A. Lichnerowicz
\paper Vari\'et\'es canoniques et transformations canoniques
\jour C. R. Acad. Sci. Paris A
\vol 280
\yr 1975
\pages 37--40
\endref
\ref \no 60
\by A. Lichnerowicz
\paper Cohomologie 1-diff\'erentiable
et d\'eformations de l'alg\`ebre de Lie dynamique d'une vari\'et\'e
canonique
\jour C. R. Acad. Sci. Paris A
\vol 280
\yr 1975
\pages 1217--1220
\endref
\ref \no 61
\by A. Lichnerowicz
\paper 
L'alg\`ebre de Lie 
des automorphismes infinit\'esimales symplectiques
\jour Symp. Math.
\vol 14
\yr 1974
\pages 11--24
\endref
\ref \no 62
\by A. Lichnerowicz
\paper 
D\'eformations de l'alg\`ebre de Lie de Poisson
d'une vari\'et\'e symplectique
\jour Pub. del V congreso de la Ag. de Matem\'aticos
de l'Expresi\'on Latina, Madrid 1977
\yr 1978
\pages 194--206
\endref
\ref \no 63
\by S. Lie
\book 
Theorie der Transformationsgruppen
\publ Teubner
\publaddr Leipzig
\yr 1890
\endref
\ref \no 64
\by K. Mackenzie
\book Lie groupoides and Lie algebroides in differential geometry
\bookinfo London Math. Soc. Lecture Note Series, vol. 124
\publ Cambridge University Press
\publaddr Cambridge, England
\yr 1987
\endref
\ref \no 65
\by S. Mac Lane
\paper Hamiltonian mechanics and geometry
\jour Amer. Math. Monthly
\vol 77
\yr 1970
\pages 570--586
\endref
\ref \no 66
\by F. Magri and C. Morosi
\paper A geometrical characterization of integrable Hamiltonian systems
through the theory of Poisson-Nijenhuis manifolds
\jour Quaderno
\vol S/19
\publaddr Dipartimento di Mathematica, Universit\`a di Milano
\yr 1984
\endref
\ref \no 67
\by F. Magri, C. Morosi, and O. Ragnisco
\paper 
Reduction techniques for infinite dimensional Hamiltonian systems,
some ideas and applications
\jour Comm. in Math. Physics
\vol 99
\yr 1985
\pages 115--140
\endref
\ref \no 68
\by B. Malgrange
\book Ideals of differentiable functions
\publ Oxford University Press
\publaddr Oxford
\yr 1966
\endref
\ref \no 69
\by M. P. Malliavin
\paper 
Alg\`ebre homologique et op\'erateurs diff\'erentiels
\paperinfo Ring theory,
Proc. of a conference held at Granada, 1986,
 eds. J. L. Buesco P. Jara and B. Torrecillas
\jour Lecture Notes in Mathematics vol. 1328
\publ Springer
\publaddr Berlin-Heidelberg-New York-Tokyo
\yr 1988
\pages  172--186
\endref
\ref \no 70
\by J. Marsden
\book Applications of Global Analysis in Mathematical
Physics
\bookinfo Mathematical Lecture Series,
 No. 2
\publ Publish or Perish, Inc.
\publaddr Berkeley
\yr 1974
\endref
\ref \no 71
\by J. Marsden ed.
\paper
Fluids and Plasmas: Geometry and Dynamics
\paperinfo
 Proc. AMS Conference Boulder
\jour Cont. Math.
\vol 28
\yr 1983
\endref
\ref \no 72
\by J. Marsden and A. Weinstein
\paper Reduction of symplectic manifolds with symmetries
\jour Rep. on Math. Phys.
\vol 5
\yr 1974
\pages 121--130
\endref
\ref \no 73
\by W. S. Massey and F. P. Petersen
\paper The cohomology structure of certain fibre spaces.I
\jour Topology
\vol 4
\yr 1965
\pages  47--65
\endref
\ref \no 74
\by K. Mikami and A. Weinstein
\paper Moments and reduction for symplectic groupoid actions
\jour Publ. RIMS Kyoto University
\vol 24
\yr 1988
\pages 121--140
\endref
\ref \no 75
\by J. W. Milnor and J. C. Moore
\paper On the structure of Hopf algebras
\jour Ann. of Math.
\vol 81
\yr 1965
\pages  211--264
\endref
\ref \no 76
\by A. Nijenhuis
\paper Jacobi-type identities for bilinear differential concomitants
of tensor fields.I
\jour Indag. Math.
\vol 17
\yr 1955
\pages 390--403
\endref
\ref \no 77
\by R. Palais
\paper The cohomology of Lie rings
\jour  Proc. Symp. Pure Math.
\vol III
\yr 1961
\pages 130--137
\paperinfo Amer. Math. Soc., Providence, R. I.
\endref
\ref \no 78
\by J. Pradines
\paper Th\'eorie de Lie pour les groupoides diff\'erentiables.
Calcul diff\'erentiel dans la cat\'egorie de groupoides infinit\'esimaux
\jour  C. R. Acad. Sci. Paris 
\vol 264
\yr 1967
\pages 245--248
\endref
\ref \no 79
\by J. Rawnsley, W. Schmid, and J. A. Wolf
\paper Singular unitary representations and indefinite harmonic theory
\jour J. Funct. Anal.
\vol 51
\yr 1983
\pages 1--114
\endref
\ref \no 80
\by G. Rinehart
\paper Differential forms for general commutative algebras
\jour  Trans. Amer. Math. Soc.
\vol 108
\yr 1963
\pages 195--222
\endref
\ref \no 81
\by J. A. Schouten
\paper \"Uber Differentialkomitanten zweier kontravarianter Gr\"o{\ss}en
\jour Proc. Kon. Ned. Akad. Wet. Amsterdam
\vol 43
\yr 1940
\pages 449--452
\endref
\ref \no 82
\by J. A. Schouten
\paper On the differential operators of first order in tensor calculus
\jour Convegno Int. Geom. Diff. Italia
\yr 1953
\pages 1--7
\publ ed. Cremonese
\publaddr Roma, 1954
\endref
\ref \no 83
\by J. T. Schwartz
\book Nonlinear Functional Analysis
\publ Gordon and Breach
\yr 1967
\endref
\ref \no 84
\by I. E. Segal
\paper Quantization of non-linear systems
\jour J. of Math. Phys.
\vol 1
\yr 1960
\pages 468--488
\endref
\ref \no 85
\by I. E. Segal
\paper Symplectic structures and the quantization problem for wave functions
\jour Symp. Math. 
\vol 14
\yr 1974
\pages 89--118
\publ Academic Press
\publaddr London
\endref
\ref \no 86
\by D. Simms
\paper On the Schr\"odinger equation given by geometric quantization
\jour Lecture Notes in Mathematics
\vol 570 
\yr 1977
\pages 351--356
\paperinfo in Diff. Geom. Meth. in Math. Physics, Bonn, 1975
\publ Springer 
\publaddr Berlin-Heidelberg-New York
\endref
\ref \no 87
\by D. J. Simms and N. M. J. Woodhouse
\book Lectures on Geometric quantization
\bookinfo Lecture Notes in Physics, No. 53
\publ Springer
\publaddr Berlin-Heidelberg-New York
\yr 1976
\endref
\ref \no 88
\by J. \'Sniatycki
\book Geometric quantization and quantum mechanics
\bookinfo Applied Mathematical Sciences
 No.~30
\publ Springer
\publaddr Berlin-Heidelberg-New York
\yr 1980
\endref
\ref \no 89
\by J. \'Sniatycki 
\paper Constraints and quantization
\paperinfo in: Nonlinear partial differential operators
and quantization procedures,
Clausthal 1981
\jour Lecture Notes in Mathematics, No.~1037
\pages 301--334
\publ Springer
\publaddr Berlin-Heidelberg-New York
\yr 1983
\endref
\ref \no 90
\by J. \'Sniatycki 
\paper On quantization of systems with constraints 
\paperinfo in: Diff. geom. meth. in math. physics,
Jerusalem 1982
\jour Math. Phys. Studies, No.~6
\pages 207--211
\publ Reidel
\publaddr Dordrecht-Boston, Mass
\yr 1984
\endref
\ref \no 91
\by J. \'Sniatycki
\paper Applications of geometric quantization in quantum mechanics
\jour Lecture Notes in Mathematics
\vol 570
\yr 1977
\pages 357--367
\paperinfo in Diff. Geom. Meth. in Math. Physics, Bonn, 1975
\publ Springer 
\publaddr Berlin-Heidelberg-New York
\endref
\ref \no 92
\by J. \'Sniatycki and A. Weinstein
\paper Reduction and quantization for singular moment mappings
\jour Lett. Math. Phys.
\vol 7
\yr 1983
\pages 155--161
\endref
\ref \no 93
\by J. M. Souriau
\paper 
Quantification g\'eom\'etrique
\jour Comm. Math. Physics
\vol 1
\yr 1966
\pages 374--398
\endref
\ref \no 94
\by J. M. Souriau
\book 
Structure des Syst\`emes dynamiques
\publ Dunod
\publaddr Paris
\yr 1970
\endref
\ref \no 95
\by J. D. Stasheff
\paper Constrained Hamiltonians
\book Conference on Elliptic Cohomology
\bookinfo Institute for Advanced Study,
Princeton, September 15 --17, 1986,
ed. P. S. Landweber,
 Lecture Notes in Mathematics, No. 1326
\publ Springer
\publaddr Berlin-Heidelberg-New York
\yr 1988
\pages 150--160
\endref
\ref \no 96
\by J. D. Stasheff
\paper Constrained Poisson algebras
\jour Bull. Amer. Math. Soc.
\vol 19
\yr 1988
\pages 287--290
\endref
\ref \no 97
\by W. M. Tulczyjew
\paper Poisson brackets and canonical manifolds
\jour Bull. Acad. Polon. Sci. Ser. Sci. Math. Astronom. Phys.
\vol 22
\yr 1974
\pages 931--935
\endref
\ref \no 98
\by W. M. Tulczyjew
\paper The graded Lie algebra of multi-vector fields and the
generalized Lie derivative of forms
\jour Bull. Acad. Polon. Sci. Ser. Sci. Math. Astronom. Phys.
\vol 22
\yr 1974
\pages 937--942
\endref
\ref \no 99
\by G. M. Tuynman and W. A. J. J. Wiegerinck
\paper Central extensions and physics
\jour J. of Geom. and Physics
\vol 4
\yr 1987
\pages 207--258
\endref
\ref \no 100
\by A. M. Vinogradov and I. S. Krasil'shchik
\paper What is the Hamiltonian formalism
\paperinfo (Russian)
\jour Uspekhi Mat. Nauk.
\vol 30:1
\pages 173--196
\yr 1975
\moreref
\jour engl. translation:
Russian Math. Surveys
\vol 30:1
\pages 177--202
\yr 1975
\endref
\ref \no 101
\by A. Weil
\book Vari\'et\'es K\"ahleriennes
\publ Herman
\publaddr Paris
\yr 1958
\endref
\ref \no 102
\by A. Weinstein
\paper The local structure of Poisson manifolds
\jour J. of Diff. Geom.
\vol 18
\yr 1983
\pages 523--557
\endref
\ref \no 103
\by A. Weinstein
\paper Poisson structures
\jour Asterisque,
\vol hors-serie, 
\yr 1985
\pages 421--434
\paperinfo in E. Cartan et les Math\'ematiciens d'aujourd hui, 
Lyon, 25--29 Juin, 1984
\endref
\ref \no 104
\by A. Weinstein
\paper Symplectic groupoides and Poisson manifolds
\jour Bull. Amer. Math. Soc.
\vol 16
\yr 1987
\pages 101--104
\endref
\ref \no 105
\by A. Weinstein
\paper Coisotropic calculus and Poisson groupoides
\jour J. Math. Soc. Japan
\vol 40
\yr 1988
\pages 705--727
\endref
\ref \no 106
\by A. Weinstein
\paper Some remarks on dressing transformations
\jour J. Fac. Sci. Univ. Tokyo Sect. 1A, Math.
\vol 36
\yr 1988
\pages 163--167
\endref
\ref \no 107 
\by E. Witten
\paper Topological quantum field theory
\jour Comm. Math. Phys.
\vol 117
\yr 1988
\pages 353--386
\endref
\ref \no 108
\by N. Woodhouse
\book Geometric quantization
\publ Clarendon Press
\publaddr Oxford
\yr 1980
\endref
\enddocument